\documentclass[onecolumn,11pt,letterpaper]{article}
\renewcommand{\include}{\input}

\usepackage{amsmath,amssymb,mathrsfs,theorem}
\usepackage{graphicx}
\usepackage{epsfig,color,psfrag}
\usepackage{xspace}
\usepackage[sort,numbers]{natbib}
\usepackage{vmargin}
\setmarginsrb{2cm}{1.5cm}{2cm}{4cm}%
             {10pt}{10pt}{15pt}{21pt}
\newtheorem{theorem}{Theorem}[section]
\newtheorem{proposition}[theorem]{Proposition}

\newtheorem{definition}[theorem]{Definition}
\newtheorem{lemma}[theorem]{Lemma}
{\theorembodyfont{\rmfamily} % These ones have roman font in their body
\newtheorem{remark}[theorem]{Remark}
\newtheorem{remarks}[theorem]{Remarks}
\newtheorem{example}[theorem]{Example}}

\renewcommand{\mod}{\mbox{mod\,}}

\newcommand{\dist}{\operatorname{dist}}
\newcommand{\distC}{\operatorname{\subscr{dist}{\texttt{c}}}}
\newcommand{\distCC}{\operatorname{\subscr{dist}{\texttt{cc}}}}
\newcommand{\real}{\ensuremath{\mathbb{R}}}
\newcommand{\complex}{\ensuremath{\mathbb{C}}}

\newcommand{\sph}{{\mathbb{S}}}

\newcommand{\Pc}{{\mathcal{P}}}
\renewcommand{\natural}{{\mathbb{N}}}
\newcommand{\naturalzero}{{\mathbb{N}_0}}
\newcommand{\map}[3]{#1\colon#2\rightarrow#3}
\newcommand{\subscr}[2]{{#1}_{\textup{#2}}}

\newcommand{\MM}{\mathcal{M}}
\newcommand{\nll}{\textup{\texttt{null}}\xspace}
\newcommand{\ctrl}{\textup{ctl}}
\newcommand{\msg}{\textup{msg}}

\newcommand{\stf}{\textup{stf}}

\newcommand{\union}{\cup}

\newcommand{\fl}[1]{\left\lfloor #1\right\rfloor}
\newcommand{\eps}{\varepsilon}
\newcommand{\e}{\mathbf{e}} %% FB: who decided this? This is bad latex
\renewcommand{\bar}{\overline}

\newcommand{\reschedule}[3]{{#1}_{(#2,#3)}}

\newcommand{\setdef}[2]{\{#1 \; | \; #2\}}
\newcommand{\bigsetdef}[2]{\big\{#1 \; | \; #2\big\}}
\newcommand{\Bigsetdef}[2]{\Big\{#1 \; | \; #2\Big\}}

\newcommand{\average}{\operatorname{avrg}}
\newcommand{\piSpace}[1]{\pi_{#1}}

\newcommand{\PP}{\mathcal{P}}
\newcommand{\pder}[2]{\frac{\partial #1}{\partial #2}}

\newcommand{\argmax}{\operatorname{argmax}}
\newcommand{\circmat}{\operatorname{Circ}}
\newcommand{\tridmat}{\operatorname{Trid}}
\newcommand{\atridmat}{\operatorname{ATrid}}
\newcommand{\yave}{\subscr{y}{ave}}
\newcommand{\xave}{\subscr{x}{ave}}
\newcommand{\zave}{\subscr{z}{ave}}

\newcommand{\Hset}{H}%
\newcommand{\nHset}{n_H}%
\newcommand{\fcase}[1]{($\mathfrak{#1}$)\xspace}

\newcommand{\timeschedule}{\mathbb{T}}
\newcommand{\agent}{A}
\newcommand{\setofagents}{\mathcal{A}}
\newcommand{\true}{\textup{\texttt{true}}\xspace}
\newcommand{\false}{\textup{\texttt{false}}\xspace}

\newcommand{\diam}{\operatorname{diam}}

\newcommand{\Dc}{\mathcal{D}}

\newcommand{\until}[1]{\{1,\dots,#1\}}
\newcommand{\HH}{{\mathcal{H}}}

\newcommand{\WW}{{\mathcal{W}}}
\newcommand{\XX}{{\mathcal{X}}}

\newcommand{\vers}{\operatorname{vers}}
\newcommand{\kprop}{k_{\text{prop}}}
\newcommand{\diag}{\operatorname{diag}}
\newcommand{\intersection}{\ensuremath{\operatorname{\cap}}}
\newcommand{\intersect}{\ensuremath{\operatorname{\cap}}}
\newcommand{\argmin}{\ensuremath{\operatorname{argmin}}}
\newcommand{\oball}[2]{B_{#1}(#2)}
\newcommand{\cball}[2]{\overline{B}_{#1}(#2)}
\renewcommand{\oball}[2]{B(#2,#1)}
\renewcommand{\cball}[2]{\overline{B}(#2,#1)}

\newcommand{\Bigcball}[2]{\overline{B}\Big(#2,#1\Big)}

\newcommand{\myboolean}{\textup{\texttt{BooleSet}}\xspace}

\newcommand{\dm}{\texttt{direction}\xspace}
\newcommand{\prior}{\texttt{priority}\xspace}
\renewcommand{\dm}{\texttt{drctn}\xspace}
\renewcommand{\prior}{\texttt{prior}\xspace}
\newcommand{\csym}{\texttt{c}\xspace}
\newcommand{\ccsym}{\texttt{cc}\xspace}

\newcommand{\norm}[2]{\|#1\|_{#2}}
\newcommand{\bignorm}[2]{\big\|#1\big\|_{#2}}
\newcommand{\Bignorm}[2]{\Big\|#1\Big\|_{#2}}
\newcommand{\Enorm}[1]{\|#1\|_{2}}
\newcommand{\Infnorm}[1]{\|#1\|_{\infty}}

\newcommand{\maximize}{\text{maximize}}
\newcommand{\subj}{\text{subject to}}

\newcommand{\network}[1][]{\Sigma_{\textup{#1}}}
\renewcommand{\network}[1][]{\mathcal{S}_{\textup{#1}}}
\newcommand{\supind}[2]{{#1}^{[#2]}}
\newcommand{\FCC}[1][]{\mathcal{CC}_{\textup{#1}}}
\newcommand{\task}[1][]{\mathcal{T}_{\textup{#1}}}
\newcommand{\msgstandard}{\textup{msg}_{\textup{std}}}
\newcommand{\sampled}[1]{\subscr{#1}{smpld}}
\newcommand{\realdisk}{$\real^d$,$r$-disk}
\newcommand{\realLD}{$\real^d$,$r$-LD}
\newcommand{\realdiskone}{$\real$,$r$-disk}
\newcommand{\realLDone}{$\real$,$r$-LD}
\newcommand{\realIdisk}{$\real^d$,$r$-$\infty$-disk}
\newcommand{\circledisk}{$\sph^1$,$r$-disk}

\newcommand{\rdisk}{$r$-\textup{disk}}
\newcommand{\rLD}{$r$-\textup{LD}}
\newcommand{\rIdisk}{$r$-$\infty$-\textup{disk}}

\newcommand{\CircumC}{\operatorname{Circum}}
\newcommand{\Centroid}{\operatorname{Centroid}}
\newcommand{\Mass}{\operatorname{Mass}}

\newcommand{\agreepursuitName}{{agree-and-pursue}\xspace}
\newcommand{\agreepursuit}{agr-pursuit}
\newcommand{\VicsekName}{{move-toward-average}\xspace}
\newcommand{\VicsekSub}{{avrg}\xspace}

\newcommand{\circumcenter}{crcmcntr}
\newcommand{\parallelcircumcenter}{pll-crcmcntr}
\newcommand{\centroid}{centrd}
\newcommand{\rendezvous}{rndzvs}
\newcommand{\equidistance}{eqdstnc}
\newcommand{\deployment}{deplmnt}
\newcommand{\sat}[1]{\operatorname{sat}_{#1}}
\newcommand{\TC}{\operatorname{TC}}
\newcommand{\MCC}{\operatorname{MCC}}
\newcommand{\TCC}{\operatorname{TCC}}
\newcommand{\cc}{\operatorname{cc}}
\newcommand{\nonnllmsgs}{{cmm$\setminus\emptyset$}}
\newcommand{\CommCost}[1][]{\subscr{\operatorname{C}}{rnd}^{#1}}
\newcommand{\newC}{\CommCost[L]\circ\subscr{E}{\nonnllmsgs}}%
\newcommand\oprocendsymbol{\hbox{$\bullet$}}
\newcommand\oprocend{\relax\ifmmode\else\unskip\hfill\fi\oprocendsymbol}
\def\eqoprocend{\tag*{$\bullet$}}

\newcommand{\margin}[1]{\marginpar{\tiny #1}}
\newcommand{\todo}[1]{\par\noindent{\raggedright\textsc{#1}\par\marginpar{\Large $\star$}}}

\renewcommand{\margin}[1]{}
\renewcommand{\todo}[1]{}

\begin{document}
\title{Synchronous robotic networks and complexity of
    control and communication laws}

\date{\today}

\author{Sonia Mart{\'\i}nez$^{*}$ \quad Francesco Bullo\thanks{Sonia
    Mart{\'\i}nez and Francesco Bullo are with the Department of Mechanical
    and Environmental Engineering, University of California at Santa
    Barbara, Santa Barbara, California 93106,
    \texttt{\{smartine,bullo\}@engineering.ucsb.edu}} \quad%
  Jorge Cort\'es\thanks{Jorge Cort\'es is with the Department of Applied
    Mathematics and Statistics, University of California at Santa Cruz,
    Santa Cruz, California 95064, \texttt{jcortes@ucsc.edu}} \quad%
  Emilio Frazzoli\thanks{Emilio Frazzoli is with the Department of
    Mechanical and Aerospace Engineering, University of California at Los
    Angeles, Los Angeles, California 90095, \texttt{frazzoli@ucla.edu}} }

\maketitle

\begin{abstract}
  This paper proposes a formal model for a network of robotic agents that
  move and communicate.  Building on concepts from distributed computation,
  robotics and control theory, we define notions of robotic network,
  control and communication law, coordination task, and time and
  communication complexity. We then analyze a number of basic motion
  coordination algorithms such as pursuit, rendezvous and deployment.
\end{abstract}

\section{Introduction}\label{se:introduction}
\paragraph*{Problem motivation}
The study of networked mobile systems presents new challenges that lie at
the confluence of communication, computing, and control.  In this paper we
consider the problem of designing joint communication protocols and control
algorithms for groups of agents with controlled mobility.  For such groups
of agents we define the notion of communication and control law by
extending the classic notion of distributed algorithm in synchronous
networks.  Decentralized control strategies are appealing for networks of
robots because they can be scalable and they provide robustness to vehicle
and communication failures.

One of our key objectives is to develop a computable theory of time and
communication complexity for motion coordination algorithms.  Hopefully,
our formal model will be suitable to analyze objectively the performance of
various coordination algorithms.  It is our contention that such a theory
is required to assess the complex trade-offs between computation,
communication and motion control or, in other words, to establish what
algorithms are \emph{scalable} and practically implementable in large
networks of mobile autonomous agents.  The need for modern models of
computation in wireless and sensor network applications is discussed in the
well-known report \cite{CNSEC:01}.

\paragraph*{Literature review}

To study complexity of motion coordination, our starting points are the
standard notions of \emph{synchronous and asynchronous networks} in
distributed and parallel computation, e.g., see Lynch~\cite{NAL:97} and,
with an emphasis on numerical methods, Bertsekas and
Tsitsiklis~\cite{DPB-JNT:97}.  This established body of knowledge, however,
is not applicable to the robotic network setting because of the agents'
mobility and the ensuing dynamic communication topology.

An important contribution towards a network model of mobile interacting
robots is introduced by Suzuki and Yamashita~\cite{IS-MY:99}, see also
\cite{HA-YO-IS-MY:99,KS-IS:96}.  The Suzuki-Yamashita model consists of a
group of ``distributed anonymous mobile robots'' that interact by sensing
each other's relative position.  A related model is presented
in~\cite{PF-GP-NS-PW:99,PF-GP-NS-PW:05}. A brief survey of models,
algorithms, and the need for appropriate complexity notions is presented
in~\cite{NS:01}.  This model is discussed in
Section~\ref{se:suzuki/yamashita}.

Recently, a notion of communication complexity for control and
communication algorithms in multi-robot systems is analyzed
in~\cite{EK:02a}, see also~\cite{EK:03b} where a formal model of
communication and control laws for multi-agent networks is proposed.  A
general modeling paradigm is discussed in~\cite{NAL-RS-FV:03}. The time
complexity of a distributed algorithm for coordinated motion planning is
computed in \cite{MS-EF:04}.

Finally, we conclude this literature review session with a general overview
of the recent large research effort dedicated to multi-vehicle systems in
the robotics and control domains.  A survey on cooperative mobile robotics
is presented in~\cite{YUC-ASF-AK:97} and an overview of control and
communication issues is contained in~\cite{EK-RMM:04}. Specific topics of
interest include formation control~\cite{TB-RA:98,JHR-HW:99,ME-XH:01b,
  IS-MY:99,EWJ-PSK:04,XD-AK:02},
rendezvous~\cite{HA-YO-IS-MY:99,PF-GP-NS-PW:05,JL-ASM-BDOA:04b,
  ZL-MB-BF:04a,JC-SM-FB:04h-tmp}, flocking~\cite{AJ-JL-ASM:02}, swarm
aggregation~\cite{VG-KMP:03}, gradient climbing~\cite{PO-EF-NEL:04}, cyclic
pursuit~\cite{JM-MB-BF:04c}, robotic exploration~\cite{IAW-ML-AMB:00},
deployment~\cite{JC-SM-TK-FB:02j,JC-SM-FB:03p-tmp}, and
foraging~\cite{KMP:04,EF-FB:03r}.  Consensus and control theoretical
problems on dynamic graphs are discussed
in~\cite{ROS-RMM:03c,LM:04,WR-RWB:03} and in~\cite{MM:04}, respectively.

\paragraph*{Statement of contributions}
We summarize our approach as follows.  A \emph{robotic network} is a group
of robotic agents moving in space and endowed with communication
capabilities.  The agents position obey a differential equation and the
communication topology is a function of the agents' relative positions.
Each agent repeatedly performs communication, computation and physical
motion as follows.  At predetermined time instants, the agents exchange
information along the communication graph and update their internal state.
Between successive communication instants, the agents move according to a
motion control law, computed as a function of the agent location and of the
available information gathered through communication with other agents.  In
short, a \emph{control and communication law} for a robotic network
consists of a message-generation function (what do the agents
communicate?), a state-transition function (how do the agents update their
internal state with the received information?), and a motion control law
(how do the agents move between communication rounds?).  We then define the
notion of \emph{time complexity} of a control and communication law (aimed
at solving a given coordination task) as the minimum number of
communication rounds required by the agents to achieve the task.  The
\emph{time complexity of a coordination task}, as opposed to that of an
algorithm, is the minimum time complexity of any algorithm achieving the
task. We also provide similar definitions for notions of mean and total
communication complexity.  In particular, we show that our model of time
and communication complexity satisfies a basic well-posedness property that
we refer to as ``invariance under reschedulings.''

Building on these notions we establish complexity estimates for a few basic
motion coordination algorithms such as pursuit, rendezvous, and deployment.
First, for a network of agents moving on the circle, we introduce a dynamic
law, called \agreepursuitName, that combines elements of the classic cyclic
pursuit and leader election problems; see~\cite{JM-MB-BF:04c,NAL:97},
respectively. We show that this law achieves consensus on the agents'
direction of motion and equidistance between the agents' positions.
Furthermore, we show that these tasks are complete in time $O(N)$ and
$O(N^2\log N)$, respectively.
Second, we analyze a simple averaging law for a network of
locally-connected agents moving on a line. This law is related to the
widely known Vicsek's model, see~\cite{AJ-JL-ASM:02,TV-AC-EBJ-IC-OS:95}.
We show that this law achieves rendezvous (without preserving connectivity)
and that its time complexity belongs to $\Omega(N)$ and $O(N^5)$.
Third, for a network of locally-connected agents moving on a line or on a
segment, we show that the well-known circumcenter algorithm
by~\cite{HA-YO-IS-MY:99} has time complexity of order $\Theta(N)$. (This
algorithm achieves rendezvous while preserving connectivity with a
communication graph with $O(N^2)$ links.) We then consider a network based
on a different communication graph, called the limited Delaunay graph, that
arises naturally in computational geometry and in the study of wireless
communication topologies. For this less dense graph with $O(N)$
communication links, we show that the time complexity of the circumcenter
algorithm grows to $\Theta(N^2\log N)$.  For a network of agents moving on
$\real^d$ (with a certain communication graph) we introduce a novel
``parallel-circumcenter algorithm'' and establish its time complexity of
order $\Theta(N)$.
Fourth and last, for a network of agents in a one-dimensional environment,
we show that the time complexity of the deployment algorithm introduced
in~\cite{JC-SM-TK-FB:02j} is $O(N^3\log N)$.
To obtain these complexity estimates, we develop some novel analysis
methods. In particular, we develop a key set of results on tridiagonal
Toeplitz and circulant matrices that characterize their convergence rates
as a function of the matrices dimensions. 

We refer the reader to the Conclusions (Section~\ref{se:conclusions}) for a
discussion about the limitations of the proposed results and about avenues
for future research.

\paragraph*{Organization}
Section~\ref{se:modeling} presents a general approach to the modeling of
robotic networks by formally introducing various notions including, for
example, those of communication graph, control and communication law, and
network evolution.  Section~\ref{se:complexities} defines the notions of
task and of time and communication complexity for a control and
communication law.  Section~\ref{se:algorithms} presents a few motion
coordination algorithms performing the basic tasks of rendezvous and
deployment. For each algorithm, we characterize time and communication
complexity, along with asymptotic correctness.  Finally, we present our
conclusions in Section~\ref{se:conclusions}.  In the appendices, we review
some basic computational geometric structures and we prove some key facts
about the discrete-time dynamical systems defined by tridiagonal Toeplitz
and circulant matrices.

\paragraph*{Notation}
We let $\myboolean$ be the set $\{\true,\false\}$.  We let $
\prod_{i\in\until{N}}S_i$ denote the Cartesian product of sets $S_1,\dots,
S_N$. We let $\real_+$ and $\bar\real_+$ denote the set of strictly
positive and non-negative real numbers, respectively.  We let $\natural$
denote the set of positive natural numbers and $\naturalzero$ denote the
set of non-negative integers.  If $S$ is a set, then
$\diag(S\times{S})=\setdef{(s,s)\in S\times S}{s\in S}$.  For $x\in\real$,
we let $\fl{x}$ denote the floor of $x$. For $x\in\real^d$, we let
$\norm{x}{2}$ and $\norm{x}{\infty}$ denote the Euclidean and the
$\infty$-norm of $x$, respectively. Recall that
$\Infnorm{x}\leq\Enorm{x}\leq\sqrt{d}\Infnorm{x}$ for all $x\in\real^d$.
For $x\in\real^d$ and $r\in\real_+$, we let $\oball{r}{x}$ and
$\cball{r}{x}$ denote the open and closed ball in $\real^d$ centered at $x$
of radius $r$, respectively.  We let $\e_1,\dots,\e_d$ be the standard
orthonormal basis of $\real^d$.  Also, we define the vectors
$\mathbf{0}=(0,\dots,0)$ and $\mathbf{1}=(1,\dots,1)$ in $\real^d$.  For
$\map{f,g}{\natural}{\real}$, we say that $f \in O(g)$ (respectively, $f
\in \Omega(g)$) if there exist $N_0 \in \natural$ and $k \in \real_+$ such
that $|f(N)| \le k |g(N)|$ for all $N \ge N_0$ (respectively, $|f(N)| \ge k
|g(N)|$ for all $N \ge N_0$).  If $f \in O(g)$ and $f \in \Omega(g)$, then
we use the notation $f \in \Theta (g)$.  Finally, we refer the reader to
Appendix~\ref{app:geometry} for some useful geometric concepts.

\section{A formal model for synchronous robotic networks}
\label{se:modeling} 

In this section we introduce a notion of robotic network, as a group of
robotic agents with the ability to move and to communicate according to a
specified communication topology.  Our notion of control and communication
law for a robotic network parallels the definition in \cite{NAL:97} of
synchronous distributed algorithm.

\subsection{The physical components of a robotic network}
\label{se:physical-components}
Here we introduce our basic definition of physical quantities such as the
agents and such as the ability of agents to communicate.  We begin by
defining a control system as our basic model for how each robotic agent
moves in space.

\begin{definition}\label{dfn:ctrl-sys}
  A \emph{control system} is a tuple $(X,U,X_0,f)$ consisting of
  \begin{enumerate}
  \item $X$ is a differentiable manifold, called the \emph{state space};
  \item $U$ is a compact subset of $\real^m$ containing $0$, called the
    \emph{input space};
  \item $X_0$ is a subset of $X$, called the \emph{set of allowable initial
      states};
  \item $\map{f}{X\times{U}}{TX}$ is a $C^\infty$-map with $f(x,u)\in T_xX$
    for all $(x,u)\in X\times U$.  \oprocend
  \end{enumerate}
\end{definition}
We refer to $x \in X$ and $u \in U$ as a \emph{state} and an \emph{input}
of the control system, respectively.  We will often consider control-affine
systems, i.e., control systems for which $f(x,u) = f_0(x) + \sum_{a=1}^m
f_a (x)\, u_a$.  In such a case, with a slight abuse of notation, we will
represent the map $f$ as the ordered family of $C^\infty$-vector fields
$(f_0, f_1, \dots, f_m)$ on $X$.  We will also sometime consider
\emph{driftless} systems, i.e., control systems for which $f(x,0)=0$.

\begin{definition} \label{dfn:network-robotic-agents}
  A \emph{network of robotic agents} (or \emph{robotic network}) $\network$
  is a tuple $(I,\setofagents,\subscr{E}{cmm})$ consisting of
  \begin{enumerate}
  \item $I=\until{N}$; $I$ is called the \emph{set of unique identifiers
      (UIDs)};
   
  \item $\setofagents=\{\supind{\agent}{i}\}_{i\in{I}}=
    \{(\supind{X}{i},\supind{U}{i},\supind{X_0}{i},\supind{f}{i})\}_{i\in{I}}$
    is a set of control systems; this set is called the \emph{set of
      physical agents};
      
  \item $\subscr{E}{cmm}$ is a map from $\prod_{i \in I} \supind{X}{i}$ to
    the subsets of $I\times{I} \setminus \diag(I\times{I})$; this map is
    called the \emph{communication edge map}.
  \end{enumerate}
  If $\supind{\agent}{i}=(X,U,X_0,f)$ for all $i \in I$, then the robotic network
  is called \emph{uniform}.  \oprocend
\end{definition}

Let us comment on this definition and on how robotic agents communicate in
a robotic network $(I,\setofagents,\subscr{E}{cmm})$.
\begin{remarks}
  \begin{enumerate}
  \item By convention, we let the superscript $[i]$ denote the variables
    and spaces which correspond to the agent with unique identifier $i$;
    for instance, $\supind{x}{i} \in \supind{X}{i}$ and $\supind{x}{i}_{0}
    \in \supind{X}{i}_{0}$ denote the state and the initial state of agent
    $\supind{\agent}{i}$, respectively.  We refer to
    $(\supind{x}{1},\dots,\supind{x}{N}) \in \prod_{i \in I} \supind{X}{i}$
    as a \emph{state} of the network.
    
  \item The map $\subscr{E}{cmm}$ models the topology of the communication
    service between the agents.  In other words, at a network state
    $x=(\supind{x}{1},\dots,\supind{x}{N})$, two agents at locations
    $\supind{x}{i}$ and $\supind{x}{j}$ can communicate if the pair $(i,j)$
    is an edge in $\subscr{E}{cmm}(\supind{x}{1},\dots,\supind{x}{N})$.
    Accordingly, we refer to the pair
    $(I,\subscr{E}{cmm}(\supind{x}{1},\dots,\supind{x}{N}))$ as the
    \emph{communication graph} at $x$.  When and what agents communicate is
    discussed below.
    
    Maps of the form
    $\map{E}{\prod_{i\in{I}}\supind{X}{i}}{2^{I\times{I}\setminus\diag(I\times{I})}}$
    are called \emph{proximity edge maps}, they arise in wireless
    communication and in computational geometry, we refer the reader to
    Appendix~\ref{app:geometry} for more details.  Excluding edges of the
    form $(i,i)$, for $i\in{I}$, means that any individual agent does not
    communicate with itself. \oprocend
  \end{enumerate}
\end{remarks}

To make things concrete, let us now present some useful examples of robotic
networks; we will use these examples throughout the paper. We start with a
fairly common example and define some interesting variations.

\begin{example}[Locally-connected first-order agents in $\real^d$]
  \label{ex:network-reald-Er}
  Consider $N$ points $\supind{x}{1},\dots,\supind{x}{N}$ in the Euclidean
  space $\real^d$, $d\geq 1$, obeying a first-order dynamics
  $\supind{\dot{x}}{i}(t) = \supind{u}{i}(t)$.  According to
  Definition~\ref{dfn:ctrl-sys}, these are identical agents of the form
  $\agent = (\real^d,\real^d,\real^d,(\mathbf{0},\e_1,\dots,\e_d))$.
  Assume that each agent can communicate to any other agent within
  Euclidean distance $r$, that is, adopt as communication edge map the
  $r$-disk proximity edge map $\subscr{E}{\rdisk}$ defined in
  Appendix~\ref{app:geometry}.  These data define the uniform robotic
  network $\network[\realdisk]=(I,\setofagents,\subscr{E}{\rdisk})$.
  \oprocend
\end{example}

\begin{example}[LD-connected first-order agents in $\real^d$]
  \label{ex:network-reald-ELDr}
  Consider the set of physical agents defined in the previous example.  For
  $r\in\real_+$, recall from Appendix~\ref{app:geometry} the $r$-limited
  Delaunay map $\subscr{E}{\rLD}$ defined by
  \begin{align*}
    (i,j) \in \subscr{E}{\rLD}(\supind{x}{1},\dots,\supind{x}{N}) 
    \quad \text{if and only if} \quad  %%
    \big(\supind{V}{i} \cap \cball{\tfrac{r}{2}}{\supind{x}{i}} \big) \intersection
    \big(\supind{V}{j} \cap \cball{\tfrac{r}{2}}{\supind{x}{j}} \big) \neq \emptyset, \; i\neq j,
  \end{align*}
  where $\{\supind{V}{1},\dots,\supind{V}{N}\}$ is the Voronoi partition of
  $\real^d$ generated by $\{x^{[1]},\dots,x^{[N]}\}$.  These data define
  the uniform robotic network
  $\network[\realLD]=(I,\setofagents,\subscr{E}{\rLD})$.  \oprocend
\end{example}

\begin{example}[Locally-$\infty$-connected first-order agents in $\real^d$]
  \label{ex:network-reald-EIr}
  Consider the set of physical agents defined in the previous two examples.
  For $r\in\real_+$, define the proximity edge map $\subscr{E}{\rIdisk}$ by
  \begin{align*}
    (i,j) \in \subscr{E}{\rIdisk}
    (\supind{x}{1},\dots,\supind{x}{N}) 
    \quad \text{if and only if} \quad  %%
    \norm{\supind{x}{i} - \supind{x}{j} }{\infty} \le r, \; i\neq j.
  \end{align*}
  These data define the uniform robotic network
  $\network[\realIdisk]=(I,\setofagents,\subscr{E}{\rIdisk})$.  \oprocend
\end{example}

In what follows, let $\sph^1$ be the unit circle, and measure positions on
the circle counterclockwise from the positive horizontal axis.  Let us be
specific about distances on the circle and related concepts.  For
$x,y\in\sph^1$, we let $\dist(x,y)=\min\{\distC(x,y),\distCC(x,y)\}$.
Here, $\distC(x,y) = (x-y)\, (\mod 2\pi)$ is the clockwise distance, that
is, the path length from $x$ to $y$ traveling clockwise.  Similarly,
$\distCC(x,y) = (y-x)\, (\mod 2\pi)$ is the counterclockwise distance,
i.e., the path length from $x$ to $y$ traveling counterclockwise.  Here
$x\, (\mod 2\pi)$ is the remainder of the division of $x$ by $2\pi$.

\begin{example}[Locally-connected first-order agents on the
  circle]
  \label{ex:network-circle-Er}
  For $r\in\real_+$, consider the uniform robotic network
  $\network[\circledisk] = (I,\setofagents, \subscr{E}{\rdisk})$ composed
  of identical agents of the form $(\sph^1,(0,\e))$. Here $\e$ is the
  vector field on $\sph^1$ describing unit-speed counterclockwise rotation.
  We define the $r$-disk proximity edge map $\subscr{E}{\rdisk}$ on the
  circle by
  \begin{align*}
    (i,j) \in \subscr{E}{\rdisk}(\theta^{[1]},\dots,\theta^{[N]}) 
    \quad \text{if and only if} \quad %%
    \dist(\supind{\theta}{i}, \supind{\theta}{j} ) \le r \, ,
  \end{align*} 
  where $\dist(x,y)$ is the geodesic distance between the two points $x,y$
  on the circle.  \oprocend
\end{example}

%%%%%%%%%%%%%%%%%%%%%%%%%%%%%%%%%%%%%%%%%%%%%%%%%%%%%%%%%%%%
%%%%%%%%%%%%%%%%%%%%%%%%%%%%%%%%%%%%%%%%%%%%%%%%%%%%%%%%%%%%

\subsection{Control and communication laws for robotic networks}
\label{se:control-and-communication-laws}

Here we present a discrete-time communication, continuous-time motion model
for the evolution of a robotic network.  In our model, the robotic agents
evolve in the physical domain in continuous-time and have the ability to
exchange information (position and/or dynamic variables) that affect their
motion at discrete-time instants.

\begin{definition}[Control and communication law]
  \label{dfn:dynamic-feedback-control-communication}%
  Let $\network$ be a robotic network.  A \emph{(synchronous, dynamic,
    feedback) control and communication law} $\FCC$ for $\network$ consists
  of the sets:
  \begin{enumerate}
  \item $\timeschedule=\{t_\ell\}_{\ell\in\naturalzero} \subset
    \bar{\real}_+$ is an increasing sequence of time instants, called
    \emph{communication schedule};
    
  \item $L$ is a set containing the $\nll$ element, called the
    \emph{communication language}; elements of $L$ are called
    \emph{messages};
    
  \item $\supind{W}{i}$, $i\in{I}$, are sets of values of some \emph{logic
      variables} $\supind{w}{i}$, $i\in{I}$;
    
  \item $\supind{W}{i}_{0} \subseteq \supind{W}{i}$, $i\in{I}$, are subsets of
    \emph{allowable initial values};
  \end{enumerate}
  and of the maps:
  \begin{enumerate}    
  \item $\map{ \supind{\msg}{i} }%
    {\timeschedule\times\supind{X}{i} \times \supind{W}{i} \times I} {L}$,
    $i\in{I}$, are called \emph{message-generation functions};
    
  \item $\map{\supind{\stf}{i}}%
    {\timeschedule\times\supind{W}{i} \times L^N } {\supind{W}{i}}$,
    $i\in{I}$, are called \emph{state-transition functions};
    
  \item $\map{\supind{\ctrl}{i}}%
    {\bar{\real}_+\times\supind{X}{i} \times \supind{X}{i} \times
      \supind{W}{i} \times L^N } {\supind{U}{i}}$, $i\in{I}$, are called
    \emph{control functions}.   \oprocend
\end{enumerate}
\end{definition}

We will sometimes refer to a control and communication law as a
\emph{motion coordination algorithm}.  Control and communication laws might
have various properties.
\begin{definition}[Properties of control and communication laws]
Let $\network$ be a robotic network and $\FCC$ be a control and
  communication law for $\network$.
\begin{enumerate}
\item If $\network$ is uniform and if $\supind{W}{i}=W$,
  $\supind{\msg}{i}=\msg$, $\supind{\stf}{i}=\stf$,
  $\supind{\ctrl}{i}=\ctrl$, for all $i\in{I}$, then $\FCC$ is said to be
  \emph{uniform} and is described by a tuple $(\timeschedule,L,W,
  \{\supind{W}{i}_{0}\}_{i\in{I}},\msg,\stf,\ctrl)$.%
  
\item If $\supind{W}{i}=\supind{W}{i}_0=\emptyset$ for all $i\in{I}$, then
  $\FCC$ is said to be \emph{static} and is described by a tuple
  $(\timeschedule,L, \{\supind{\msg}{i}\}_{i\in{I}},
  \{\supind{\ctrl}{i}\}_{i\in{I}})$,  with  $\map{ \supind{\msg}{i}}%
  {\timeschedule\times\supind{X}{i} \times I} {L}$, and
  $\map{\supind{\ctrl}{i}}
  {\timeschedule\times\supind{X}{i}\times\supind{X}{i}\times L^N }
  {\supind{U}{i}}$.%
  
\item $\FCC$ is said to be \emph{time-independent} if the
  message-generation, state-transition and control functions are of the
  form
  $\map{ \supind{\msg}{i} }%
  {\supind{X}{i} \times \supind{W}{i} \times I} {L}$,
  $\map{\supind{\stf}{i}}%
  {\supind{W}{i} \times L^N } {\supind{W}{i}}$,
  $\map{\supind{\ctrl}{i}}%
  {\supind{X}{i} \times \supind{X}{i} \times \supind{W}{i} \times L^N }
  {\supind{U}{i}}$, $i\in{I}$, respectively.  \oprocend
\end{enumerate}
\end{definition}

Roughly speaking this definition has the following meaning: for all
$i\in{I}$, to the $i$th physical agent corresponds a logic process, labeled
$i$, that performs the following actions. First, at each time instant
$t_\ell\in\timeschedule$, the $i$th logic process sends to each of its
neighbors in the communication graph a message (possibly the $\nll$
message) computed by applying the message-generation function to the
current values of $\supind{x}{i}$ and $\supind{w}{i}$.  After a negligible
period of time (therefore, still at time instant $t_\ell \in
\timeschedule$), the $i$th logic process resets the value of its logic
variables $\supind{w}{i}$ by applying the state-transition function to the
current value of $\supind{w}{i}$, and to the messages received at time
$t_\ell$.  Between communication instants, i.e., for $t \in
[t_\ell,t_{\ell+1})$, the motion of the $i$th agent is determined by
applying the control function to the current value of $\supind{x}{i}$, the
value of $\supind{x}{i}$ at $t_\ell$, and the current value of
$\supind{w}{i}$.  This idea is formalized as follows.

\begin{definition}[Evolution of a robotic network subject to a
  control and communication law] \label{dfn:evolution}%
  Let $\network$ be a robotic network and $\FCC$ be a control and
  communication law for $\network$.  The \emph{evolution} of
  $(\network,\FCC)$ from initial conditions $\supind{x}{i}_{0}\in
  \supind{X_0}{i}$ and $\supind{w}{i}_{0}\in \supind{W}{i}_0$, $i\in{I}$,
  is the set of curves
  $\map{x^{[i],\ell}}{[t_\ell,t_{\ell+1}]}{\supind{X}{i}}$, $i\in{I}$,
  $\ell \in \naturalzero$, and
  $\map{\supind{w}{i}}{\timeschedule}{\supind{W}{i}}$, $i\in{I}$,
  satisfying
  \begin{align*}
    \dot{x}^{[i],\ell}(t) &= f\big( x^{[i],\ell}(t), \,
    \supind{\ctrl}{i}(t,x^{[i],\ell}(t),x^{[i],\ell}(t_\ell),
    \supind{w}{i}(t_{\ell}), y^{[i]}(t_\ell) ) \big) ,
  \end{align*}
  where, for $\ell\in\naturalzero$, and $i\in{I}$,
  \begin{align*}
    x^{[i],\ell}(t_\ell) = x^{[i],\ell -1}(t_\ell) \, , \quad
    \supind{w}{i}(t_\ell) =
    \supind{\stf}{i}(t_\ell,\supind{w}{i}(t_{\ell-1}), y^{[i]}(t_\ell)) \,
    ,
  \end{align*}
  with the conventions that $x^{[i],-1}(t_0)=\supind{x}{i}_{0}$ and
  $\supind{w}{i}(t_{-1})=\supind{w}{i}_{0}$, $i\in{I}$. Here, the function
  $\map{\supind{y}{i}}{\timeschedule}{L^N}$ (describing the messages
  received by agent $i$) has components $\supind{y}{i}_j(t_\ell)$, for
  $j\in{I}$, given by
  \begin{equation*}
    \supind{y}{i}_j(t_\ell) =
    \begin{cases}
      \msg^{[j]}(t_\ell,x^{[j],\ell-1}(t_\ell),w^{[j]}(t_{\ell-1}),i), &
      \text{if} \enspace (i,j)\in
      \subscr{E}{cmm}(x^{[1],\ell-1}(t_\ell),\dots,x^{[N],\ell-1}(t_\ell)),
      \\
      \nll, & \text{otherwise}.
  \end{cases}   \eqoprocend
\end{equation*}
\end{definition}

\medskip 

Let us emphasize two limitations regarding the proposed communication
model.

\begin{remarks}[Idealized aspects of communication model]
  \begin{enumerate}
  \item We refer to $\FCC$ as a \emph{synchronous} control and
    communication law, because the communications between all agents takes
    place always at the same time for all agents.  We do not discuss here
    the important setting of asynchronous laws. 
    
  \item The set $L$ is used to exchange information between two robotic
    agents. The message \nll indicates no communication.  We assume that
    the messages in the communication language $L$ allow us to encode
    logical expressions such as \true and \false, integers, and real
    numbers.  A realistic assumption on the communication language would be
    to adopt a finite-precision representation for integers and real
    numbers in the messages. Instead, in what follows, we neglect any
    inaccuracies due to quantization and we discuss this topic as an open
    problem in Section~\ref{se:conclusions}. \oprocend
  \end{enumerate}
\end{remarks}

In the following remarks we introduce additional useful notations and
emphasize an important assumption.

\begin{remarks}[Related concepts and notations]
  \begin{enumerate}    
  \item To distinguish between the $\nll$ and the non-$\nll$ messages
    received by an agent at a given time instant, it is convenient to
    define the \emph{natural projection} $\map{\piSpace{L}}{L^N}{2^{L}}$
    that maps an array of messages $y$ to the subset of $L$ containing only
    the non-$\nll$ messages in $y$.
       
  \item In many uniform control and communication laws, the messages
    interchanged among the network agents are (quantized representations
    of) the agents' states and dynamic states.  The corresponding
    communication language is $L=X\times{W}$ and message generation
    function
    $\map{\msgstandard}{\timeschedule\times{X}\times{W}\times{I}}{X\times{W}}$
    is referred to as the \emph{standard message-generation function} and
    is defined by $\msgstandard (t,x,w,j) = (x,w)$.
    
  \item By concatenating the curves $x^{[i],\ell}$ and $w^{[i],\ell}$, for
    $\ell\in\naturalzero$, we can define the evolution of the $i$th robotic
    agent $\bar{\real}_+ \ni t \mapsto (\supind{x}{i}(t),\supind{w}{i}(t))
    \in\supind{X}{i}\times\supind{W}{i}$. Additionally we can define the
    curves
    \begin{align*}
      \bar\real_+ \ni t & \mapsto x(t) = (x^{[1]}(t),\dots,x^{[N]}(t)) \in
      \prod_{i \in I} \supind{X}{i} ,
      \\
      \bar\real_+ \ni t & \mapsto w(t) = (w^{[1]}(t),\dots,w^{[N]}(t)) \in
      \prod_{i \in I} \supind{W}{i}.
    \end{align*}  
    
  \item In the proposed notion of synchronous robotic network, we assume
    that the states
    $(\supind{x}{1},\dots,\supind{x}{N})\in\prod_{i\in{I}}\supind{X}{i}$
    are not reset after a communication round, i.e., we are setting
    $x^{[i],\ell}(t_\ell)= x^{[i],\ell-1}(t_\ell)$.  This is reasonable for
    instance when $(\supind{x}{1},\dots,\supind{x}{N})$ are positions or
    other physical variables. It would be possible, though, to allow for
    more general situations by defining state-transition functions of the
    form $\map{\supind{\stf}{i}}%
    {\timeschedule\times\supind{X}{i}\times\supind{W}{i}\times{L^N}}
    {\supind{X}{i}\times\supind{W}{i}}$, for $i\in{I}$, and setting
    $(x^{[i],\ell},w^{[i],\ell})(t_\ell) = \supind{\stf}{i}
    (t_\ell,x^{[i],\ell -1}(t_\ell),
    w^{[i],\ell-1}(t_\ell),y^{[i]}(t_\ell))$.  Other possible extensions
    would be to consider (1) non-deterministic continuous-time evolutions
    of the state variables (by means of controlled differential inclusions,
    as opposed to the controlled differential equations) or (2)
    non-deterministic discrete-time updates of the dynamical variables (by
    means of set-valued state-transitions functions).  Further
    generalizations are also possible, but we do not consider them here in
    the sake of simplicity.  \oprocend
  \end{enumerate}
\end{remarks}

%%%%%%%%%%%%%%%%%%%%%%%%%%%%%%%%%%%%%%%%%%%%%%%%%%%%%%%%%%%%
%%%%%%%%%%%%%%%%%%%%%%%%%%%%%%%%%%%%%%%%%%%%%%%%%%%%%%%%%%%%

\subsection{Example control and communication laws}
It is now a good time to present various examples of control and
communication laws. We start by considering a very simple averaging
algorithm as a static control and communication law.  We then present an
interesting dynamic control and communication law on the circle.

\begin{example}[The \VicsekName control and communication law]
  \label{ex:rendezvous}%
  From Example~\ref{ex:network-reald-Er}, consider the uniform network
  $\network[\realdisk]$ of locally-connected first-order agents in
  $\real^d$.  We now define a static, uniform and time-independent law that
  we refer to as the \VicsekName law and that we denote by
  $\FCC[\VicsekSub]$.  We loosely describe it as follows:
  \begin{quote}
    \emph{[Informal description]} Communication rounds take place at each
    natural instant of time. At each communication round each agent
    transmits its position. Between communication rounds, each agent moves
    towards and reaches the point that is the average of its neighbors'
    positions; the average point is computed including the agent's own
    position.
  \end{quote}
  Next, we \emph{formally} define the law as follows.  First, we take
  $\timeschedule=\naturalzero$ and we assume that each agent operates with
  the standard message-generation function, i.e., we set $L=\real^d$ and
  $\msg(x,j)=\msgstandard(x,j)=x$.  Second, we define the control function
  $\map{\ctrl}{\real^d\times\real^d\times L^N}{\real^d}$ by
  \begin{align*}
    \ctrl(x,\sampled{x},y) &= - \kprop \, \vers \big(x -
    \average(y\union\{\sampled{x}\}) \big) ,
  \end{align*}
  where $\kprop \ge r$, $\map{\vers}{\real^d}{\real^d}$ is defined by
  $\vers(0)=0$ and $\vers(v)=v/\Enorm{v}$ for $v\neq 0$, and the map
  $\average$ computes the average of a finite point set in $\real^d$:
  \begin{equation*}
    \average(S) =  \frac{1}{\displaystyle \sum_{p \in
        \piSpace{\real}(S)} 1}  \; \sum_{p \in
      \piSpace{\real}(S)} \!  p  \, .
  \end{equation*}
  In summary we set
  $\FCC[\VicsekSub]=(\naturalzero,\real^d,\msgstandard,\ctrl)$.  An
  implementation of this control and communication law is shown in
  Fig.~\ref{fig:rendezvous-disk} for $d=1$.  Note that, along the
  evolution, (1) several agents \emph{rendezvous}, i.e., agree upon a
  common location, and (2) some agents are connected at the simulation's
  beginning and not connected at the simulation's end. Finally, we remark
  that this law is related to the Vicsek's model discussed
  in~\cite{AJ-JL-ASM:02,TV-AC-EBJ-IC-OS:95}.  \oprocend
   \begin{figure}[htbp]
     \centering
     \psfrag{0}{}\psfrag{1}{1}\psfrag{2}{2}\psfrag{3}{3}\psfrag{4}{4}
     \psfrag{8}{}\psfrag{7}{7}\psfrag{6}{6}\psfrag{5}{5}
     \includegraphics[width=.75\linewidth]%
     {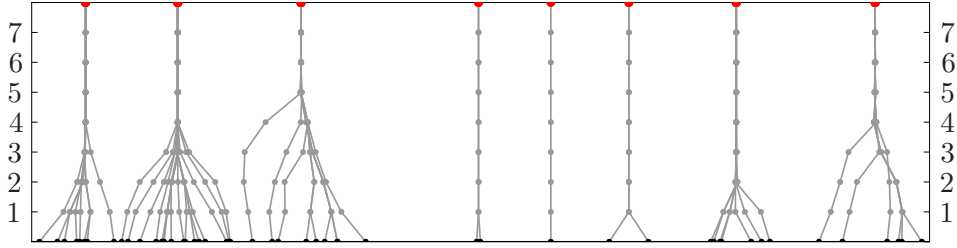}
     \caption{Evolution of a robotic network under the \VicsekName control
       and communication law in Example~\ref{ex:rendezvous} implemented
       over the $r$-disk graph, with $r=1.5$. The vertical axis corresponds
       to the elapsed time, and the horizontal axis to the positions of the
       agents in the real line. The $51$ agents are initially randomly
       deployed over the interval $[-15,15]$.}
     \label{fig:rendezvous-disk}
   \end{figure}
\end{example}

\begin{example}[The \agreepursuitName dynamic control and communication law]
  \label{ex:agree-pursuit}
  From Example~\ref{ex:network-circle-Er}, consider the uniform network
  $\network[\circledisk]$ of locally-connected first-order agents in
  $\sph^1$.  We now define the \agreepursuitName law, denoted by
  $\FCC[\agreepursuit]$, as the uniform and time-independent law loosely
  described as follows:
  \begin{quote}
    \emph{[Informal description]} The dynamic variables are $\dm$ taking
    values in $\{\csym,\ccsym\}$ and $\prior$ taking values in $I$.  At
    each communication round, each agent transmits its position and its
    dynamic variables and sets its dynamic variables to those of the
    incoming message with the largest value of \prior.  Between
    communication rounds, each agent moves in the counterclockwise or
    clockwise direction depending on whether its dynamic variable \dm is
    \ccsym or \csym.  For $\kprop\in ]0,\frac{1}{2}[$, each agent moves
    $\kprop$ times the distance to the immediately next neighbor in the
    chosen direction, or, if no neighbors are detected, $\kprop$ times the
    communication range $r$.
  \end{quote}
  Next, we \emph{formally} define the law as follows.  Each agent has logic
  variables $w=(\dm,\prior)$, where $w_1=\dm \in \{\ccsym,\csym\}$, with
  arbitrary initial value, and $w_2=\prior\in{I}$, with initial value equal
  to the agent's identifier $i$. In other words, we define
  $W=\{\ccsym,\csym\}\times{I}$, and we set
  $\supind{W}{i}_0=\{\ccsym,\csym\}\times\{i\}$.  Each agent $i \in {I}$
  operates with the standard message-generation function, i.e., we set $L =
  \sph^1\times W$ and $\supind{\msg}{i}=\msgstandard$, where
  $\msgstandard(\theta,w,j)=(\theta,w)$.  The state-transition function is
  defined by
  \begin{equation*}
    \stf(w,y) = \argmax\setdef{z_2}{z\in (\pi_L(y))_2\union\{w\}}.
  \end{equation*}
  For $\kprop\in\real_+$, the control function is
  \begin{equation*}
    \ctrl(\theta,\subscr{\theta}{smpld},w,y)= 
    \kprop 
    \begin{cases}
      \min \{r\}\union\setdef{ \distCC(\subscr{\theta}{smpld},\subscr{\theta}{rcvd}) }%
      {\subscr{\theta}{rcvd} \in (\pi_L(y))_1}, 
      & \text{if} \enspace \dm = \ccsym ,
      \\
      -\min \{r\}\union\setdef{ \distC(\subscr{\theta}{smpld},\subscr{\theta}{rcvd}) }%
      {\subscr{\theta}{rcvd} \in (\pi_L(y))_1}, & \text{if} \enspace \dm = \csym.
    \end{cases}
  \end{equation*}
  
  Finally, we sketch the control and communication in equivalent pseudocode
  language. This is possible for this example, and necessary for more
  complicated ones.  For example, the state-transition function is written
  as:
  \begin{center}{\begin{minipage}[c]{.7\textwidth}
        \texttt{function}\, \stf\,$((\dm,\prior),\, y)$
        \\
        \texttt{for} each non-\nll message
        $(\subscr{\theta}{rcvd},(\subscr{\dm}{rcvd},\subscr{\prior}{rcvd}))$
        in $y$:\\
        \null\quad \texttt{if} $(\subscr{\prior}{rcvd} > \prior)$, 
        \;\texttt{then} \\
        \null\qquad $\dm := \subscr{\dm}{rcvd}$  \\
        \null\qquad $\prior := \subscr{\prior}{rcvd}$\\
        \null\quad \texttt{endif}\\
        \texttt{endfor}\\
        \texttt{return}\, (\dm,\prior)
      \end{minipage}}\end{center}\medskip
  Similarly, the control function $\ctrl$ is written as:
  \begin{center}{\begin{minipage}[c]{.7\textwidth}
        \texttt{function} \ctrl\,$(\theta,\,\subscr{\theta}{smpld},\, (\dm,\prior),\, y)$\\
        $\subscr{d}{tmp}$ := r
        \\
        \texttt{for} each non-\nll message 
        $(\subscr{\theta}{rcvd},(\subscr{\dm}{rcvd},\subscr{\prior}{rcvd}))$
        in $y$:
        \\
        \null\quad \texttt{if} (\dm = \ccsym) \;AND\;
        $(\distCC(\subscr{\theta}{smpld},\subscr{\theta}{rcvd}) < \subscr{d}{tmp})$, \;\texttt{then} 
        \\
        \null\qquad $\subscr{d}{tmp} := \distCC(\subscr{\theta}{smpld},\subscr{\theta}{rcvd})$  
        \\
%%        \null\quad \texttt{endif}
%%        \\
        \null\quad \texttt{else if} (\dm = \csym) \;AND\;
        $(\distC(\subscr{\theta}{smpld},\subscr{\theta}{rcvd}) < \subscr{d}{tmp})$, \;\texttt{then} 
        \\
        \null\qquad $\subscr{d}{tmp} := \distC(\subscr{\theta}{smpld},\subscr{\theta}{rcvd})$  
        \\
        \null\quad \texttt{endif}
        \\
        \texttt{endfor}
        \\
        \texttt{if} (\dm = \ccsym), \;\texttt{then} 
        \texttt{return}\, $\kprop \subscr{d}{tmp}$,
        \;\texttt{else}
        \texttt{return}\, $-\kprop \subscr{d}{tmp}$\;  \texttt{endif}
      \end{minipage}}\end{center}
  \smallskip%
  An implementation of this control and communication law is shown in
  Fig.~\ref{fig:roundabout-1}.  Note that, along the evolution, all agents
  agree upon a common direction of motion and, after suitable time, they
  reach a uniform distribution. Finally, we remark that this law is related
  to the leader election algorithm discussed in~\cite{NAL:97}. \oprocend
   \begin{figure}[htbp]
     \centering
     \includegraphics[width=.19\linewidth]{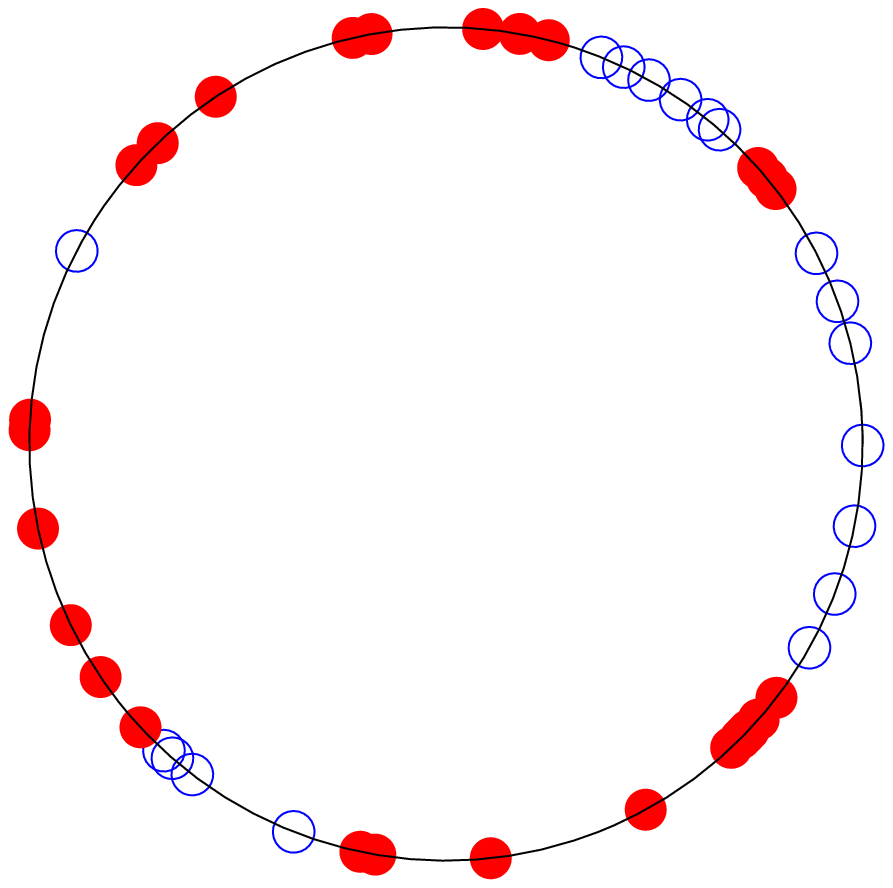}
     \includegraphics[width=.19\linewidth]{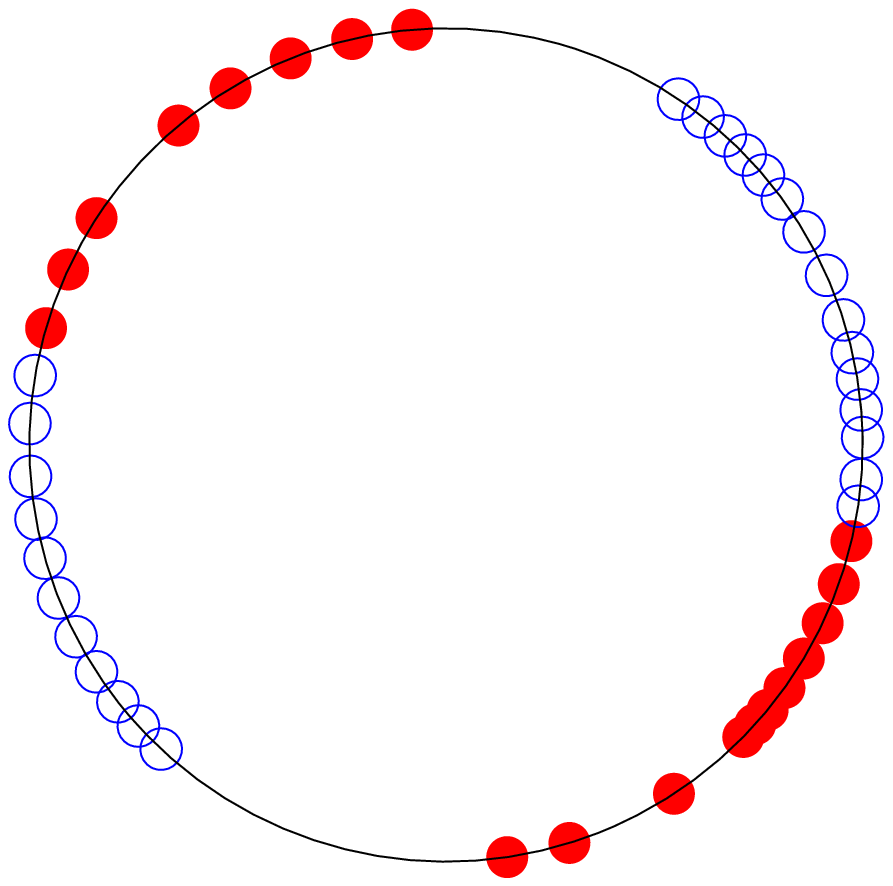}
     \includegraphics[width=.19\linewidth]{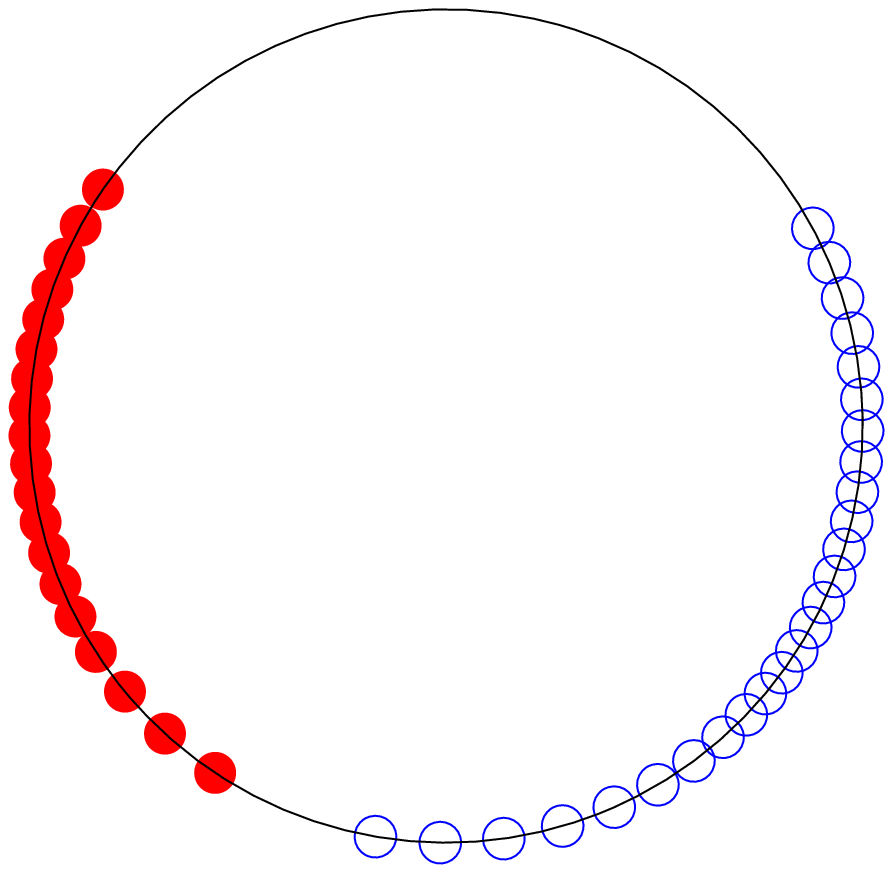}
     \includegraphics[width=.19\linewidth]{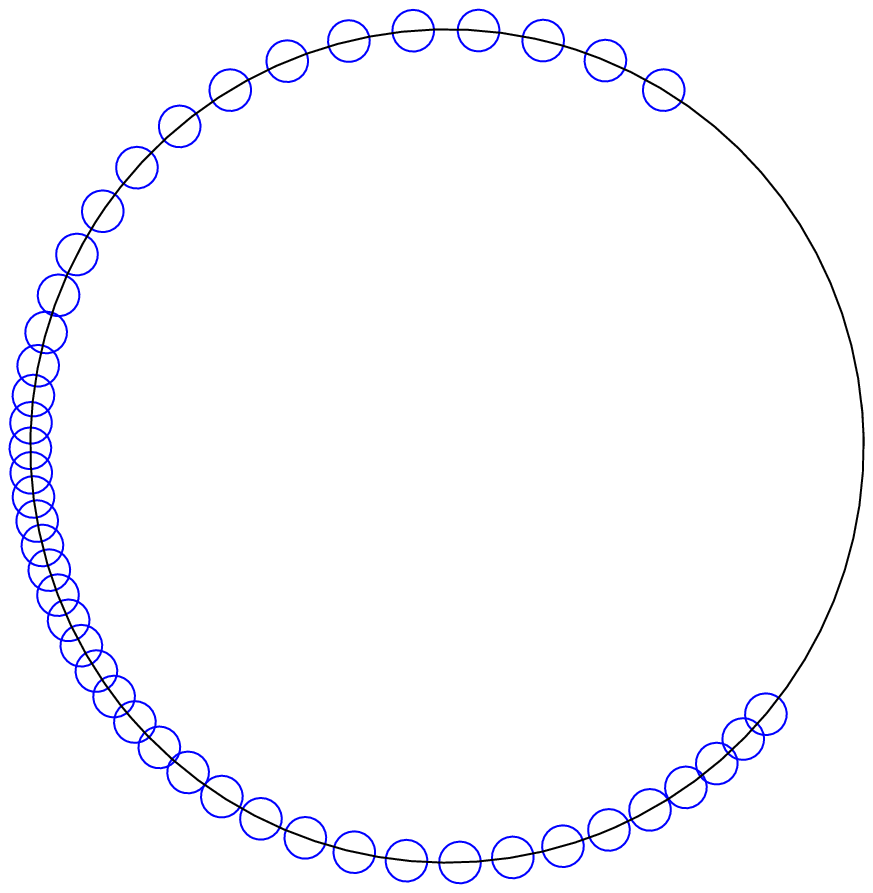}
     \includegraphics[width=.19\linewidth]{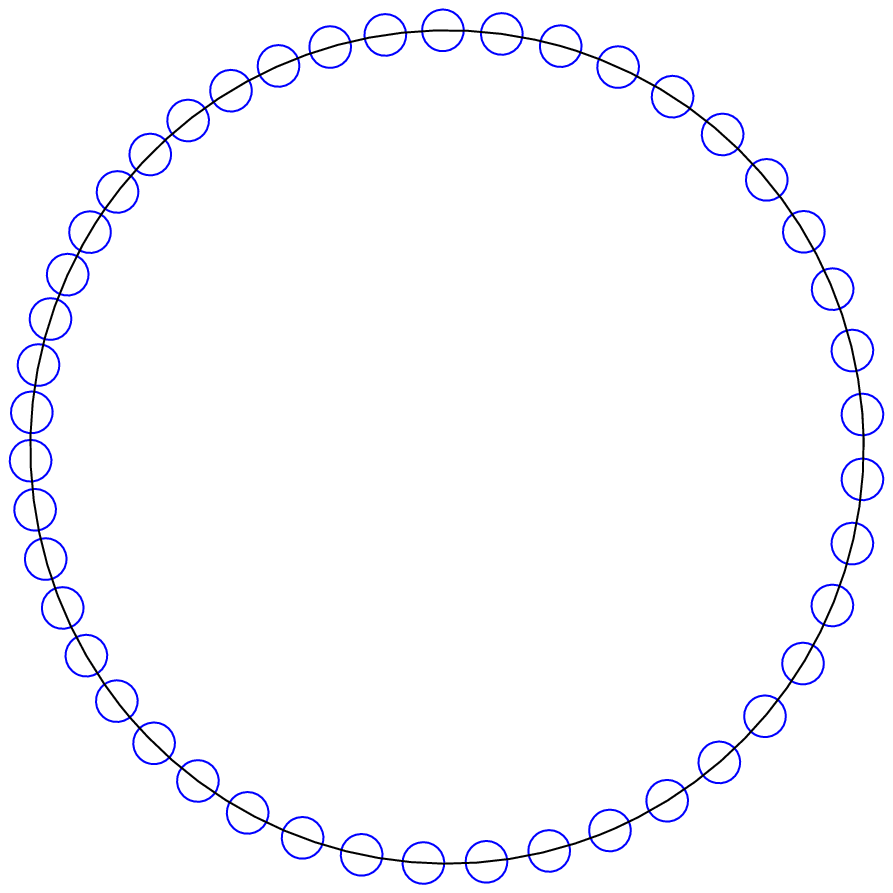}
     \caption{The \agreepursuitName control and communication law in
       Example~\ref{ex:agree-pursuit} with $N=45$, $r=2\pi/40$, and
       $\kprop=1/4$.  Disks and circles correspond to agents moving
       counterclockwise and clockwise, respectively.  The initial positions
       and the initial directions of motion are randomly generated. The
       five pictures depict the network state at times $0,12,37,100,400$.}
     \label{fig:roundabout-1}
   \end{figure}
\end{example}

\subsection{Groups of robotic agents with relative-position sensing}
\label{se:suzuki/yamashita}
In this last subsection on modeling, we discuss in some detail the
Suzuki-Yamashita model mentioned in the Introduction, see~\cite{IS-MY:99}.
This model consists of a group of identical mobile robots characterized as
follows: no explicit communication takes place between the agents, at each
instant of an ``activation schedule,'' each robot measure the relative
position of all other robots and moves according to a specified algorithm.
In this model, robots are referred to as ``anonymous'' and ``oblivious'' in
precisely the same way in which we defined control and communication laws
to be uniform and static, respectively.

As compared with our notion of robotic network, the Suzuki-Yamashita model
is more general in that the robots' activations schedules do not
necessarily coincide (i.e., this model is asynchronous), and at the same
time it is less general in that (1) robots cannot communicate any
information other than their respective positions, and (2) each robot
observes every other robot's position (i.e., the complete communication
graph is adopted; this limitation is not present for example in
\cite{HA-YO-IS-MY:99}).
Note that a control and communication law, as in our definition, can be
implemented on a synchronous Suzuki-Yamashita model if the law (1) is
static and uniform, (2) only relies on communicating the agents' positions
(e.g., the message-generation function is the standard one), and (4)
entails a control function that only depends on relative positions (as
opposed to absolute positions).

\section{Coordination tasks and complexity measures}\label{se:complexities}

In this section we introduce concepts and tools useful to analyze a
communication and control law.  We address the following issues: What is a
coordination task for a robotic network? When does a control and
communication law achieve a task? And with what time and communication
complexity?

\subsection{Coordination tasks}
Our first analysis step is to characterize the correctness properties of a
communication and control law.  We do so by defining the notion of task and
of task achievement by a robotic network.

\begin{definition}[Coordination task] \label{dfn:task}
  Let $\network$ be a robotic network and let $\WW$ be a set.
  \begin{enumerate}
  \item A \emph{coordination task} for $\network$ is a map $\map{\task}{
      \prod_{i\in{I}}\supind{X}{i} \times \WW^N }{\myboolean}$.
  \item If $\WW=\emptyset$, then the coordination task is said to be
    \emph{static} and is described by a map
    $\map{\task}{\prod_{i\in{I}}\supind{X}{i}}{\myboolean}$.
  \end{enumerate}
  Additionally, let $\FCC$ a control and communication law for $\network$.
  \begin{enumerate}
  \item The law $\FCC$ is \emph{compatible} with the task $\map{\task}{
      \prod_{i\in{I}}\supind{X}{i}\times\WW^N}{\myboolean}$ if its logic
    variables take values in $\WW$, that is, if $\supind{W}{i}=\WW$, for
    all $i\in{I}$.
  \item The law $\FCC$ \emph{achieves} the task $\task$ if it is compatible
    with it and if, for all initial conditions
    $\supind{x}{i}_{0}\in\supind{X}{i}_0$ and $\supind{w}{i}_0\in
    \supind{W}{i}_{0}$, $i\in I$, the corresponding network evolution
    $t\mapsto(x(t),w(t))$ has the property that there exists $T\in\real_+$
    such that $\task(x(t),w(t)) = \true$ for all $t\geq T$.  \oprocend
  \end{enumerate}
\end{definition} 

Loosely speaking, achieving a task might mean obtaining a specified pattern
in the position of the agents or of their dynamic variables.  We now give some
examples of interesting coordination tasks.
\begin{example}[Rendezvous tasks]
  \label{ex:rendezvous-tasks}
  First, let $\network=(I,\setofagents,\subscr{E}{cmm})$ be a uniform
  robotic network. The \emph{(exact) rendezvous task}
  $\map{\task[\rendezvous]}{X^N}{\myboolean}$ for $\network$ is the static
  task defined by
  \begin{equation*}
    \task[\rendezvous] (\supind{x}{1},\dots,\supind{x}{N}) =
    \begin{cases}
      \true, & \text{if} \enspace 
      \supind{x}{i}=\supind{x}{j}, \; \text{for all} \;
      (i,j)\in \subscr{E}{cmm}(\supind{x}{1},\dots,\supind{x}{N}),\\
      \false, & \text{otherwise} .
    \end{cases}
  \end{equation*} 
  
  Second, let $\network=(I,\setofagents,\subscr{E}{cmm})$ be a uniform
  robotic network with agents' state space $X\subset\real^d$. Examples
  networks of this form are $\network[\realdisk]$, see
  Examples~\ref{ex:network-reald-Er} and~\ref{ex:rendezvous}, and
  $\network[\realLD]$, see Examples~\ref{ex:network-reald-ELDr}.  For
  $\eps>0$, the \emph{$\eps$-rendezvous task}
  $\map{\task[$\eps$-\rendezvous]}{X^N}{\myboolean}$ for $\network$ is
  defined by
  \begin{equation*}
    \task[$\eps$-\rendezvous] (x) =
    \begin{cases}
      \true, & \text{if} \enspace 
      \Bignorm{
        \supind{x}{i} - 
        \average\big( 
        \{\supind{x}{i}\} \union 
        \setdef{\supind{x}{j}}{(i,j)\in \subscr{E}{cmm}(x)} 
        \big)
      }{2} < \eps , \;  \text{for all} \; i \in  {I},\\
      \false, & \text{otherwise} ,
    \end{cases}
  \end{equation*} 
  where $x =(\supind{x}{1},\dots,\supind{x}{N})\in X^N\subset(\real^d)^N$.
  In other words, $\task[$\eps$-\rendezvous]$ is \true at $x\in(\real^d)^N$
  if, for all $i\in{I}$, $\supind{x}{i}$ is at distance less than $\eps$
  from the average of its own position with the position of its
  $\subscr{E}{cmm}$-neighbors. \oprocend
\end{example}

\begin{example}[Agreement and equidistance tasks]
  \label{ex:roundabout-agreement}
  From Example~\ref{ex:network-circle-Er}, consider the uniform network
  $\network[\circledisk]$ of locally-connected first-order agents in
  $\sph^1$.  From Example~\ref{ex:agree-pursuit}, recall the
  \agreepursuitName control and communication law $\FCC[\agreepursuit]$
  with dynamic variables taking values in $W=\{\ccsym,\csym\}\times{I}$.
  There are two tasks of interest.  First, we define the \emph{agreement
    task} $\map{\task[\dm]}{(\sph^1)^N\times{W^N}}{\myboolean}$ by
  \begin{align*}
    \task[\dm](\theta,w) =
    \begin{cases}
      \true, & \text{if} \enspace
      \dm^{[1]}=\dots=\dm^{[N]}, \\
      \false, & \text{otherwise} ,
    \end{cases}
  \end{align*}
  where $\theta = (\theta^{[1]},\dots,\theta^{[N]})$,
  $w=(w^{[1]},\dots,w^{[N]})$, and $w^{[i]}=(\dm^{[i]},\prior^{[i]})$, for
  $i\in{I}$.  Furthermore, for $\eps>0$, we define the static
  \emph{$\eps$-equidistance task}
  $\map{\task[\equidistance]}{(\sph^1)^N}{\myboolean}$ by
  \begin{align*}
    \task[$\eps$-\equidistance] (\theta) =
    \begin{cases}
      \true, & \text{if} \enspace \big|\min_{j\neq i}
      \distC(\theta^{[i]},\theta^{[j]}) - \min_{j\neq i}
      \distCC(\theta^{[i]},\theta^{[j]})
      \big| < \eps, \; \text{for all} \; i\in{I}, \\
      \false, & \text{otherwise} .
    \end{cases} 
  \end{align*}
  In other words, $\task[$\eps$-\equidistance]$ is true when, for every
  agent, the clockwise distance to the closest clockwise neighbor and the
  counterclockwise distance to the closest counterclockwise neighbor are
  approximately equal. \oprocend
\end{example}

\begin{example}[Deployment tasks]
  \label{ex:deployment-tasks}
  By optimal deployment on the convex simple polytope $Q\subset\real^d$
  with density function $\map{\phi}{Q}{\real_+}$, we mean the following
  objective: place the agents on $Q$ so that the expected square Euclidean
  distance from any point in $Q$ to one of the agents is minimized.  To
  define this task formally, let us review some known preliminary notions;
  we will require some computational geometric notions from
  Appendix~\ref{app:geometry}.  We consider the following network objective
  function $\map{ \subscr{\HH}{\deployment} }{Q^N}{\real}$,
  \begin{equation}
    \label{eq:scenario-1}
    \subscr{\HH}{\deployment} (\supind{x}{1},\dots,\supind{x}{N}) 
    = \int_Q \min_{i\in{I}} \Enorm{ q - \supind{x}{i} }^2 \, \phi(q)dq \, .
  \end{equation}
  This function and variations of it are studied in the facility location
  and resource allocation research literature;
  see~\cite{AO-BB-KS-SNC:00,JC-SM-TK-FB:02j}. It is
  convenient~\cite{JC-SM-FB:03p-tmp} to study a generalization of this
  function.  For $r\in\real_+$, define the saturation function
  $\map{\sat{r}}{\real}{\real}$ by $\sat{r}(x)=x$ if $x\le{r}$ and
  $\sat{r}(x)=r$ otherwise.  For $r\in\real_+$, define the new objective
  function $\map{ \subscr{\HH}{$r$-\deployment} }{Q^N}{\real}$ by
  \begin{equation}
    \label{eq:scenario-2}
    \subscr{\HH}{$r$-\deployment} (\supind{x}{1},\dots,\supind{x}{N}) 
    = \int_Q \min_{i\in{I}} \sat{\frac{r}{2}}(\Enorm{ q - \supind{x}{i}
    }^2) \, \phi(q) dq \, .
  \end{equation}
  Note that if $r \ge 2 \diam(Q)$, then
  $\subscr{\HH}{\deployment}=\subscr{\HH}{$r$-\deployment}$.  Let
  $\{\supind{V}{1},\dots,\supind{V}{N} \}$ be the Voronoi partition of $Q$
  associated with $\{\supind{x}{1},\dots,\supind{x}{N}\}$.  The partial
  derivative of the cost function takes the following meaningful form
  \begin{align*}
    \pder{ \subscr{\HH}{$r$-\deployment} }{\supind{x}{i}}
    (\supind{x}{1},\dots,\supind{x}{N}) = %
    2 \Mass(\supind{V}{i} \intersect \cball{\tfrac{r}{2}}{\supind{x}{i}})
    \big( \Centroid(\supind{V}{i} \intersect \cball{\tfrac{r}{2}}{\supind{x}{i}}) -
    \supind{x}{i} \big) \, , \quad i \in I \, .
  \end{align*} 
  (Here, as in Appendix~\ref{app:geometry}, $\Mass(S)$ and $\Centroid(S)$
  are, respectively, the mass and the centroid of $S \subset \real^d$.)
  Clearly, the critical points of $\subscr{\HH}{$r$-\deployment}$ are
  network states where $\supind{x}{i} = \Centroid(\supind{V}{i} \intersect
  \cball{\frac{r}{2}}{\supind{x}{i}})$. We call such configurations
  $\frac{r}{2}$-centroidal Voronoi configurations.  For $r \ge 2 \diam(Q)$,
  they coincide with the standard centroidal Voronoi configurations on~$Q$.
  Fig.~\ref{fig:centroidals} illustrates these notions.
  \begin{figure}[htbp]
    \centering
    \includegraphics[width=.3\linewidth]{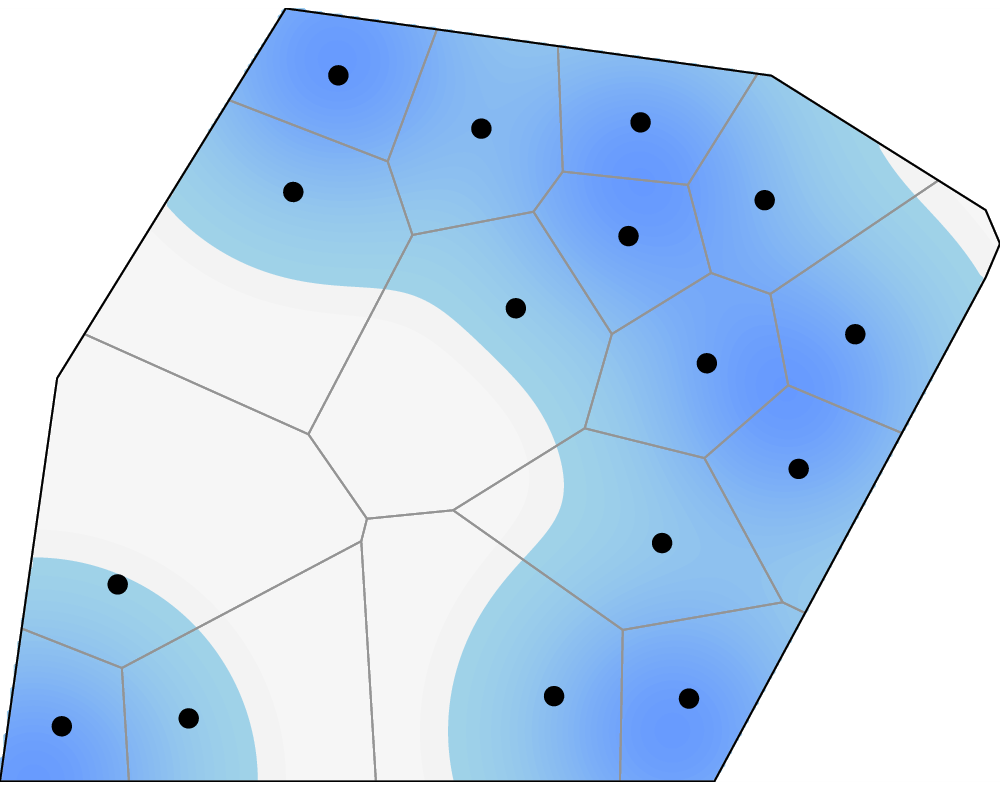}    
    \includegraphics[width=.3\linewidth]{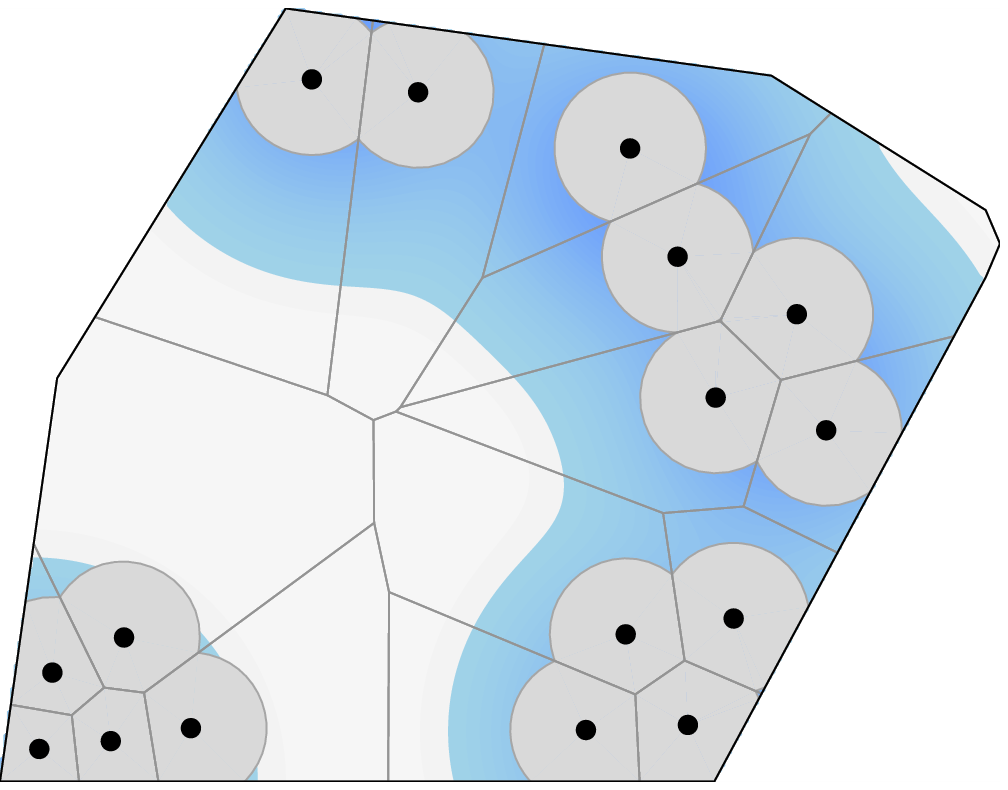}
    \caption{Centroidal and $\frac{r}{2}$-centroidal
      Voronoi configurations.  The density function $\phi$ is depicted by a
      contour plot.  For each agent $i$, the set $\supind{V}{i} \intersect
      \cball{\tfrac{r}{2}}{p_i}$ is plotted in light gray.}
    \label{fig:centroidals}
  \end{figure}
  
  Motivated by these observations, we define the following deployment task.
  For $r,\eps \in \real_+$, define the \emph{$\eps$-$r$-deployment task}
  $\map{\task[$\eps$-$r$-\deployment]}{Q^N}{\myboolean}$ by
  \begin{equation*}
    \task[$\eps$-$r$-\deployment](x) =
    \begin{cases}
      \true, & \text{if} \; %%%
      \bignorm{ \supind{x}{i} - 
        \Centroid( \supind{V}{i} \intersect
          \cball{\tfrac{r}{2}}{\supind{x}{i}}  ) }{2}
      \le \eps , 
      \; \text{for all} \; i \in I \, , \\
      \false, & \text{otherwise.}
    \end{cases}
  \end{equation*}
  Roughly speaking, $\task[$\eps$-$r$\deployment]$ is $\true$ for those
  network configurations where each agent is sufficiently close to the
  centroid of an appropriate region $\supind{V}{i} \intersect
  \cball{\tfrac{r}{2}}{\supind{x}{i}}$. Note that, with a natural
  definition of Voronoi partitions of the circle, the two tasks
  $\task[$\eps$-\equidistance]$ and $\task[$\eps$-$r$-\deployment]$ are
  closely related. \oprocend
\end{example}

%%%%%%%%%%%%%%%%%%%%%%%%%%%%%%%%%%%%%%%%%%%%%%%%%%%%%%%%%%%%%%%%
%%%%%%%%%%%%%%%%%%%%%%%%%%%%%%%%%%%%%%%%%%%%%%%%%%%%%%%%%%%%%%%%
%%%%%%%%%%%%%%%%%%%%%%%%%%%%%%%%%%%%%%%%%%%%%%%%%%%%%%%%%%%%%%%%

\subsection{Complexity notions for control and communication laws and for
  coordination tasks}

We are finally ready to define the key notions of time and communication
complexity. These notions describe the cost that a certain control and
communication law incurs while completing a certain coordination task.  We
also define the complexity of a task to be the infimum of the costs
incurred by all laws that achieve that task.

First we define the time complexity of an achievable task as the minimum
number of communication rounds needed by the agents to achieve the task
$\task$.

\begin{definition}[Time complexity] \label{dfn:time-complex}
  Let $\network$ be a robotic network and let $\task$ be a coordination
  task for $\network$.  Let $\FCC$ be a control and communication law for
  $\network$ compatible with $\task$.
  \begin{enumerate}
  \item The \emph{time complexity to achieve $\task$ with $\FCC$ from
      $(x_0,w_0) \in \prod_{i\in{I}}\supind{X}{i}_0\times\prod_{i\in
        I}\supind{W}{i}_{0}$} is
    \begin{equation*}
      \TC(\task,\FCC,x_0,w_0) = \inf \,
      \{\ell\;|\;\task(x(t_k),w(t_k))=\true\,,\;
      \text{for all} \; k \ge \ell\} \,,
    \end{equation*}
    where $t\mapsto(x(t),w(t))$ is the evolution of $(\network,\FCC)$ from
    the initial condition $(x_0,w_0)$.
    
  \item The \emph{time complexity to achieve $\task$ with $\FCC$} is
    \begin{equation*}
      \TC (\task,\FCC) = \sup 
      \Bigsetdef
      {\TC(\task,\FCC,x_0,w_0)}
      {(x_0,w_0) \in \prod_{i\in{I}}\supind{X}{i}_0\times\prod_{i\in I}\supind{W}{i}_{0}} 
       \, .  
    \end{equation*}
    
  \item The \emph{time complexity of $\task$} is
    \begin{equation*}
      \TC(\task) = 
      \inf \setdef{ \TC (\task,\FCC)}
      {\FCC \; \text{is compatible with $\task$}} \,.
       \eqoprocend
    \end{equation*}    
  \end{enumerate}
\end{definition}

Next, we define the notion of mean and total communication complexities for
a task.  As usual, we assume that the network $\network$ has a
communication edge map $\subscr{E}{cmm}$ and that the control and
communication law $\FCC$ has language $L$ and message-generation functions
$\supind{\msg}{i}$, $i\in{I}$.  With these data we can discuss the
communication cost of realizing one communication round.  At time
$t\in\timeschedule$ from state $(x,w)\in
\prod_{i\in{I}}\supind{X}{i}\times\prod_{i\in I}\supind{W}{i}$, an element
of $L$ needs to be transmitted for each edge of the directed graph
$(I,\subscr{E}{\nonnllmsgs}(t,x,w))$ defined by
\begin{equation*}
  (i,j)\in\subscr{E}{\nonnllmsgs}(t,x,w) 
  \quad \text{if and only if} \quad  %%
  (i,j)\in \subscr{E}{cmm}(x) \enspace\text{and}\enspace
  \supind{\msg}{i}(t,\supind{x}{i},\supind{w}{i},j) \neq \nll.
\end{equation*}
Next, we need a model for the cost of sending a message for each directed
edge in $\subscr{E}{\nonnllmsgs}$.
\begin{definition}[One-round cost] \label{dfn:one-round}
  \begin{enumerate}
  \item For $I=\until{N}$, a function
    $\map{\CommCost}{2^{I\times{I}}}{\bar\real_+}$ is a \emph{one-round
      cost function} if $\CommCost(\emptyset)=0$, and
    $S_1\subset{S_2}\subset I\times{I}$ implies
    $\CommCost(S_1)\leq\CommCost(S_2)$.    
  \item A one-round cost function $\CommCost$ is \emph{additive} if, for
    all $S_1,S_2\subset{I\times{I}}$, $S_1\intersect S_2=\emptyset$ implies
    $\CommCost(S_1\union{S_2})=\CommCost(S_1)+\CommCost(S_2)$. \oprocend
  \end{enumerate}
\end{definition}
We postpone our discussion about specific functions $\CommCost$ to the next
subsection.  Here we only emphasize that, for a given control and
communication law $\FCC$ with language $L$, the one-round cost depends on
$L$; we therefore write it as
$\map{\CommCost[L]}{2^{I\times{I}}}{\bar\real_+}$.

\begin{definition}[Communication complexity] \label{dfn:comm-complex}
  Let $\network$ be a robotic network and let $\task$ be a coordination
  task for $\network$.  Let $\FCC$ be a control and communication law for
  $\network$ compatible with $\task$, and let
  $\map{\CommCost[L]}{2^{I\times{I}}}{\bar\real_+}$ be a one-round cost
  function.
  \begin{enumerate}    
  \item Let $(x_0,w_0) \in \prod_{i\in{I}}\supind{X}{i}_0\times\prod_{i\in
      I}\supind{W}{i}_{0}$.  The \emph{mean communication complexity to
      achieve $\task$ with $\FCC$ from $(x_0,w_0)$} and the \emph{total
      communication complexity to achieve $\task$ with $\FCC$ from
      $(x_0,w_0) \in \prod_{i\in{I}}\supind{X}{i}_0\times\prod_{i\in
        I}\supind{W}{i}_{0}$} are, respectively,
    \begin{align*}
      \MCC(\task,\FCC,x_0,w_0) &= 
      \frac{1}{\TC(\FCC,\task,x_0,w_0)} 
      \sum_{\ell=0}^{\TC(\task,\FCC,x_0,w_0)-1} 
      \CommCost[L]\circ
      \subscr{E}{\nonnllmsgs}(t_\ell,x(t_\ell),w(t_\ell)),
      \\
      \TCC(\task,\FCC,x_0,w_0) &= 
      \sum_{\ell=0}^{\TC(\task,\FCC,x_0,w_0)-1} 
      \CommCost[L]\circ \subscr{E}{\nonnllmsgs}(t_\ell,x(t_\ell),w(t_\ell)),
    \end{align*}
    where $t\mapsto(x(t),w(t))$ is the evolution of $(\network,\FCC)$ from
    the initial condition $(x_0,w_0)$.  (Here $\MCC$ is defined only for
    $(x_0,w_0)$ with the property that $\task(x_0,w_0)=\false$.)
    
  \item The \emph{mean communication complexity to achieve $\task$ with
      $\FCC$} and the \emph{total communication complexity to achieve
      $\task$ with $\FCC$} are, respectively,
    \begin{align*}
      \MCC (\task,\FCC) &= \sup\bigsetdef
      {\MCC(\task,\FCC,x_0,w_0)}
      {(x_0,w_0) \in \prod_{i\in{I}}\supind{X}{i}_0\times\prod_{i\in I}\supind{W}{i}_{0}} 
       \, ,
       \\
      \TCC (\task,\FCC) &= \sup\bigsetdef
      {\TCC(\task,\FCC,x_0,w_0)}
      {(x_0,w_0) \in \prod_{i\in{I}}\supind{X}{i}_0\times\prod_{i\in I}\supind{W}{i}_{0}} 
       \, .  
    \end{align*}
    
  \item The \emph{mean communication complexity of $\task$} and the
    \emph{total communication complexity of $\task$} are, respectively,
    \begin{align*}
      \MCC(\task) &= 
      \inf \setdef{ \MCC (\task,\FCC)}
      {\FCC \; \text{is compatible with $\task$}} \,,
      \\
      \TCC(\task) &= 
      \inf \setdef{ \TCC (\task,\FCC)}
      {\FCC \; \text{is compatible with $\task$}} \,.
       \eqoprocend
    \end{align*}    
  \end{enumerate}
\end{definition}

We conclude this subsection with some remarks.
\begin{remarks}
  \begin{enumerate}
  \item The total communication complexity is equal to the average
    transmission cost during the execution multiplied by the number of
    rounds required to carry out the desired task.  That is, for $(x_0,w_0)
    \in \prod_{i\in{I}}\supind{X}{i}_0\times\prod_{i\in
      I}\supind{W}{i}_{0}$,
    \begin{equation*}
      \TCC (\task,\FCC,x_0,w_0) = 
      \MCC (\task,\FCC,x_0,w_0) \cdot \TC (\task,\FCC,x_0,w_0).
    \end{equation*}
    In turn this implies that $\TCC(\task,\FCC) \leq \MCC(\task,\FCC) \cdot
    \TC(\task,\FCC)$.

  \item According to this notation, given a robotic network $\network$ and
    a control and communication law $\FCC$, the time complexity of
    achieving a task $\task$ with $\FCC$ is $\TC(\task,\FCC)\in O(f)$
    (resp.~$\TC(\task,\FCC) \in \Omega(f)$), if there exist $N_0 \in
    \natural$ and $k \in \real_+$ such that $\TC(\task,\FCC,x_0,w_0) \le k
    f(N)$ for \emph{all} initial conditions $(x_0,w_0)$ for each $N \ge
    N_0$ (resp.~if $\TC (\task,\FCC,x_0,w_0) \ge k f(N)$ for at least
    \emph{an} initial condition $(x_0,w_0)$ for each $N \ge N_0$).
    
  \item A different notion of communication complexity is defined
    in~\cite{EK:02a} for a different robotic network model.  Transcribed to
    the current setting, this notion of communication complexity of the
    execution of a control and communication law $\FCC$ from initial
    conditions $(x_0,w_0)$ would read as
    \begin{equation}\label{eq:cc-Klavins}
      \cc (\FCC,x_0,w_0) = \lim_{k \to +\infty} \frac{1}{k} 
      \sum_{\ell=0}^{k} 
      \CommCost[L]\circ\subscr{E}{\nonnllmsgs}(t_\ell,x(t_\ell),w(t_\ell))
      \, .
    \end{equation}
    Note that this definition does not make reference to the completion of
    a task.  We will later come back to this notion in
    Section~\ref{se:invariance-rescheduling}. \oprocend
  \end{enumerate}
\end{remarks}

\subsection{Communication costs in unidirectional and omnidirectional
  wireless channels}

In this subsection we discuss some modeling aspects of the one-round
communication cost function $I\times{I}\supset E \mapsto \CommCost(E)$
described in Definition~\ref{dfn:one-round}.  First, let us mention that
the definition is motivated by the assumptions that (1) the cost of
exchanging any message is bounded, and that (2) this cost is zero only for
the \nll message.  More specific detail about the communication cost
depends necessarily on the type of communication service available between
the agents.

In \emph{unidirectional} models of communication messages are sent in a
point-to-point-wise fashion.  Certain forms of communication, such as those
based on the TCP-IP protocol, and certain technologies, such as wireless
networks equipped with unidirectional
antennas~\cite{YBK-VS-NHV:00,TK-GJ-LT:03}, fall into this category.  On the
other hand, in an \emph{omnidirectional} model of communication (e.g.,
wireless networks equipped with omnidirectional antennas), a single
transmission made by a node can be heard by several other nodes at the same
time.  This has the advantage that, by choosing a sufficiently large
transmission power, a signal can reach all the neighboring nodes in a
single time instant.

Broadly speaking, it is very difficult to come up with an abstract model
that captures adequately the cost of all possible communication
technologies.  For example, networking protocols for omnidirectional
wireless networks rely on a many nested layers to handle, for example,
media access, power control, congestion control, and routing.  The presence
of these layers and the non-trivial interactions between them make it
difficult to assess communication costs of individual messages. Let us
elaborate on this point in the following remark.

\begin{remark}[Omnidirectional wireless communication] 
  \label{re:layers-omnidir-comm}
  The \emph{Minimum Power Broadcast} (MPB) problem, the \emph{Medium Access
    Control} (MAC) problem, and their relationship are subjects of vigorous
  research in the wireless communications literature, see for
  instance~\cite{PRK:01} and references therein.  Loosely speaking, the MPB
  problem consists of finding, for each agent $i$, the minimum broadcast
  radius $\supind{R}{i}$ such that if agent $i$ sends a message with
  communication radius $\supind{R}{i}$, then its neighbors in a given graph
  $E$ receive it.  The MAC problem consists of determining a minimum number
  of broadcasting turns required for all agents to communicate their
  messages without interference.  A schematic approach to these problems is
  as follows: first, from the communication graph $(I,E)$, one constructs
  the \emph{neighbor-induced} graph $(I,E_{\mathcal{N}})$ by
  \begin{align*}
    (i,j) \in E_{\mathcal{N}} \quad \text{if and only if} \quad%
    (i,j) \in E \; \text{or} \; (i,k), (j,k) \in E \, , \; \text{for some
      $k\in{I}$}
  \end{align*}
  In the new graph $(I,E_{\mathcal{N}})$, the set of neighbors of the agent
  $i$ is composed by its neighbors in the graph $(I,E)$, together with the
  their respective neighbors. As a second step, one has to compute the
  \emph{chromatic number} of the graph, i.e., the minimum number of colors
  $\chi(E_{\mathcal{N}})$ needed to color the agents in such a way that
  there are no two neighboring agents with the same color. (This is also
  referred to as the \emph{coloring-graph problem}.)  Theorem~5.2.4
  in~\cite{RD:00} asserts that if a connected graph is neither complete,
  nor an odd cycle, then $\chi(E_{\mathcal{N}})$ is less than or equal to
  the maximum valency of the graph. Once the chromatic number has been
  determined, broadcasting turns can be established according to an ordered
  sequence of the agents' colors. Although this approach is clearly
  inadequate, it provides some basic pointers with regards to communication
  costs. \oprocend
\end{remark}

Motivated by the difficulty of obtaining a detailed model, the rest of this
paper relies on the following simplified models that capture some broad
relevant aspects:
\begin{enumerate}
\item For a unidirectional communication model, $E\mapsto\CommCost(E)$ is
  proportional to the total number of non-\nll messages sent over the
  directed edges in $E$, that is, $\CommCost(E)= c_0 \cdot
  \operatorname{cardinality}(E)$, where $c_0\in\real_+$ is the cost of
  sending a single message.  This one-round cost function is additive.
  This number is trivially bounded by twice the number of edges of the
  complete graph, which is $N(N-1)$.  Therefore, for unidirectional models
  of communication, we have $\MCC(\FCC,\task) \in O(N^2)$.
  
  \margin{for FB: here and in what follows, define subscript ``unidir'' and
    use it with $\CommCost, \MCC, \TCC$}

\item For an omnidirectional communication model, $E\mapsto\CommCost(E)$ is
  proportional to the number of turns employed to complete a communication
  round without interference between the agents (see
  Remark~\ref{re:layers-omnidir-comm}).  This number is trivially upper bounded
  by $N$. Therefore, for omnidirectional models of communication, we have
  $\MCC(\FCC,\task) \in O(N)$. 
\end{enumerate}
\margin{for all: for omnidirectional comm,  this should be a communication
  or energy cost, and instead here we talk about communication turns!}

\subsection{Invariance under rescheduling of control and
  communication laws}\label{se:invariance-rescheduling}

In this section, we discuss the invariance properties of the notions of
time and communication complexity under the \emph{rescheduling} of a
control and communication law. The idea behind rescheduling is to
``spread'' the execution of the law over time without affecting the
trajectories described by the robotic agents of the network. There are at
least two natural ways of doing this. One possible way consists of the
network slowing down its motion, and letting some communication rounds pass
without effectively interchanging any messages. Another possible way is to
schedule the messages originally sent at a single time instant to be sent
over multiple consecutive time instants, and adapt the motion of the
network accordingly.  Our objective is here is to formalize these ideas and
to examine the effect that these processes have on the notions of
complexity introduced earlier.  For simplicity we consider the setting of
static laws; similar results can be obtained for the general setting.

Let $\network=(I,\setofagents,\subscr{E}{cmm})$ be a robotic network with
driftless physical agents, that is, a robotic network where each physical
agent is a driftless control system. Let $\FCC =
(\naturalzero,L,\{\supind{\ctrl}{i}\}_{i\in{I}},\{\supind{\msg}{i}\}_{i\in{I}})$
be a static control and communication law.  It is out intention to define a
new control and communication law by modifying $\FCC$; to do so we
introduce some notation.  Let $s\in\natural$, with $s \le N$, and let
$\Pc_I = \{I_0,\dots,I_{s-1}\}$ be an \emph{$s$-partition} of $I$, that is,
$I_0,\dots,I_{s-1}$ are disjoint and nonempty subsets of $I$ and $I =
\union_{k=0}^{s-1} I_k$.

\begin{figure}[htb] 
  \begin{center}
    \resizebox{.9\linewidth}{!}{\input{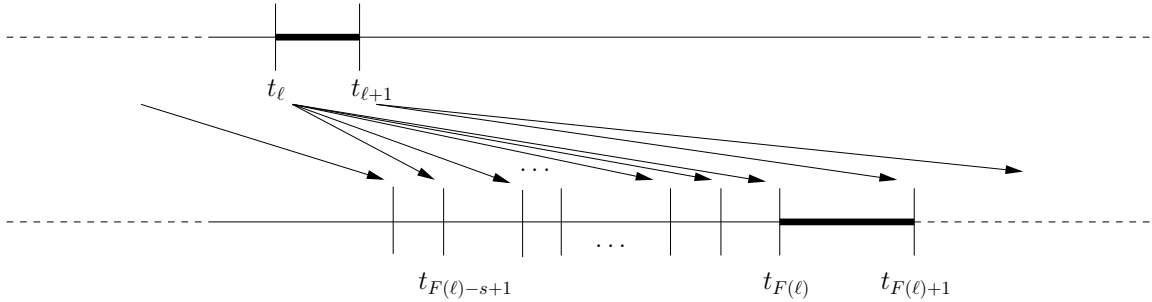}}
  \end{center}
  \caption{Under the rescheduling, the
    messages that are sent at the time instant $t_\ell$ under the control
    and communication law $\FCC$ are rescheduled to be sent over the time
    instants $t_{F(\ell)-s+1}, \dots, t_{F(\ell)}$ under the control and
    communication law $\reschedule{\FCC}{s}{\Pc_I}$.  Accordingly, the
    evolution of the robotic network under the original law during the time
    interval $[t_\ell,t_{\ell+1}]$ is now executed when all the
    corresponding messages have been transmitted, i.e., along the time
    interval $[t_{F(\ell)},t_{F(\ell)+1}]$.}\label{fig:stretching}
\end{figure}

For $i \in I$, define the message-generation functions
$\map{ \supind{\msg}{i}_{(s,\Pc_I)} }%
{\naturalzero\times\supind{X}{i} \times I} {L}$ by
\begin{align}
  \label{eq:rescheduled-message}
  \supind{\msg}{i}_{(s,\Pc_I)} (t_\ell, x,j) & =
  \begin{cases}
    \supind{\msg}{i}(t_{\fl{\ell/s}},x,j), \quad & \text{if}\enspace i \in I_k
    \enspace \text{and} \enspace
    k= \ell (\mod s)\, ,\\
    \nll, & \text{otherwise} \, .
  \end{cases}
\end{align}
According to this new message-generation function, only the agents with
unique identifier in $I_k$ will send messages at time $t_{\ell}$, with
$\ell\in\setdef{k+as}{a\in\naturalzero}$.  Equivalently, this can be stated
as follows.  Define the increasing function
$\map{F}{\naturalzero}{\naturalzero}$ by $F(\ell)=s(\ell+1)-1$.  According
to the message-generation functions specified
by~\eqref{eq:rescheduled-message}, the messages originally sent at the time
instant $t_\ell$ are now rescheduled to be sent at the time instants
$t_{F(\ell)-s+1}, \dots, t_{F(\ell)}$.  Fig.~\ref{fig:stretching}
illustrates this idea.

For $i \in I$, define the control functions
$\map{\supind{\ctrl}{i}}%
{\bar{\real}_+\times\supind{X}{i}\times\supind{X}{i}\times L^N }
{\supind{U}{i}}$ by
\begin{align}\label{eq:rescheduled-control}
  \supind{\ctrl}{i}_{(s,\Pc_I)} (t,x,\sampled{x},y)&=
  \begin{cases}
    \frac{t_{F^{-1}(\ell)+1}-t_{F^{-1}(\ell)}}{t_{\ell+1} -t_{\ell}}
    \, \supind{\ctrl}{i} \left( h_{\ell}(t),x,\sampled{x},y \right), \quad &
    \text{if} \; t \in [t_{\ell},t_{\ell+1}] \;
    \text{and} \; \ell = -1 (\mod s) \, ,\\
    0, & \text{otherwise} \, ,
  \end{cases}
\end{align} 
where $\map{F^{-1}}{\naturalzero}{\naturalzero}$ is the inverse of $F$,
defined by $F^{-1}(\ell) = \frac{\ell+1}{s} -1$, and for $\ell = -1 (\mod
s)$, the function $\map{h_{\ell}} {[t_{\ell},t_{\ell+1}]}
{[t_{F^{-1}(\ell)},t_{F^{-1}(\ell)+1}]}$ is the time re-parameterization
function defined by
\begin{align*}
  h_{\ell}(t) = \frac{(t_{F^{-1}(\ell)+1}-t_{F^{-1}(\ell)}) t +
    t_{\ell+1} t_{F^{-1}(\ell)} - t_{\ell} \,
    t_{F^{-1}(\ell)+1}}{t_{\ell+1}-t_{\ell}} \, , \quad t \in
  [t_{\ell},t_{\ell+1}] \, .
\end{align*}
Roughly speaking, the control law $\supind{\ctrl}{i}_{(s,\Pc_I)}$ makes the
agent $i$ wait for the time intervals $[t_\ell,t_{\ell+1}]$, with
$\ell\in\setdef{as-1}{a\in\natural}$, to execute any motion.  Accordingly,
the evolution of the robotic network under the original law $\FCC$ during
the time interval $[t_\ell,t_{\ell+1}]$ now takes place when all the
corresponding messages have been transmitted, i.e., along the time interval
$[t_{F(\ell)},t_{F(\ell)+1}]$.

We gather the above construction in the following definition.
\begin{definition}[Rescheduling of control and communication laws]
  \label{dfn:rescheduling}
  Let $\network=(I,\setofagents,\subscr{E}{cmm})$ be a robotic network with
  driftless physical agents, and let $\FCC =
  (\naturalzero,L,\{\supind{\ctrl}{i}\}_{i\in{I}},\{\supind{\msg}{i}\}_{i\in{I}})$
  be a static control and communication law.  Let $s \in \natural$, with $s
  \le N$, and let $\Pc_I$ be an $s$-partition of $I$.  The control and
  communication law $\reschedule{\FCC}{s}{\Pc_I} = (\naturalzero,L,
  \{\supind{\ctrl}{i}_{(s,\Pc_I)}\}_{i\in{I}},\{\supind{\msg}{i}_{(s,\Pc_I)}
  \}_{i\in I})$ defined by equations~\eqref{eq:rescheduled-message}
  and~\eqref{eq:rescheduled-control} is called a
  \emph{$(s,\Pc_I)$-rescheduling of $\FCC$}.  \oprocend
\end{definition}

Next, we examine the relation between the evolutions and the time and
communication complexities of a control and communication law, and of those
of its reschedulings.

\begin{proposition}
  With the same assumptions as in Definition~\ref{dfn:rescheduling},
  let $t \mapsto x(t)$ and $t \mapsto \tilde{x}(t)$ denote the network
  evolutions starting from $x_0 \in \prod_{i\in{I}}\supind{X}{i}_0$ under
  $\FCC$ and $\reschedule{\FCC}{s}{\Pc_I}$, respectively, and let
  $\map{\task}{\prod_{i\in{I}}\supind{X}{i}}{\myboolean}$ be a coordination
  task for $\network$.
  \begin{enumerate}
  \item For all $k\in\naturalzero$,
    \begin{equation}
      \label{eq:relation-trajectories}
      \supind{\tilde{x}}{i}(t)
      = 
      \begin{cases}
        \supind{\tilde{x}}{i} (t_{F(k-1)+1}) \,, \quad & \text{for all} \;
        t \in \bigcup_{\ell=F(k-1)+1}^{F(k)-1} [t_{\ell},t_{\ell+1}] \, ,
        \\
        \supind{x}{i} (h_{F(k)}(t)) \, , \quad & \text{for all} \; t \in
        [t_{F(k)},t_{F(k)+1}] \, .
      \end{cases}
    \end{equation}

  \item For all $x_0\in\prod_{i\in{I}}\supind{X}{i}_0$,
    \begin{equation*}
      \TC (\reschedule{\FCC}{s}{\Pc_I},\task,x_0) = s \cdot
      \TC(\FCC,\task,x_0) \, .      
    \end{equation*}
  \item If $\CommCost$ is additive, then, for all
    $x_0\in\prod_{i\in{I}}\supind{X}{i}_0$
    \begin{equation*}
      \MCC (\reschedule{\FCC}{s}{\Pc_I},\task,x_0) = \frac{1}{s} \cdot
      \MCC(\FCC,\task,x_0) \,  , 
    \end{equation*}
    and, therefore, the total communication cost of $\FCC$ is invariant
    under rescheduling.  \oprocend
  \end{enumerate}
\end{proposition}
\begin{proof}
  The relationships~\eqref{eq:relation-trajectories} are direct
  consequences of the definition of rescheduling.  We leave the bookkeeping
  to the interested reader.  By definition of $\TC(\FCC,\task,x_0)$, we
  have that $\task(x(t_k)) = \true$, for all $k \ge \TC(\FCC,\task,x_0)$,
  and $ \task(x(t_{\TC(\FCC,\task,x_0)-1})) = \false$. Let us rewrite these
  equalities in terms of the trajectories corresponding to the rescheduled
  control and communication law. From
  equation~\eqref{eq:relation-trajectories}, one can write
  $\supind{x}{i}(t_k) = \supind{x}{i}(h_{F(k)}(t_{F(k)})) =
  \supind{\tilde{x}}{i} (t_{F(k)})$, for all $i\in{I}$ and
  $k\in\naturalzero$.  Therefore, we have
  \begin{align*}
    & \task(\tilde{x} (t_{F(k)})) = \task(x(t_k)) = \true \, , \qquad
    \text{for all} \; F(k) \ge F( \TC(\FCC,\task,x_0)) \, ,
    \\
    &\task(\tilde{x} (t_{F( \TC(\FCC,\task,x_0) -1)})) = \task(x(t_{
      \TC(\FCC,\task,x_0)-1})) = \false \, ,
  \end{align*}
  where we have used the definition~\eqref{eq:rescheduled-message} of the
  rescheduled message-generation function. Now, note that by
  equation~\eqref{eq:relation-trajectories}, one has
  \begin{align*}
    \supind{\tilde{x}}{i} (t_{\ell}) = \supind{\tilde{x}}{i} (t_{F(\fl{\ell/s}-1)+1})
    \, , \quad \text{for all} \; \ell\in\naturalzero \; \text{and all} \;
    i \in I \, .
  \end{align*}
  Therefore, $\task(\tilde{x} (t_{F(\TC(\FCC,\task,x_0)-1)+1 })) =
  \task(\tilde{x} (t_{ F( \TC(\FCC,\task,x_0))}))$ and we can
  rewrite the previous identities as
  \begin{align*}
    \task(\tilde{x} (t_{k})) & = \true \, , \quad \text{for all} \;
    k \ge F(\TC(\FCC,\task,x_0)-1)+1 \, ,\\
    \task(\tilde{x} (t_{F( \TC(\FCC,\task,x_0) -1)})) &= \false \, ,
  \end{align*}
  which imply that
  \begin{align*}
    \TC (\reschedule{\FCC}{s}{\Pc_I},\task,x_0) & =
    F(\TC(\FCC,\task,x_0)-1)+1 = s \TC(\FCC,\task,x_0) \, .
  \end{align*}  
  As for the mean communication complexity, additivity of $\CommCost$
  implies
  \begin{align*}
    \newC(t_\ell,x(t_\ell)) =
    \newC(t_{F(\ell)-s+1},\tilde{x}(t_{F(\ell)-s+1})) + \dots + \newC
    (t_{F(\ell)},\tilde{x}(t_{F(\ell)})) \, ,
  \end{align*}
  where we have used $F(\ell-1)+1 = F(\ell)-s+1$.  Now, we compute
  \begin{multline*}
    \sum_{\ell=0}^{\TC ( \reschedule{\FCC}{s}{\Pc_I},\task,x_0) -1}
    \newC(t_\ell,\tilde{x}(t_\ell)) = \sum_{\ell=0}^{F(\TC (\FCC,\task,x_0)
      -1)} \newC(t_\ell,\tilde{x}(t_\ell))  \\
    = \sum_{\ell=0}^{\TC (\FCC,\task,x_0) -1} \sum_{k =
      F(\ell)-s+1}^{F(\ell)} \newC(t_k,\tilde{x}(t_k)) = \sum_{\ell =
      0}^{\TC (\FCC,\task,x_0) -1} \newC(t_\ell,x(t_\ell)) \, ,
  \end{multline*}
  which completes the proof of part~(iii).
\end{proof}

\begin{remark}
  It is worth noting that the notion of communication complexity defined
  in~\eqref{eq:cc-Klavins} is not invariant under rescheduling. Indeed,
  reasoning as before, one computes
  \begin{align*}
    \lim_{k\to+\infty} \frac{1}{k} \sum_{\ell=0}^{k}
    \newC(t_\ell,\tilde{x}(t_\ell)) &= \lim_{\tilde{k}\to+\infty}
    \frac{1}{\tilde{k}s-1} \sum_{\ell=0}^{\tilde{k}s-1}
    \newC(t_\ell,\tilde{x}(t_\ell))
    \\
    &=\lim_{\tilde{k}\to+\infty} \frac{1}{\tilde{k}s-1}
    \sum_{\tilde{\ell}=0}^{\tilde{k}-1} \sum_{\ell =
      F(\tilde{\ell}-1)+1}^{F(\tilde{\ell})} \newC(t_\ell,\tilde{x}(t_\ell)) \\
    &= \lim_{\tilde{k}\to+\infty} \frac{1}{\tilde{k}s-1}
    \sum_{\tilde{\ell}=0}^{\tilde{k}-1} \newC(t_\ell,x(t_{\tilde{\ell}})) \, .
  \end{align*}
  Therefore, $\cc ( \reschedule{\FCC}{s}{\Pc_I},x_0) = \frac{1}{s} \cc
  (\FCC,x_0)$. This means that, by performing a rescheduling of the control
  and communication law, one can indeed lower the measure of communication
  complexity $\cc$, although the trajectory described by the robotic
  network will continue to be the same. \oprocend
\end{remark}

\section{Motion coordination algorithms and their time complexity}
\label{se:algorithms}
In this section we provide examples of motion coordination algorithms for
robotic networks performing a variety of distributed tasks.  For each
algorithm and task, we present some results on the corresponding time and
communication complexity.

\subsection{Agreement on direction of motion and equidistance}

From Examples~\ref{ex:network-circle-Er},~\ref{ex:agree-pursuit}
and~\ref{ex:roundabout-agreement}, recall the definition of uniform network
$\network[\circledisk]$ of locally-connected first-order agents in
$\sph^1$, the \agreepursuitName control and communication law
$\FCC[\agreepursuit]$, and the two coordination tasks $\task[\dm]$ and
$\task[$\eps$-\equidistance]$.
\begin{theorem}
  For $\kprop\in]0,\frac{1}{2}[$, $r\in]0,2\pi]$, $\alpha = Nr-2\pi$
  and $\eps\in]0,1[$, the network $\network[\circledisk]$, the law
  $\FCC[\agreepursuit]$, and the tasks $\task[\dm]$ and
  $\task[$\eps$-\equidistance]$ together satisfy:
  \begin{enumerate}
  \item the upper bound $\TC(\task[\dm],\FCC[\agreepursuit]) \in
    O(Nr^{-1})$ and the lower bound
    \begin{equation*}
      \TC(\task[\dm],\FCC[\agreepursuit]) \in 
      \begin{cases}
        \Omega(r^{-1}) & \text{if} \; \alpha \ge 0 ,\\
        \Omega(N) & \text{if} \; \alpha \le 0 ;
      \end{cases}
    \end{equation*}
  \item if $\alpha>0$, then the upper bound
    $\TC(\task[$\eps$-\equidistance],\FCC[\agreepursuit]) \in
    O(N^2\log(N\eps^{-1})+N\log\alpha^{-1})$ and the lower bound
    $\TC(\task[$\eps$-\equidistance],\FCC[\agreepursuit]) \in \Omega
    (N^2 \log(\eps^{-1}))$. If $\alpha \le 0$, then
    $\FCC[\agreepursuit]$ does not achieve
    $\task[$\eps$-\equidistance]$ in general.
  \end{enumerate}
\end{theorem}

\begin{proof}
  Let us start by proving fact~(i).  Without loss of generality, assume
  $\dm^{[N]}(0)=\csym$, and that $\task[\dm]$ is \false at time $0$.
  Therefore, at least one agent is moving counterclockwise at time $0$, and
  we can define $k = \max\setdef{i\in I}{\dm^{[i]}(0)=\ccsym}$.  Define
  $t_k = \inf ( \setdef{\ell\in\naturalzero}{\dm^{[k]}(\ell)=\csym}\union
  \{+\infty\} )$.  In what follows we provide an upper bound on $t_k$.
  
  For $\ell<t_k$, define
  \begin{equation*}
    j(\ell) = \argmin \setdef{ 
      \distC(\theta^{[i]}(\ell)   ,\theta^{[k]}(\ell) ) }
    { \prior^{[i]}= N, \, i\in I  }.
  \end{equation*}
  In other words, for all instants of time when agent $k$ is moving
  counterclockwise, the agent $j(l)$ has $\prior$ equal to $N$, is moving
  clockwise, and is the agent closest to agent $k$ with these two
  properties.  Clearly,
  \begin{equation*}
    2 \pi > \distC(\theta^{[N]}(0),\theta^{[k]}(0))
    = \distC(\theta^{[j(0)]}(0),\theta^{[k]}(0)).
  \end{equation*}
  Additionally, for $\ell<t_k-1$, we claim that
  \begin{equation*}
    \distC( \theta^{[j(\ell)]}(\ell)  , \theta^{[j(\ell+1)]}(\ell+1) )
    > \kprop r.
  \end{equation*}
  This happens because either (1) there is no agent clockwise-ahead of
  $\theta^{[j(\ell)]}(\ell)$ within clockwise distance $r$ and, therefore,
  the claim is obvious, or (2) there are such agents. In case (2), let $m$
  denote the agent whose clockwise distance to agent $j(\ell)$ is maximal
  within the set of agents with clockwise distance $r$ from
  $\theta^{[j(\ell)]}(\ell)$.  Then,
  \begin{align*}
    \distC( \theta^{[j(\ell)]}(\ell) , \theta^{[j(\ell+1)]}(\ell+1) )
    &= \distC( \theta^{[j(\ell)]}(\ell) , \theta^{[m]}(\ell+1) )
    \\
    &= \distC( \theta^{[j(\ell)]}(\ell) , \theta^{[m]}(\ell) ) +
    \distC( \theta^{[m]}(\ell) , \theta^{[m]}(\ell+1) )
    \\
    &\geq \distC( \theta^{[j(\ell)]}(\ell) , \theta^{[m]}(\ell) ) +
    \kprop \big( r - \distC( \theta^{[j(\ell)]}(\ell) ,
    \theta^{[m]}(\ell) ) \big)
    \\
    &= \kprop r + ( 1-\kprop ) \distC( \theta^{[j(\ell)]}(\ell) ,
    \theta^{[m]}(\ell) ) \; \geq \kprop r,
  \end{align*}
  where the first inequality follows from the fact that at time $\ell$
  there can be no agent whose clockwise distance to agent $m$ is less than
  $( r - \distC( \theta^{[j(\ell)]}(\ell) , \theta^{[m]}(\ell) ) )$.

  In summary, either agent $k$ changes direction of motion or at each
  instant of time its distance to the closest agent with $\prior$ equal to
  $N$ decreases by a constant $\kprop r$.  This shows that
  \begin{equation*}
    t_k \in O(r^{-1}).
  \end{equation*}
  Because at least one more agent moves in the clockwise direction
  after $O(r^{-1})$ time, it follows that all agents will move
  clockwise after $O(Nr^{-1})$ time. This completes the proof of the
  upper bound in~(i).
  Let us prove the lower bound in~(i). Note that the number of agents
  that can fit into the circle $\sph^1$, spaced a distance $r$ apart
  one from each other, is at most $\fl{\frac{2 \pi}{r}}$.  Consider an
  initial configuration where $\dm^{[i]}(0)=\ccsym$ for $i \in
  \until{N-1}$, $\dm^{[N]}(0)=\csym$, and
  \begin{enumerate}
  \item for $\alpha>0$, $\fl{\frac{2 \pi}{r}}$ agents (including the
    agent $N$) lie at a distance $r$ one from each other in $\sph^1$,
    and the remaining $N- \fl{\frac{2 \pi}{r}}$ agents are located
    within a clockwise distance $\frac{r}{2}$ of the agent $i^*$ lying
    in the position in $\sph^1$ symmetric to the location of agent
    $N$;
  \item for $\alpha \le 0$, all agents lie at a distance $r$ in
    $\sph^1$ in counterclockwise order incrementally according to
    their unique identifier (note that they might not cover the whole
    circle).
  \end{enumerate}
  Note that the displacement of each agent is upper bounded by $\kprop
  r \le \frac{r}{2}$.  When $\alpha \ge 0$, the number of time steps
  that takes agent $i^*$ to receive the message with priority $N$ is
  lower bounded by $\frac{1}{2} \fl{\frac{2\pi}{r}}$. Therefore,
  $\TC(\task[\dm],\FCC[\agreepursuit]) \in \Omega(r^{-1})$.  When
  $\alpha \le 0$, the analysis of the network evolution is parallel to
  the one discussed in~\cite[Chapter 1]{NAL:97} for the leader
  election algorithm in static networks with the ring topology.  The
  number of time steps that takes agent $1$ to receive the message
  with priority $N$ is given by $N-1$. Therefore,
  $\TC(\task[\dm],\FCC[\agreepursuit]) \in \Omega(N)$.

  To prove fact~(ii), we assume that $\task[\dm]$ has been achieved (so
  that all agents are moving clockwise), and we first prove a fact
  regarding connectivity.  At time $\ell\in\naturalzero$, define
  \begin{equation*}
    \Hset(\ell) = %
    \setdef{x\in\sph^1}{
      \min_{i\in I}
      \distC(x,\theta^{[i]}(\ell))
      \,+\, \min_{j\in I} \distCC(x,\theta^{[j]}(\ell)) >r}.
  \end{equation*} 
  In other words, any point in $\Hset(\ell)$ is at least a distance $r$,
  clockwise or counterclockwise, from an agent.  Therefore, $\Hset(\ell)$
  does not contain any point between two agents separated by a distance
  less than $r$, and each connected component of $\Hset(\ell)$ has length
  at least $r$.  Let $\nHset(\ell)$ be the number of connected components
  of $\Hset(\ell)$, if $\Hset(\ell)$ is empty, then we take the convention
  that $\nHset(\ell)=0$.  Clearly, $\nHset(\ell)\leq N$.  We claim that, if
  $\nHset(\ell)>0$, then $t\mapsto \nHset(\ell+t)$ is non-increasing.  Let
  $d(\ell)<r$ be the distance between any two consecutive agents at time
  $\ell$. Because both agents move in the same direction, a simple
  calculation shows that
   \begin{equation*}
     d(\ell+1) \leq d(\ell) + \kprop (r - d(\ell)) 
     = (1-\kprop)d(\ell) + \kprop r
     < (1-\kprop)r + \kprop r = r.
   \end{equation*}
   This means that the two agents remain within distance $r$ and, therefore
   connected, at the following time instant. Because the number of
   connected components of $E_r(\theta^{[1]},\dots,\theta^{[N]})$ does not
   increase, it follows that the number of connected components of $\Hset$
   cannot increase.
       
   Next we claim that, if $\nHset(\ell)>0$, then there exists $t>\ell$ such
   that $\nHset(t)<\nHset(\ell)$. By contradiction, assume
   $\nHset(\ell)=\nHset(t)$ for all $t>\ell$.  Without loss of generality,
   let $\{1,\dots,m\}$ be a set of agents with the properties that
   $\distCC\big(\theta^{[i]}(\ell),\theta^{[i+1]}(\ell)\big)\leq r$, for
   $i\in\until{m}$, that $\theta^{[1]}(\ell)$ and $\theta^{[m]}(\ell)$
   belong to the boundary of $\Hset(\ell)$, and that there is no other set
   with the same properties and more agents.
   %% In other words, the agents $1,\dots,m$ are isolated.  
   One can show that, for $t\geq \ell$,
  \begin{align*}
    \theta^{[1]}(t+1) &= \theta^{[1]}(t)-\kprop r , \\
    \theta^{[i]}(t+1) &= \theta^{[i]}(t) - \kprop
    \distC(\theta^{[i]}(t),\theta^{[i-1]}(t)), \quad i\in\{2,\dots,m\}.
  \end{align*}
  If we define $d(t)= \big(\distCC(\theta^{[1]}(t),\theta^{[2]}(t))
  ,\dots, \distCC(\theta^{[m-1]}(t),\theta^{[m]}(t)) \big) \in
  \real_+^{m-1}$, then one can show that
  \begin{equation*}
    d(t+1) = \tridmat_{m-1}(\kprop,1-\kprop,0) \, d(t) 
    -  r 
    \begin{bmatrix}
      \kprop \\ 0 \\ \vdots \\ 0
    \end{bmatrix},
  \end{equation*}
  where the linear map
  $(a,b,c)\mapsto\tridmat_{m-1}(a,b,c)\in\real^{(m-1)\times(m-1)}$ is
  defined in Appendix~\ref{app:toeplitz}.  This is a discrete-time linear
  time-invariant dynamical system with unique equilibrium point
  $r(1,\dots,1)$.  By Theorem~\ref{thm:convergence-trid-circ}(ii) in
  Appendix~\ref{app:toeplitz}, for $\eta\in]0,1[$, the solution $t\mapsto
  d(t)$ to this system reaches a ball of radius $\eta$ centered at the
  equilibrium point in time $O(m\log{m} +\log\eta^{-1})$. (Here we used the
  fact that the initial condition of this system is bounded.) In turn, this
  implies that $t\mapsto \sum_{i=1}^m d_i(t)$ is larger than
  $(m-1)(r-\eta)$ in time $O(m \log m +\log\eta^{-1})$.
  
  We are now ready to find the contradiction and show that $\nHset(t)$
  cannot remain equal to $\nHset(\ell)$ for all time $t$. After time $O(m
  \log m +\log\eta^{-1}) = O(N \log N +\log\eta^{-1})$, we have:
  \begin{equation*}
    2\pi \geq \nHset(\ell) r + \sum_{j=1}^{\nHset(\ell)} (r-\eta) (m_j-1)
    =
    \nHset(\ell) r + 
    (N- \nHset(\ell)) (r-\eta)
    =  \nHset(\ell) \eta + N(r-\eta) .
  \end{equation*}
  Here $m_1,\dots,m_{\nHset(\ell)}$ are the number of agents in each
  isolated group, and each connected component of $\Hset(\ell)$ has length
  at least $r$. Now, take $\eta=\frac{Nr-2\pi}{N}=\frac{\alpha}{N}$, and
  the contradiction follows from
  \begin{gather*}
    2\pi \geq 
    \nHset(\ell) \eta + N r-N \eta 
    = \nHset(\ell)\eta + N r+ 2\pi - Nr
    = \nHset(\ell)\eta + 2\pi.
  \end{gather*}
  
  In summary this shows that, in time $O(N\log N+\log\eta^{-1}) = O(N\log N
  +\log\alpha^{-1})$, the number of connected components of $\Hset$ will
  decrease by one.  Therefore, in time $O(N^2\log N+N\log\alpha^{-1})$ the
  set $\Hset$ will become empty. At that time, the resulting network will
  obey the discrete-time linear time-invariant dynamical system:
  \begin{equation*}
    d(t+1) = \circmat_{N}(\kprop,1-\kprop,0) \, d(t).
  \end{equation*}
  Here $d(t)= \big(\distCC(\theta^{[1]}(t),\theta^{[2]}(t)) ,\dots,
  \distCC(\theta^{[N]}(t),\theta^{[N+1]}(t)) \big) \in \real_+^{N}$,
  with the convention $\theta^{[N+1]}=\theta^{[1]}$.  By
  Theorem~\ref{thm:convergence-trid-circ}(iii) in
  Appendix~\ref{app:toeplitz}, the solution $t\mapsto d(t)$ reaches
  the desired configuration in time $O(N^2\log\eps^{-1})$ with an
  error whose $2$-norm, and therefore, its $\infty$-norm is of order
  $\eps$.  In summary, the desired configuration is achieved in time
  $O(N^2\log(N\eps^{-1})+N\log\alpha^{-1})$.

  For the lower bound, consider an initial configuration with the
  properties that (i) agents are counterclockwise-ordered according to
  their unique identifier, (ii) the set $\Hset$ is empty, and (iii) the
  inter-agent distances $d(0)=
  \big(\distCC(\theta^{[1]}(0),\theta^{[2]}(0)) ,\dots,
  \distCC(\theta^{[N]}(0), \theta^{[1]}(0)) \big)$ are given by
  \begin{equation*}
    d(0) =
    \begin{bmatrix}
      \frac{2 \pi}{N} \\
      \vdots \\
      \frac{2 \pi}{N}
    \end{bmatrix}
    + k (\mathbf{v}_N + \overline{\mathbf{v}}_N) ,
  \end{equation*}
  where $\mathbf{v}_N$ is the eigenvector of
  $\circmat_{N}(\kprop,1-\kprop,0)$ corresponding to the eigenvalue $1
  - \kprop +\kprop \cos\left(\frac{2\pi}{N}\right) - \kprop \sqrt{-1}
  \sin\left(\frac{2\pi}{N}\right)$ (see Appendix~\ref{app:toeplitz}),
  and $k >0$ is chosen sufficiently small so that $d(0) \in
  \real_+^{N}$.  By Theorem~\ref{thm:convergence-trid-circ}(iii) in
  Appendix~\ref{app:toeplitz}, the solution $t\mapsto d(t)$ reaches
  the desired configuration in time $\Theta (N^2\log\eps^{-1})$ with
  an error whose $2$-norm, and therefore, its $\infty$-norm is of
  order $\eps$.  This concludes the result.
\end{proof}

\subsection{Rendezvous without connectivity constraint}

From Examples~\ref{ex:network-reald-Er}, \ref{ex:rendezvous}
and~\ref{ex:rendezvous-tasks}, recall the definition of uniform network
$\network[\realdisk]$ of locally-connected first-order agents in $\real^d$,
the \VicsekName control and communication law $\FCC[\VicsekSub]$, and the
coordination task $\task[\rendezvous]$.

\begin{theorem}
  For $d=1$, the network $\network[\realdisk]$, the law $\FCC[\VicsekSub]$,
  and the task $\task[\rendezvous]$ together satisfy
  $\TC(\task[\rendezvous],\FCC[\VicsekSub]) \in O(N^5)$ and
  $\TC(\task[\rendezvous],\FCC[\VicsekSub]) \in \Omega(N)$.
\end{theorem}

\begin{proof}
  One can easily prove that, along the evolution of the network, the
  ordering of the agents is preserved, i.e., if $x^{[i]}(\ell) \le
  x^{[j]}(\ell)$, then $x^{[i]}(\ell+1) \le x^{[j]}(\ell+1)$. However,
  links between agents are not necessarily preserved (see e.g.
  Figure~\ref{fig:rendezvous-disk}). Indeed, connected components may
  split along the evolution. However, mergings are not possible.
  Consider two contiguous connected components $C_1$ and $C_2$, with
  $C_1$ to the left of $C_2$.  By definition, the rightmost agent of
  $C_1$ and the leftmost agent of $C_2$ are at a distance strictly
  bigger than $r$. Now, by executing the algorithm, they can only but
  increase that distance, since the rightmost agent of $C_1$ will move
  to the left, and the leftmost agent of $C_2$ will move to the right.
  Therefore, connected components do not merge.

  Consider first the case of an initial configuration of the network
  for which the communication graph remains connected throughout the
  evolution.  Without loss of generality, assume that the agents are
  ordered from left to right according to their identifier, that is,
  $x^{[1]}(0) = (x_{0})_1 \le \dots \le x^{[N]}(0) = (x_{0})_N$.  Let
  $\alpha \in \{3,\dots,N\}$ have the property that agents $\{2,\dots,
  \alpha-1\}$ are neighbors of agent $1$, and agent $\alpha$ is not.
  (If instead all agents are within an interval of length $r$, then
  rendezvous is achieved in $1$ time instant, and the statement in
  theorem is easily seen to be true.)  Note that we can assume that
  agents $\{2,\dots, \alpha-1\}$ are also neighbors of agent $\alpha$.
  If this is not the case, then those agents that are neighbors of
  agent $1$ and not of agent $\alpha$, rendezvous with agent $1$ at
  the next time instant.  At the time instant $\ell=1$, the new
  updated positions satisfy
  \begin{align*}
    x^{[1]}(1) = \frac{1}{\alpha-1} \sum_{k=1}^{\alpha-1} x^{[k]}(0) ,
    \quad \quad x^{[\gamma]}(1) \in \Big[ \frac{1}{\alpha}
    \sum_{k=1}^{\alpha} x^{[k]}(0) ,* \Big] , \; \gamma \in \{ 2,
    \dots, \alpha-1\} ,
  \end{align*}
  where $*$ denotes certain unimportant point.

  Now, we show that
  \begin{align}\label{eq:inequality}
    x^{[1]}(\alpha-1) - x^{[1]}(0) \ge \frac{r}{\alpha (\alpha-1)} .
  \end{align}
  Let us first show the inequality for $\alpha=3$. Note that the fact
  that the communication graph remains connected implies that agent
  $2$ is still a neighbor of agent $1$ at the time instant $\ell=1$.
  Therefore $ x^{[1]}(2) \ge \frac{1}{2} ( x^{[1]}(1) + x^{[2]}(1) )$,
  and from here we deduce
  \begin{align*}
    x^{[1]}(2) - x^{[1]}(0) & \ge \frac{1}{2} \big( x^{[2]}(1) -
    x^{[1]}(0) \big)\\
    & \ge \frac{1}{2} \big( \frac{1}{3} \big( x^{[1]}(0) +x^{[2]}(0)
    +x^{[3]}(0) \big) - x^{[1]}(0) \big) \ge \frac{1}{6}
    \big(x^{[3]}(0) - x^{[1]}(0) \big) \ge \frac{r}{6} .
  \end{align*}
  Let us now proceed by induction.  Assume that
  inequality~\eqref{eq:inequality} is valid for $\alpha-1$, and let us
  prove it for $\alpha$. Consider first the possibility when at the
  time instant $\ell=1$, the agent $\alpha-1$ is still a neighbor of
  agent $1$.  In this case, $ x^{[1]}(2) \ge \frac{1}{\alpha-1}
  \sum_{k=1}^{\alpha-1} x^{[k]}(1)$, and from here we deduce
  \begin{align*}
    x^{[1]}(2) - x^{[1]}(0) & \ge \frac{1}{\alpha-1} \Big(
    x^{[\alpha-1]}(1) - x^{[1]}(0) \Big) \ge \frac{1}{\alpha-1} \Big(
    \frac{1}{\alpha} \sum_{k=1}^{\alpha} x^{[k]}(0) - x^{[1]}(0)
    \Big) \\
    & \ge \frac{1}{\alpha (\alpha-1)} \Big( x^{[\alpha]}(0) -
    x^{[1]}(0) \Big) \ge \frac{r}{\alpha (\alpha-1)} ,
  \end{align*}
  which in particular implies~\eqref{eq:inequality}.  Consider then
  the case when agent $\alpha-1$ is not a neighbor of agent $1$ at the
  time instant $\ell=1$. Let $\beta <\alpha$ such that agent $\beta-1$
  is a neighbor of agent $1$ at $\ell=1$, but agent $\beta$ is not.
  Since $\beta < \alpha$, we have by induction $ x^{[1]}(\beta) -
  x^{[1]}(1) \ge \frac{r}{\beta (\beta-1)}$. From here, we deduce that
  $x^{[1]}(\alpha-1) - x^{[1]}(0) \ge \frac{r}{\alpha (\alpha-1)}$.
  
  Inequality~\eqref{eq:inequality} implies that, at most in
  $\alpha-1\le N-1$ time instants, the leftmost agent traverses a
  distance greater than or equal to $\frac{r}{N (N-1)}$ (provided that
  at each step there exists at least another agent which is not its
  neighbor). Since $\diam (x_0,I) \le (N-1) r$, we deduce that in $N
  (N-1)^3$ time instants there cannot be any agent which is not a
  neighbor of the agent $1$. Hence, all agents rendezvous at the next
  time instant. Consequently,
  \begin{align*}
    \TC (\task[\rendezvous], \FCC[\VicsekSub],x_0) \le N (N-1)^3 + 1.
  \end{align*}
  Finally, for a general initial configuration $x_0$, because there
  are a finite number of agents, only a finite number of splittings
  (at most $N-1$) of the connected components of the communication
  graph can take place along the evolution.  Therefore, we conclude
  $\TC (\task[\rendezvous], \FCC[\VicsekSub]) = O( N^5)$.
  
  Let us now prove the lower bound. Consider an initial configuration $x_0
  \in \real^N$ where all agents are positioned in increasing order
  according to their identity, and exactly at a distance $r$ apart, say
  $(x_0)_{i+1} - (x_0)_i = r$, $i \in \until{N-1}$. Assume for simplicity
  that $N$ is odd - when $N$ is even, one can reason in an analogous way.
  Because of the symmetry of the initial condition, in the first time step,
  only agents $1$ and $N$ move. All the remaining agents remain in their
  position because it coincides with the average of its neighbors' position
  and its own. At the second time step, only agents $1$, $2$, $N-1$ and $N$
  move, and the others remain still because of the symmetry. Applying this
  idea iteratively, one deduces the time step when agents $\frac{N-1}{2}$
  and $\frac{N+3}{2}$ move for the first time is lower bounded by
  $\frac{N-1}{2}$. Since both agents have still at least a neighbor (agent
  $\frac{N+1}{2}$), the task $\task[\rendezvous]$ has not been achieved yet
  at this time step. Therefore, $\TC (\task[\rendezvous],
  \FCC[\VicsekSub],x_0) \ge \frac{N-1}{2}$, and the result follows.
\end{proof}

\subsection{Rendezvous with connectivity constraint}
\label{se:rendezvous}

In this section we shall consider both networks $\network[\realdisk]$ and
$\network[\realLD]$ presented in Examples~\ref{ex:network-reald-Er}
and~\ref{ex:network-reald-ELDr}.

\subsubsection*{Circumcenter control and communication law} 
Here we define the \emph{circumcenter} control and communication law
$\FCC[\circumcenter]$ for both networks $\network[\realdisk]$ and
$\network[\realLD]$.  This is a uniform, static, time-independent law
originally introduced by~\cite{HA-YO-IS-MY:99} and later studied
in~\cite{JL-ASM-BDOA:04b,JC-SM-FB:04h-tmp}.  Loosely speaking, the evolution
of the network under the circumcenter control and communication law can be
described as follows:
\begin{quote}
  \emph{[Informal description]} Communication rounds take place at each
  natural instant of time.  At each communication round each agent performs
  the following tasks: (i) it transmits its position and receives its
  neighbors' positions; (ii) it computes the circumcenter of the point set
  comprised of its neighbors and of itself, and (iii) it moves toward this
  circumcenter while maintaining connectivity with its neighbors.
\end{quote}

Let us present this description in more formal terms. We set
$\timeschedule=\naturalzero$, $L=\real^d$, and
$\supind{\msg}{i}=\msgstandard$, $i \in I$.  In order to define the control
function, we need to introduce some preliminary constructions.  First,
connectivity is maintained by restricting the allowable motion of each
agent in the following appropriate manner.  If agents $i$ and $j$ are
neighbors at time $\ell\in\naturalzero$, then we require their subsequent
positions to belong to
$\cball{\frac{r}{2}}{\frac{\supind{x}{i}(\ell)+\supind{x}{j}(\ell)}{2}}$.
If an agent $i$ has its neighbors at locations $\{q_1,\dots,q_l\}$ at time
$\ell$, then its \emph{constraint set}
$\Dc_{\supind{x}{i}(\ell),r}(\{q_1,\dots,q_l\})$ is
\begin{equation*}
  \Dc_{\supind{x}{i}(\ell),r}(\{q_1,\dots,q_l\})
  = \bigcap_{q\in \{q_1,\dots,q_l\}}
  \Bigcball{\frac{r}{2}}{\frac{\supind{x}{i}(\ell)+q}{2}}. 
\end{equation*}
Second, in order to maximize the displacement toward the circumcenter of
the point set comprised of its neighbors and of itself, each agent solves a
convex optimization problem that can be stated in general as follows.  For
$q_0$ and $q_1$ in $\real^d$, and for a convex closed set $Q\subset\real^d$
with $q_0\in Q$, let $\lambda (q_0,q_1,Q)$ denote the solution to the
strictly convex problem:
\begin{equation*}
  \begin{split}
    \maximize &\enspace \lambda \\
    \subj
    &\enspace \lambda \leq 1,\; (1-\lambda) q_0+\lambda q_1 \in Q.
  \end{split}
\end{equation*}
Under the stated assumptions the solution exists and is unique.  Third,
note that since the agents operate with the standard message-generation
function, it is clear that the natural projection $\piSpace{\real^d}$ maps
the messages $\supind{y}{i}(\ell)$ received at time $\ell\in\naturalzero$
by the agent $i \in I$ onto the positions of its neighbors.  We are now
ready to define the last constitutive element of $\FCC[\circumcenter]$.
Define the control function $\map{\ctrl}{\real^d\times\real^d\times{L^N}}
{\real^d}$ by
\begin{equation}
  \label{eq:rendezvous-control-law}
  \ctrl(x,\sampled{x},y) = 
  \lambda_* \cdot  (\CircumC(\piSpace{\real^d} (y) \union \{\sampled{x}\})
  - \sampled{x} )\, ,
\end{equation}
where $\lambda_* = \lambda (\sampled{x}, (\CircumC(\piSpace{\real^d} (y) \union
\{\sampled{x}\}), \Dc_{\sampled{x},r}(\piSpace{\real^d} (y)))$. 
Evolving under this control law, it is clear that, at time $\fl{t}+1$, each
agent $i$ reaches the point $(1- \lambda_*) \supind{x}{i}(\fl{t}) +
\lambda_* \CircumC( \piSpace{\real^d} (\supind{y}{i}(\fl{t}) ) \union \{
\supind{x}{i}(\fl{t}) \} )$.

\bigskip

Next, we consider the network $\network[\rIdisk]$ of
locally-$\infty$-connected first-order agents in $\real^d$, see
Example~\ref{ex:network-reald-EIr}.  For this network we define the
\emph{parallel circumcenter law}, $\FCC[\parallelcircumcenter]$, by
designing $d$ decoupled circumcenter laws running in parallel on each
coordinate axis of $\real^d$.  As before, this law is uniform, static and
time-independent. As before, we set $\timeschedule=\naturalzero$,
$L=\real^d$, and $\supind{\msg}{i}=\msgstandard$, $i \in I$.  We define the
control function $\map{\ctrl}{\real^d\times\real^d\times{L^N}}{\real^d}$ by
\begin{equation}
  \label{eq:FCC-parallel-circumcenter-control-law}
  \ctrl(x,\sampled{x},y) 
  = \left(\, \CircumC(\tau_1(\MM)) - \tau_1(\sampled{x})
    ,\dots,\, 
    \CircumC(\tau_d(\MM)) - \tau_d(\sampled{x})
  \right) ,
\end{equation}
where $\MM = \pi_{\real^d}(y)\union\{\sampled{x}\}$, and
$\map{\tau_1,\dots,\tau_d}{\real^d}{\real}$ denote the canonical
projections of $\real^d$ onto~$\real$. See
Fig.~\ref{fig:rendezvous-parallel} for an illustration of this law in
$\real^2$.  
\begin{figure}[htb] 
  \centering
  \resizebox{.4\linewidth}{!}{\input{rendezvous-parallel-disk.tex}}
  \caption{Parallel circumcenter control and  communication law in
    $\real^2$. The target point for the agent $i$ is plotted in light gray
    and has coordinates $(\CircumC (\tau_1(\supind{\MM}{i})) , \CircumC
    (\tau_2(\supind{\MM}{i})) )$.}\label{fig:rendezvous-parallel}
\end{figure}

\subsubsection*{Asymptotic behavior and complexity analysis}
The following theorem summarizes the results known in the literature
about the asymptotic properties of the circumcenter control and
communication law.

\begin{theorem}[Correctness of the circumcenter law]
  \label{th:FCC-rendezvous-correctness}
  For $d\in\natural$, $r\in\real_+$ and $\eps\in\real_+$, the following
  statements hold:
  \begin{enumerate} 
  \item on the network $\network[\realdisk]$, the law $\FCC[\circumcenter]$
    achieves the exact rendezvous task $\task[\rendezvous]$;
    
  \item on the network $\network[\realLD]$, the law $\FCC[\circumcenter]$
    achieves the $\eps$-rendezvous task $\task[$\eps$-\rendezvous]$;
    
  \item on the network $\network[\realIdisk]$, the law
    $\FCC[\parallelcircumcenter]$ achieves the exact rendezvous task
    $\task[\rendezvous]$;
    
  \item the evolutions of $(\network[\realdisk],\FCC[\circumcenter])$, of
    $(\network[\realLD],\FCC[\circumcenter])$, and of
    $(\network[\realIdisk],\FCC[\parallelcircumcenter])$ have the property
    that, if two agents belong to the same connected component of the
    communication graph at some time $\ell\in\naturalzero$, then they
    continue to belong to the same connected component of the communication
    graph for all subsequent times $k\geq\ell$. \oprocend
  \end{enumerate}
\end{theorem}
\begin{proof}
  The results on $\network[\rdisk]$ appeared originally
  in~\cite{HA-YO-IS-MY:99}.  The proof for the results on $\network[\rLD]$
  is provided in~\cite{JC-SM-FB:04h-tmp}. We postpone the proof for
  $\network[\rIdisk]$ to the proof of
  Theorem~\ref{th:rendezvous-time-complexity} below.
\end{proof}

Next we analyze the time complexity of $\FCC[\circumcenter]$.  We provide
complete results only for the case $d=1$.  As we see next, the complexity
properties of $\FCC[\circumcenter]$ differ dramatically when applied to the
two robotic networks with different communication graphs.

\begin{theorem}[Time complexity of circumcenter law]
  \label{th:rendezvous-time-complexity}
  For $r\in\real_+$ and $\eps\in]0,1[$, the following statements hold:
  \begin{enumerate}
  \item for $d=1$, on the network $\network[\realdiskone]$,
    $\TC(\task[\rendezvous],\FCC[\circumcenter]) \in \Theta (N)$;
  \item for $d=1$, on the network $\network[\realLDone]$,
    $\TC(\task[$(r\eps)$-\rendezvous],\FCC[\circumcenter]) \in \Theta (N^2
    \log (N \eps^{-1}))$;
  \item for $d\in\natural$, on the network $\network[\realIdisk]$,
    $\TC(\task[\rendezvous],\FCC[\parallelcircumcenter]) \in \Theta (N)$.
    \oprocend
  \end{enumerate}
\end{theorem} 

\begin{proof}
  Let $x_0\in\real^N$.  Throughout the proof, neighboring relationships are
  understood with respect to the $r$-disk graph.  First of all, let us show
  that, for $n=1$, the connectivity constraints on each agent $i \in I$
  imposed by the constraint set $\Dc_{\supind{x}{i},r}(\piSpace{\real}
  (y)))$ are superfluous, i.e., the solution of the convex optimization
  problem is $\lambda_* = 1$ (cf.
  equation~\eqref{eq:rendezvous-control-law}). To see this, assume that
  agents $i$ and $j$ are neighbors at time instant $\ell$, define
  $\supind{\MM}{i}$ as $\piSpace{\real^d} (\supind{y}{i}(\ell)) \union
  \{\supind{x}{i}(\ell) \}$, and let us show that $\CircumC (\supind{\MM}{i})$
  belongs to $\cball{\frac{r}{2}}{\frac{\supind{x}{i}(\ell) +
      \supind{x}{j}(\ell)}{2}}$.  Without loss of generality, let
  $\supind{x}{i}(\ell) \le \supind{x}{j}(\ell)$.  Let $x_{-}^{[i]}(\ell)$,
  $x_{+}^{[i]}(\ell)$ denote the positions of the leftmost and rightmost
  agents among the neighbors of agent~$i$.  Note that $\supind{x}{i}(\ell)
  \le \supind{x}{j}(\ell) \le x_+^{[i]}(\ell) $ and $\CircumC (\supind{\MM}{i})
  = \frac{1}{2} (x_{-}^{[i]}(\ell) + x_{+}^{[i]}(\ell))$.  Then,
  \begin{align*}
    \big| \CircumC (\supind{\MM}{i}) - \frac{1}{2}(\supind{x}{i}(\ell) + \supind{x}{j}(\ell))
    \big| & = \frac{1}{2} \big| x_{-}^{[i]}(\ell) - \supind{x}{i}(\ell) +
    x_{+}^{[i]}(\ell) - \supind{x}{j}(\ell) \big| \\
    & \le \frac{1}{2} \max \{ | x_{-}^{[i]}(\ell) - \supind{x}{i}(\ell) | ,
    | x_{+}^{[i]}(\ell) - \supind{x}{j}(\ell) |\} \le \frac{r}{2} \, ,
  \end{align*}
  as claimed. Therefore, we have that $\supind{x}{i}(\ell+1) = \CircumC
  (\supind{\MM}{i})$. Likewise, one can deduce $\CircumC (\supind{\MM}{i}) \le
  \CircumC (\supind{\MM}{j})$, and therefore, the order of the agents is
  preserved.
 
  Consider first the case when $\subscr{E}{\rdisk}(x_0)$ is connected.
  Without loss of generality, assume that the agents are ordered from left
  to right according to their identifier, that is, $x^{[1]}(0) = (x_{0})_1
  \le \dots \le x^{[N]}(0) = (x_{0})_N$.  Let $\alpha \in \{3,\dots,N\}$
  have the property that agents $\{2,\dots, \alpha-1\}$ are neighbors of
  agent $1$, and agent $\alpha$ is not. (If instead all agents are within
  an interval of length $r$, then rendezvous is achieved in $1$ time
  instant, and the statement in theorem is easily seen to be true.)  See
  Fig.~\ref{fig:rendezvous-disk-line} for an illustration of these
  definitions.
   \begin{figure}[htb] 
     \centering
     \resizebox{.8\linewidth}{!}{\input{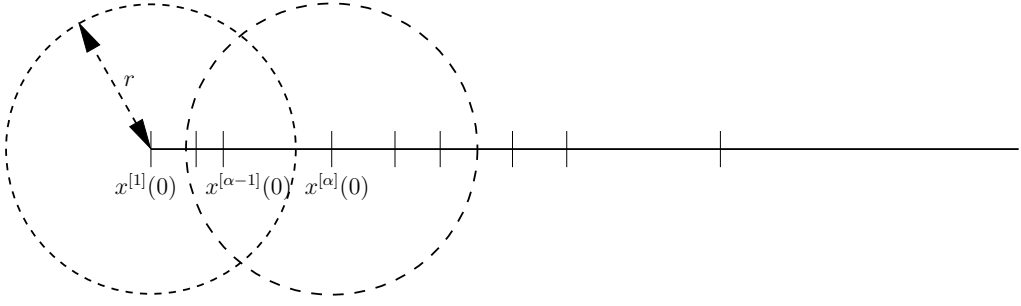}}
     \caption{Definition of $\alpha \in \{3,\dots,N\}$  for an initial network
       configuration.}\label{fig:rendezvous-disk-line}
   \end{figure}
   Note that we can assume that agents $\{2,\dots, \alpha-1\}$ are also
   neighbors of agent $\alpha$.  If this is not the case, then those agents
   that are neighbors of agent $1$ and not of agent $\alpha$, rendezvous
   with agent $1$ at the next time instant.  At the time instant $\ell=1$,
   the new updated positions satisfy
  \begin{align*}
    x^{[1]}(1) = \frac{x^{[1]}(0) + x^{[\alpha-1]}(0)}{2}\,, \quad
    x^{[\gamma]}(1) \in \left[\frac{x^{[1]}(0) +
        x^{[\alpha]}(0)}{2},\frac{x^{[1]}(0) + x^{[\gamma]}(0) +r}{ 2}
    \right] \!, \; \gamma \in \{ 2, \dots, \alpha-1\} \,.
  \end{align*}
  These equalities imply that $x^{[1]}(1) - x^{[1]}(0) = \frac{1}{2} \big(
  x^{[\alpha-1]}(0) - x^{[1]}(0) \big) \le \frac{1}{2} r$.  Analogously, we
  deduce $x^{[1]}(2) - x^{[1]}(1) \le \frac{1}{2} r$, and therefore
  \begin{align}\label{eq:upper-bound-step}
    x^{[1]}(2) - x^{[1]}(0) \le r \,.
  \end{align}
  On the other hand, from $x^{[1]}(2) \in \left[\frac{1}{2} \big(x^{[1]}(1)
    + x^{[\alpha-1]}(1) \big), * \right]$ (where the symbol $*$ represents
  a certain unimportant point in $\real$), we deduce that
  \begin{align}\label{eq:lower-bound-step}
    x^{[1]}(2) - x^{[1]}(0) & \ge \frac{1}{2} \big(x^{[1]}(1) +
    x^{[\alpha-1]}(1) \big) - x^{[1]}(0) \ge \frac{1}{2}
    \big(x^{[\alpha-1]}(1) - x^{[1]}(0) \big) \nonumber \\
    & \ge \frac{1}{2} \Big( \frac{x^{[1]}(0) + x^{[\alpha]}(0)}{2} -
    x^{[1]}(0) \Big) = \frac{1}{4} \big( x^{[\alpha]}(0) - x^{[1]}(0)
    \big) \ge \frac{1}{4} r \, .
  \end{align}
  Inequalities~\eqref{eq:upper-bound-step} and~\eqref{eq:lower-bound-step}
  mean that, after at most two time instants, agent $1$ has traveled an
  amount larger than $r/4$. In turn this implies that
  \begin{align*}
    \frac{\diam (x_0,I)}{r} \le%
    \TC(\task[\rendezvous],\FCC[\circumcenter],x_0)%
    \le \frac{4 \diam (x_0,I)}{r}\,.
  \end{align*}
  If $\subscr{E}{\rdisk} (x_0)$ is not connected, note that along the
  network evolution, the connected components of the $r$-disk graph do not
  change.  Therefore, using the previous characterization on the amount
  traveled by the leftmost agent of each connected component in at most two
  time instants, we deduce that
  \begin{align*}
    \frac{1}{r} \max_{C \in \mathcal{C}_{\subscr{E}{\rdisk}}(x_0)} \diam
    (x_0,C) \le \TC(\task[\rendezvous],\FCC[\circumcenter],x_0) \le
    \frac{4}{r} \max_{C \in \mathcal{C}_{\subscr{E}{\rdisk}}(x_0)} \diam
    (x_0,C) \, .
  \end{align*}
  Note that the connectedness of each $C \in
  \mathcal{C}_{\subscr{E}{\rdisk}}(x_0)$ implies that $\diam (x_0,C) \le
  (N-1) r$, and therefore $\TC (\task[\rendezvous],\FCC[\circumcenter]) \in
  O(N)$. Moreover, for $x_0 \in \real^N$ such that $(x_0)_{i+1} - (x_0)_i =
  r$, $i \in \until{N-1}$, we have $\diam (x_0,I) = (N-1)r$, and therefore
  $\TC (\task[\rendezvous], \FCC[\circumcenter] ,x_0) \ge N-1$. This
  concludes the proof of fact~(i):
  \begin{align*}
    \TC (\task[\rendezvous], \FCC[\circumcenter]) \in \Theta (N) \, .
  \end{align*}
  
  Next we prove fact~(ii).  In the $r$-limited Delaunay graph, two agents
  on the line that are at most at a distance $r$ from each other are
  neighbors if and only if there are no other agents between them. Also,
  note that the $r$-limited Delaunay graph and the $r$-disk graph have the
  same connected components (cf.~\cite{JC-SM-FB:03p-tmp}). Using an
  argument similar to the one above, one can show that the connectivity
  constraints imposed by the constraint sets set
  $\Dc_{\supind{x}{i}(\fl{t}),r}(\piSpace{\real} (y)))$ are again
  superfluous.
  
  Consider first the case when $\subscr{E}{\rLD} (x_0)$ is connected.  Note
  that this is equivalent to $\subscr{E}{\rdisk} (x_0)$ being connected.
  Without loss of generality, assume that the agents are ordered from left
  to right according to their identifier, that is, $x^{[1]}(0) = (x_{0})_1
  \le \dots \le x^{[N]}(0) = (x_{0})_N$.  The evolution of the network
  under $\FCC[\circumcenter]$ can then be described as the discrete-time
  dynamical system
  \begin{multline*}
    x^{[1]}(\ell + 1) = \frac{1}{2} (x^{[1]}(\ell)+ x^{[2]}(\ell)) \,
    , \quad x^{[2]}(\ell + 1) = \frac{1}{2} (x^{[1]}(\ell)+
    x^{[3]}(\ell)) \,, \;   \dots \, ,  \\
    , \; \dots \, , \; x^{[N-1]}(\ell + 1) = \frac{1}{2} (
    x^{[N-2]}(\ell) + x^{[N]}(\ell)) \,, \quad x^{[N]}(\ell + 1) =
    \frac{1}{2} (x^{[N-1]}(\ell)+ x^{[N]}(\ell) ) \,.
  \end{multline*}
  Note that this evolution respects the ordering of the agents.
  Equivalently, we can write $x(\ell+1) = A \, x(\ell)$, where $A$ is
  the $N\times N$ matrix given by
  \begin{equation*}
    A = 
    \begin{bmatrix}
      \frac{1}{2} & \frac{1}{2} & 0 & \dots & \dots & 0 \\
      \frac{1}{2} & 0 & \frac{1}{2} & \dots & \dots & 0 \\
      0 & \frac{1}{2} & 0 & \frac{1}{2} & \dots &0 \\
      \vdots & & \ddots & \ddots & \ddots & \vdots \\
      0 & \dots & \dots & \frac{1}{2} & 0 & \frac{1}{2} \\
      0 & \dots & \dots & 0 & \frac{1}{2} & \frac{1}{2}
    \end{bmatrix} \! .
  \end{equation*}
  Now, note that $A=\atridmat_N^{+}\big(\frac{1}{2},0\big)$ as defined in
  Appendix~\ref{app:toeplitz}.  Theorem~\ref{thm:convergence-atrid}(i)
  implies that, for $\xave=\frac{1}{N}\mathbf{1}^Tx(0)$, we have that
  $\lim_{\ell\to+\infty} x(\ell)=\xave\mathbf{1}$, and that the maximum
  time required for $\Enorm{x(\ell)-\xave\mathbf{1}\big}\leq
  \eta\Enorm{x(0)-\xave\mathbf{1}}$ (over all initial conditions
  $x(0)\in\real^N$) is $\Theta\big(N^2\log\eta^{-1}\big)$.  (As an aside,
  this also implies that the agents rendezvous at the location given by the
  average of their initial positions. In other words, we can forecast the
  asymptotic rendezvous position for this case, as opposed to the case with
  the $r$-disk communication graph.)
  
  Next, let us convert the contraction inequality on $2$-norms into an
  appropriate inequality on $\infty$-norms.  Note that $\diam(x_0,I) \le
  (N-1)r$ because $\subscr{E}{\rLD} (x_0)$ is connected. Therefore
  \begin{equation*}
    \Infnorm{x(0) - \subscr{x}{ave}  \mathbf{1}} 
    = \max_{i\in{I}} |\supind{x}{i}(0) -\subscr{x}{ave}| 
    \leq | x^{[1]}_0 - x^{[N]}_0| \leq (N-1) r.
  \end{equation*}
  For $\ell$ of order $N^2\log\eta^{-1}$, we use this bound on
  $\Infnorm{x(0) - \subscr{x}{ave} \mathbf{1}} $ and the basic inequalities
  $\norm{v}{\infty}\leq\norm{v}{2}\leq \sqrt{N} \norm{v}{\infty}$ for all
  $v\in\real^N$, to obtain: 
  \begin{align*}
    \Infnorm{x(\ell)-\subscr{x}{ave}\mathbf{1}}
    &\leq 
    \Enorm{x(\ell)-\subscr{x}{ave}\mathbf{1}}
    \leq \eta 
    \Enorm{x(0)-\subscr{x}{ave}\mathbf{1}}
    \leq \eta\sqrt{N}
    \Infnorm{x(0)-\subscr{x}{ave}\mathbf{1}}
    \leq \eta \sqrt{N}(N-1)r.
  \end{align*}
  This means that $(r\eps)$-rendezvous is achieved for $\eta
  \sqrt{N}(N-1)r=r\eps$, that is, in time
  $O(N^2\log\eta^{-1})=O(N^2\log(N\eps^{-1}))$.
 
  Next, we show the lower bound.  From equation~\eqref{eq:vN} in the proof
  of Theorem~\ref{thm:convergence-trid-circ}, we recall the unit-length
  vector $\mathbf{v}_{N-1}\in\real^{N-1}$.  This vector is an eigenvector
  of $\tridmat_{N-1}(\frac{1}{2},0,\frac{1}{2})$ corresponding to the
  largest singular value $\cos(\frac{\pi}{N})$.  For
  $\mu=\frac{-1}{10\sqrt{2}}rN^{5/2}$, we then define the initial
  condition  $x_0 = \mu P_+\begin{bmatrix}0\\
    \mathbf{v}_{N-1}\end{bmatrix}\in\real^N$.  One can show that
  $(x_0)_i<(x_0)_{i+1}$ for $i\in\until{N-1}$, that
  $\subscr{(x_0)}{ave}=0$, and that
  $\max\setdef{(x_0)_{i+1}-(x_0)_i}{i\in\until{N-1}}\leq r$.  By
  Lemma~\ref{lem:bounds-P} and because $\norm{w}{\infty}\leq\norm{w}{2}\leq
  \sqrt{N} \norm{w}{\infty}$ for all $w\in\real^N$, we compute
  \begin{equation*}
    \Infnorm{x_0} = \frac{rN^{5/2}}{10\sqrt{2}}
    \Bignorm{P_+ \begin{bmatrix}0\\  \mathbf{v}_{N-1}\end{bmatrix} }{\infty}
    \geq 
    \frac{rN^2}{10\sqrt{2}}
    \Bignorm{P_+ \begin{bmatrix}0\\  \mathbf{v}_{N-1}\end{bmatrix} }{2}
    \geq
    \frac{rN}{10\sqrt{2}}
    \Enorm{\mathbf{v}_{N-1}}  =
    \frac{rN}{10\sqrt{2}}. 
  \end{equation*}
  The trajectory $x(\ell)= (\cos(\frac{\pi}{N}))^\ell x_0$ therefore
  satisfies
  \begin{align*}
    \Infnorm{x(\ell)} & = \Big( \cos\Big(\frac{\pi}{N}\Big) \Big)^\ell
    \Infnorm{x_0} \geq \frac{rN}{10\sqrt{2}}
    \Big(\cos\Big(\frac{\pi}{N}\Big)\Big)^\ell.
  \end{align*}
  Therefore, $\Infnorm{x(\ell)}$ is larger than $\frac{1}{2}r\eps$ so long
  as $\frac{1}{10\sqrt{2}} N (\cos(\frac{\pi}{N}))^\ell > \frac{1}{2}\eps$,
  that is, so long as
  \begin{equation*}
    \ell < \frac{\log(\eps^{-1}N)
      -\log(5\sqrt{2})}{-\log\big(\cos(\frac{\pi}{N})\big)}.
  \end{equation*}
  The rest of the proof is the same as in
  Theorem~\ref{thm:convergence-trid-circ}(i) for the lower bound result.
   
  If $\subscr{E}{\rLD} (x_0)$ is not connected, note that along the network
  evolution, the connected components do not change.  Therefore, the
  previous reasoning can be applied to each connected component.  Since the
  number of agents in each connected component is strictly less that $N$,
  the time complexity can only but improve. Therefore, we conclude that
  \begin{align*}
    \TC(\task[\rendezvous],\FCC[\circumcenter]) 
    \in \Theta(N^2\log(N\eps^{-1})) \, .
  \end{align*}
  This completes the proof of fact~(ii).
  
  Finally, we prove the statements regarding $\network[\realIdisk]$ and
  $\FCC[\parallelcircumcenter]$ in fact~(iii) and in the previous
  Theorem~\ref{th:FCC-rendezvous-correctness}.  By definition, agents $i$
  and $j$ are neighbors at time $\ell\in\naturalzero$ if and only if
  $\Infnorm{\supind{x}{i}(\ell)-\supind{x}{j}(\ell)}\le r$, which is
  equivalent to
  \begin{equation*} 
    | \tau_k ( \supind{x}{i}(\ell)) - \tau_k (\supind{x}{j}(\ell))| \le r \,, \quad  k
      \in \until{d}\,. 
  \end{equation*}
  Recall from the proof of fact~(i) that the connectivity constraints of
  $\FCC[\circumcenter]$ on each agent are trivially satisfied in the
  $1$-dimensional case. This fact has the following important consequence:
  from the expression for the control function in
  $\FCC[\parallelcircumcenter]$, we deduce that the evolution under
  $\FCC[\parallelcircumcenter]$ of the robotic network
  $\network[\realIdisk]$ (in $d$ dimensions) can be alternatively described
  as the evolution under $\FCC[\circumcenter]$ of $d$ robotic networks
  $\network[\realdiskone]$ in $\real$.  The correctness and the time
  complexity results now follows from the analysis of $\FCC[\circumcenter]$
  at $d=1$.
\end{proof}

Finally we proceed to characterize the mean communication complexity of the
circumcenter control and communication law.  We consider the case of a
unidirectional communication model with one-round cost function depending
linearly on the cardinality of the communication graph.

\begin{theorem}[Mean communication complexity of circumcenter law for
  unidirectional communication]
  \label{th:MCC-rendezvous}
  For $r\in\real_+$ and $\eps\in]0,1[$, the following statements hold:
  \begin{enumerate}
  \item for $d=1$, on the network $\network[\realdisk]$,
    $\MCC(\task[\rendezvous], \FCC[\circumcenter]) \in \Theta(N^2)$;
  \item for $d=1$, on the network $\network[\realLD]$, $\MCC
    (\task[\rendezvous], \FCC[\circumcenter]) \in \Theta(N)$;
  \item for $d\in\natural$, on the network $\network[\realIdisk]$,
    $\MCC(\task[\rendezvous], \FCC[\parallelcircumcenter]) \in \Theta(N^2)$.
    \oprocend
  \end{enumerate}
\end{theorem}
\begin{proof}
  We start by proving fact~(i).  Let $x_0 \in (\real^d)^N$ be such that all
  $(x_0)_i$, $i \in I$, belong to a closed ball of radius $\sqrt{2}r/4$.
  In such a case, one can deduce that (i) $\subscr{E}{\rdisk}(x_0)$ is the
  complete graph, and therefore all agents compute the same goal point
  $\CircumC$, and (ii) this circumcenter belongs to the constrained set of each
  agent $i \in I$.  As a consequence, all mobile agents rendezvous at the
  same location $\CircumC$ in one time instant.  For such initial conditions, we
  have $\TC (\task[\rendezvous], \FCC[\circumcenter], x_0) = 1$, and
  therefore
  \begin{align*}
    \MCC(\task[\rendezvous], \FCC[\circumcenter], x_0) &= \frac{1}{\TC
      (\task[\rendezvous], \FCC[\circumcenter], x_0)} \sum_{\ell =0}^{\TC
      (\task[\rendezvous], \FCC[\circumcenter], x_0)-1}
    \CommCost[L]\circ\subscr{E}{\nonnllmsgs} (\ell,x(\ell)) \\
    &= N (N-1) \, .
  \end{align*}
  This proves that $\MCC (\task[\rendezvous], \FCC[\circumcenter]) \in
  \Theta(N^2)$, that is, fact~(i).  Additionally, fact~(iii) is proved by
  analyzing the parallel circumcenter law as $d$ decoupled versions of the
  circumcenter law.  Next we show fact~(ii).  Given that the $r$-limited
  Delaunay graph is a subgraph of the Delaunay
  graph~\cite{JC-SM-FB:03p-tmp}, and the number of edges of the latter is
  bounded by $3N-6$ (see, for instance,~\cite{AO-BB-KS-SNC:00}), we deduce
  that $\MCC (\task[\rendezvous], \FCC[\circumcenter]) \in O(N)$. On the
  other hand, for any initial configuration $x_0$ where all the agents are
  aligned and non-coincident, there are $N-1$ neighboring relationships
  which are preserved throughout the network evolution.  In this case,
  $\MCC (\task[\rendezvous], \FCC[\circumcenter], x_0) = 2 (N-1)$, and the
  result in fact~(ii) follows.
\end{proof}

\begin{remark}
  Theorems~\ref{th:rendezvous-time-complexity} and~\ref{th:MCC-rendezvous}
  induce lower bounds on the time and mean communication complexity of the
  circumcenter law for the higher-dimensional case. Indeed, as a
  consequence of these results, we have
  \begin{enumerate}
  \item for $d \in \natural$, on the network $\network[\realdiskone]$,
    $\TC(\task[\rendezvous],\FCC[\circumcenter]) \in \Omega (N)$ and
    $\MCC(\task[\rendezvous], \FCC[\circumcenter]) \in \Omega(N^2)$;
  \item for $d \in \natural$, on the network $\network[\realLDone]$,
    $\TC(\task[$(r\eps)$-\rendezvous],\FCC[\circumcenter]) \in \Omega (N^2
    \log (N \eps^{-1}))$ and $\MCC (\task[\rendezvous],
    \FCC[\circumcenter]) \in \Omega(N)$.
  \end{enumerate}
  We have performed extensive numerical simulations for the case $d=2$ and
  the network $\network[\realdisk]$.  We have ran the algorithm starting
  from generic initial configurations (where, in particular, agents'
  positions are not aligned) contained in a bounded region of $\real^2$.
  We have consistently obtained that the time complexity to achieve
  $\task[\rendezvous]$ with $\FCC[\circumcenter]$ starting from these
  initial configurations is independent of the number of agents.  This
  leads us to conjecture that, in fact, initial configurations where all
  agents are aligned (i.e., the $1$-dimensional case) give rise to the
  worst possible performance of the algorithm. In more formal terms, we
  conjecture that, for $d \ge 2$,
  $\TC(\task[\rendezvous],\FCC[\circumcenter]) = \Theta (N)$. \oprocend
\end{remark}

\subsection{Deployment}
\label{se:deployment}

In this section we consider the uniform robotic network $\network[\realLD]$
presented in Example~\ref{ex:network-reald-ELDr} with parameter
$r\in\real_+$.  We assume we are given a convex simple polytope $Q \subset
\real^d$, with an integrable density function $\map{\phi}{Q}{\real_+}$.  We
assume that the initial positions of the agents belong to $Q$ and we intend
to design a control law that keeps them in $Q$ for subsequent times.  To
achieve the $\eps$-$r$-deployment task discussed in
Example~\ref{ex:deployment-tasks}, we define the \emph{centroid} control
and communication law $\FCC[\circumcenter]$.  This is a uniform, static,
time-independent law studied in~\cite{JC-SM-TK-FB:02j,JC-SM-FB:03p-tmp}.
Loosely speaking, the evolution of the network under the centroid control
and communication law can be described as follows:
\begin{quote}
  \emph{[Informal description]} Communication rounds take place at each
  natural instant of time.  At each communication round each agent performs
  the following tasks: (i) it transmits its position and receives its
  neighbors' positions; (ii) it computes the centroid of an appropriate
  region (the region is the intersection between the agent's Voronoi cell
  and a closed ball centered at its position and of radius $\frac{r}{2}$),
  and (iii) it moves toward this centroid.
\end{quote} 
Let us present this description in more formal terms. We set
$\timeschedule=\naturalzero$, $L=\real^d$, and
$\supind{\msg}{i}=\msgstandard$, $i \in I$.  We define the control function
$\map{\ctrl}{\real^d\times\real^d\times{L^N}}{\real^d}$ by
\begin{align*}
  \ctrl(x,\sampled{x},y) = \Centroid(\XX) - \sampled{x} \, ,
\end{align*}
where $\XX = Q \cap \cball{\frac{r}{2}}{\sampled{x}} \cap \big(\cap_{p \in
  \piSpace{L} (y)} H_{\sampled{x},p}\big)$ and $H_{\sampled{x},p}$ is the
half-space $\setdef{q\in \real^d}{ \Enorm{q-\sampled{x}} \le \Enorm{q-p}}$.
One can show that $Q^N$ is a positively-invariant set for this control law.

The following theorem on the centroid control and communication law
summarizes the known results about the asymptotic properties and the novel
results on the complexity of this law.  In characterizing complexity, we
assume $\diam(Q)$ is independent of $N$, $r$ and $\eps$, and we do not
calculate how the bounds depend on $r$.  As for the circumcenter law, we
provide complete time-complexity results only for the case $d=1$.

\begin{theorem}[Time and mean communication complexity of centroid law]
  \label{th:deployment-time-complexity}
  For $r\in\real_+$ and $\eps\in\real_+$, consider the network
  $\network[\realLD]$ with initial conditions in $Q$. The following
  statements hold:
  \begin{enumerate}
  \item for $d\in\natural$, the law $\FCC[\centroid]$ achieves the
    $\eps$-$r$-deployment task $\task[$\eps$-$r$-\deployment]$;
  \item for $d=1$ and $\phi=1$,
    $\TC(\task[$\eps$-$r$-\deployment],\FCC[\centroid]) \in
    O(N^3\log(N\eps^{-1}))$;
  \item for $d=1$, $\phi=1$ and for unidirectional communication,
    $\MCC(\task[$\eps$-$r$-\deployment],\FCC[\centroid]) \in \Theta(N)$.
    \oprocend
  \end{enumerate}
\end{theorem}
\begin{proof}
  Fact~(i) is proved in~\cite{JC-SM-FB:03p-tmp} for $d\in\{1,2\}$ and it is
  clear that the same proof technique can be generalized to any dimension.
  Fact~(iii) is proved in an analogous way to that of $\FCC[\circumcenter]$
  in Theorem~\ref{th:MCC-rendezvous}. In what follows we sketch the proof
  of fact~(ii).  For $d=1$, $Q$ is a compact interval on $\real$, say
  $Q=[q_-,q_+]$.
  
  We start with a brief discussion about connectivity.  Note that in the
  $r$-limited Delaunay graph, two agents on the line that are at most at a
  distance $r$ from each other are neighbors if and only if there are no
  other agents between them.  Additionally we claim that, if agents $i$ and
  $j$ are neighbors at time instant $\ell$, then
  $|\Centroid(\supind{\XX}{i}(\ell)) - \Centroid(\supind{\XX}{j}(\ell))|
  \le r$.  To see this, assume without loss of generality that
  $\supind{x}{i}(\ell) \le \supind{x}{j}(\ell)$.  Let us consider the case
  where the agents have neighbors on both sides (the other cases can be
  treated analogously).  Let $x_{-}^{[i]}(\ell)$ (respectively,
  $x_{+}^{[j]}(\ell)$) denote the position of the neighbor of agent $i$ to
  the left (respectively, of agent $j$ to the right). Now, we have
  \begin{align*}
    \Centroid(\supind{\XX}{i}(\ell)) = \frac{1}{4} (x_-^{[i]}(\ell) + 2
    \supind{x}{i}(\ell) + \supind{x}{j}(\ell)) \, , \quad
    \Centroid(\supind{\XX}{j}(\ell)) = \frac{1}{4} (\supind{x}{i}(\ell) + 2
    \supind{x}{j}(\ell) + x_+^{[j]}(\ell)) \, .
  \end{align*}
  Therefore, $|\Centroid(\supind{\XX}{i}(\ell)) -
  \Centroid(\supind{\XX}{j}(\ell))| \le \frac{1}{4} \big( | x_-^{[i]}(\ell)
  - \supind{x}{i}(\ell) | + 2| \supind{x}{i}(\ell) - \supind{x}{j}(\ell) |
  + |\supind{x}{j}(\ell) - x_+^{[j]}(\ell) | \big) \le r$. This implies
  that agents $i$ and $j$ are in the same connected component of the
  $r$-limited Delaunay graph at time instant $\ell +1$.
  
  Next, let us consider the case that $\subscr{E}{\rLD} (x_0)$ is
  connected.  Without loss of generality, assume that the agents are
  ordered from left to right according to their identifier, that is,
  $x^{[1]}(0) = (x_{0})_1 \le \dots \le x^{[N]}(0) = (x_{0})_N$.  We
  distinguish three cases depending on the proximity of the leftmost and
  rightmost agents $1$ and $N$, respectively, to the boundary of the
  environment: \fcase{a} both agents are within a distance
  $\frac{r}{2}$ of $\partial Q$; \fcase{b} none of the two is within
  a distance $\frac{r}{2}$ of $\partial Q$; and \fcase{c} only one
  of the agents is within a distance $\frac{r}{2}$ of~$\partial Q$.  Here
  is an important observation: from one time instant to the next one, the
  network configuration can fall into any of the cases described above.
  However, because of the discussion on connectivity, transitions can only
  occur from case \fcase{b} to either case \fcase{a} or
  \fcase{c}; and from case \fcase{c} to case
  \fcase{a}.  As we show in the following, for each of these cases,
  the network evolution under $\FCC[\centroid]$ can be described as a
  discrete-time linear dynamical system which respects agents' ordering.
  
  Let us consider case~\fcase{a}. In this case, we have
  \begin{multline*}
    x^{[1]}(\ell + 1) = \frac{1}{4} (x^{[1]}(\ell)+ x^{[2]}(\ell)) +
    \frac{1}{2}q_- \, , \quad x^{[2]}(\ell + 1) = \frac{1}{4}
    (x^{[1]}(\ell) + 2 x^{[2]}(\ell) + x^{[3]}(\ell)) \,, \; \dots
    \, ,  \\
    \dots \, , \; x^{[N-1]}(\ell + 1) = \frac{1}{4} (x^{[N-2]}(\ell) + 
    2 x^{[N-1]}(\ell) + x^{[N]}(\ell) ) \,, \quad x^{[N]}(\ell + 1) =
    \frac{1}{4} (x^{[N-1]}(\ell)+ x^{[N]}(\ell) ) + \frac{1}{2}q_+\,.
  \end{multline*}
  Equivalently, we can write $x(\ell+1) = A_{\mathfrak{a}} \cdot x(\ell) +
  b_{\mathfrak{a}}$, where the $N\times N$-matrix $A_{\mathfrak{a}}$ and
  the vector $b_{\mathfrak{a}}$ are given by
  \begin{equation*}
    A_{\mathfrak{a}} = 
    \begin{bmatrix}
      \frac{1}{4} & \frac{1}{4} & 0 & \dots & \dots & 0 \\
      \frac{1}{4} & \frac{1}{2} & \frac{1}{4} & \dots & \dots & 0 \\
      0 & \frac{1}{4} & \frac{1}{2} & \frac{1}{4} & \dots &0 \\
      \vdots & & \ddots & \ddots & \ddots & \vdots \\
      0 & \dots & \dots & \frac{1}{4} & \frac{1}{2} & \frac{1}{4} \\ 
      0 & \dots & \dots & 0 & \frac{1}{4} & \frac{1}{4}
    \end{bmatrix} , \quad 
    b_{\mathfrak{a}} = 
    \begin{bmatrix}
      \frac{1}{2}q_-\\
      0\\
      \vdots\\
      0\\
      \frac{1}{2}q_+
    \end{bmatrix}
    \!.
  \end{equation*}
  Note that the only equilibrium network configuration $x_*$
  respecting the ordering of the agents is given by
  \begin{align*}
    x_*^{[i]} = q_- + \frac{1}{2N} (1+ 2(i-1)) (q_+ - q_-) \, , \quad
    i \in I \, ,
  \end{align*}
  and note that this is a $\frac{r}{2}$-centroidal Voronoi configuration
  (under the assumption of case \fcase{a}).  We can therefore write
  $(x(\ell)-x_*)=A_{\mathfrak{a}}(x(\ell-1)-x_*)$.  Now, note that
  $A_{\mathfrak{a}} =\atridmat_N^{-}\big(\frac{1}{4},\frac{1}{2}\big)$ as
  defined in Appendix~\ref{app:toeplitz}.
  Theorem~\ref{thm:convergence-atrid}(ii) implies that
  $\lim_{\ell\to+\infty} \big(x(\ell)-x_*\big) =\mathbf{0}$, and that the
  maximum time required for $\Enorm{x(\ell)-x_*\big}\leq
  \eps\Enorm{x(0)-x_*}$ (over all initial conditions $x(0)\in\real^N$) is
  $\Theta\big(N^2\log\eps^{-1}\big)$.  It is not obvious, but it can be
  verified, that the initial condition providing the lower bound in the
  time complexity estimate does indeed have the property of respecting the
  agents' ordering; this fact holds for all three cases \fcase{a},
  \fcase{b} and \fcase{c}.
  
  The case \fcase{b} can be treated in the same way. The network
  evolution takes now the form $x(\ell+1) = A_{\mathfrak{b}} \cdot x(\ell)
  + b_{\mathfrak{b}}$, where the $N\times N$-matrix $A_{\mathfrak{b}}$ and
  the vector $b_{\mathfrak{b}}$ are given by
  \begin{equation*}
    A_{\mathfrak{b}} = 
    \begin{bmatrix}
      \frac{3}{4} & \frac{1}{4} & 0 & \dots & \dots & 0 \\
      \frac{1}{4} & \frac{1}{2} & \frac{1}{4} & \dots & \dots & 0 \\
      0 & \frac{1}{4} & \frac{1}{2} & \frac{1}{4} & \dots &0 \\
      \vdots & & \ddots & \ddots & \ddots & \vdots \\
      0 & \dots & \dots & \frac{1}{4} & \frac{1}{2} & \frac{1}{4} \\ 
      0 & \dots & \dots & 0 & \frac{1}{4} & \frac{3}{4}
    \end{bmatrix} \! , \quad 
    b_{\mathfrak{b}} = 
    \begin{bmatrix}
      - \frac{1}{4}r\\
      0\\
      \vdots\\
      0\\
      \frac{1}{4}r
    \end{bmatrix}
    \!.
  \end{equation*}
  In this case, a (non-unique) equilibrium network configuration
  respecting the ordering of the agents is of the form
  \begin{align*}
    \supind{x}{i}_* = i r -\frac{1+N}{2}r \, , \quad i \in I \, .
  \end{align*}
  Note that this is a $\frac{r}{2}$-centroidal Voronoi configuration (under
  the assumption of case \fcase{b}).  We can therefore write
  $(x(\ell)-x_*)=A_{\mathfrak{b}}(x(\ell-1)-x_*)$.  Now, note that
  $A_{\mathfrak{b}} =\atridmat_N^{+}\big(\frac{1}{4},\frac{1}{2}\big)$ as
  defined in Appendix~\ref{app:toeplitz}.  As in the notation of
  Theorem~\ref{thm:convergence-atrid}(i) we compute
  $\xave=\frac{1}{N}\mathbf{1}^T(x_0-x_*) = \frac{1}{N}\mathbf{1}^Tx_0$.
  With this calculation, Theorem~\ref{thm:convergence-atrid}(i) implies
  that $\lim_{\ell\to+\infty} \big(x(\ell)-x_* -\xave\mathbf{1} \big)
  =\mathbf{0}$, and that the maximum time required for
  $\Enorm{x(\ell)-x_*-\xave\mathbf{1}\big}\leq
  \eps\Enorm{x(0)-x_*-\xave\mathbf{1}}$ (over all initial conditions
  $x(0)\in\real^N$) is $\Theta\big(N^2\log\eps^{-1}\big)$.
  
  Case \fcase{c} needs to be handled differently. Without loss of
  generality, assume that agent $1$ is within distance $\frac{r}{2}$
  of~$\partial Q$ and agent $N$ is not (the other case is treated
  analogously). Then, the network evolution takes now the form $x(\ell+1) =
  A_{\mathfrak{c}} \cdot x(\ell) + b_{\mathfrak{c}}$, where the $N\times
  N$-matrix $A_{\mathfrak{c}}$ and the vector $b_{\mathfrak{c}}$ are given
  by
  \begin{equation*}
    A_{\mathfrak{c}} = 
    \begin{bmatrix}
      \frac{1}{4} & \frac{1}{4} & 0 & \dots & \dots & 0 \\
      \frac{1}{4} & \frac{1}{2} & \frac{1}{4} & \dots & \dots & 0 \\
      0 & \frac{1}{4} & \frac{1}{2} & \frac{1}{4} & \dots &0 \\
      \vdots & & \ddots & \ddots & \ddots & \vdots \\
      0 & \dots & \dots & \frac{1}{4} & \frac{1}{2} & \frac{1}{4} \\
      0 & \dots & \dots & 0 & \frac{1}{4} & \frac{3}{4}
    \end{bmatrix} \! , \quad 
    b_{\mathfrak{c}} = 
    \begin{bmatrix}
      \frac{1}{2}q_-\\
      0\\
      \vdots\\
      0\\
      \frac{1}{4}r
    \end{bmatrix}
    \!.
  \end{equation*}
  Note that the only equilibrium network configuration $x_*$ respecting the
  ordering of the agents is given by
  \begin{align*}
    x_*^{[i]} = q_- + \frac{1}{2} (2i-1) r \, , \quad i \in I
    \, ,
  \end{align*}
  and note that this is a $\frac{r}{2}$-centroidal Voronoi configuration
  (under the assumption of case \fcase{c}).  In order to analyze
  $A_{\mathfrak{c}}$, we recast the $N$-dimensional discrete-time dynamical
  system as a $2N$-dimensional one.  To do this, we define a
  $2N$-dimensional vector $y$ by
  \begin{equation}\label{eq:trick}
    \supind{y}{i}=\supind{x}{i}, i\in{I}, \quad \text{and} \quad
    \supind{y}{N+i}=\supind{x}{N-i+1}, i\in{I},
  \end{equation}
  Now, one can see that the network evolution can be alternatively
  described in the variables $(y^{[1]},\dots,y^{[2N]})$ as a linear
  dynamical system determined by the $2N \times 2N$ matrix
  $\atridmat_{2N}^{-}(\frac{1}{4},\frac{1}{2})$.  Using analogous arguments
  to the ones used before in Theorems~\ref{thm:convergence-trid-circ} and
  \ref{thm:convergence-atrid} and exploiting the chain of
  equalities~\eqref{eq:trick}, we can characterize the eigenvalues and
  eigenvectors of $\tridmat_{2N-1}(\frac{1}{4},\frac{1}{2},\frac{1}{4})$,
  and infer that, even for case~\fcase{c}, the maximum time required for
  $\Enorm{x(\ell)-x_*\big}\leq \eps\Enorm{x(0)-x_*}$ (over all initial
  conditions $x(0)\in\real^N$) is $\Theta\big(N^2\log\eps^{-1}\big)$.
  %% \margin{for all: this is only a sketch.}
  
  In summary, for all three cases \fcase{a}, \fcase{b} and \fcase{c}, our
  calculations show that, in time $O\big(N^2\log\eps^{-1}\big)$, the error
  $2$-norm satisfies the contraction inequality
  $\Enorm{x(\ell)-x_*\big}\leq \eps\Enorm{x(0)-x_*}$.  We convert this
  inequality on $2$-norms into an appropriate inequality on $\infty$-norms
  as follows.  Note that $\Infnorm{x(0) - x_*} = \max_{i\in{I}}
  |\supind{x}{i}(0) - x_*^{[i]}| \leq (q_+-q_-)$.  For $\ell$ of order
  $N^2\log\eta^{-1}$, we have:
  \begin{align*}
    \Infnorm{x(\ell)-x_*} &\leq \Enorm{x(\ell)-x_*} \leq \eta
    \Enorm{x(0)-x_*} \leq \eta\sqrt{N} \Infnorm{x(0)-x_*} \leq \eta
    \sqrt{N}(q_+-q_-).
  \end{align*}
  This means that $\eps$-$r$-deployment is achieved for $\eta
  \sqrt{N}(q_+-q_-)=\eps$, that is, in time
  $O(N^2\log\eta^{-1})=O(N^2\log(N\eps^{-1}))$.
 
  Up to here we have proved that, if the graph $(I,\subscr{E}{\rLD}(x_0))$
  is connected, then $\TC (\task[$\eps$-$r$-\deployment],\FCC[\centroid])
  \in O(N^2\log(N\eps^{-1}))$.  If $(I,\subscr{E}{\rLD}(x_0))$ is not
  connected, note that along the network evolution there can only be a
  finite number of time instants, at most $N-1$ where a merging of two
  connected components occurs. Therefore, the time complexity is at most
  $O(N^3\log(N\eps^{-1}))$.
\end{proof}
%% {Incomplete proof of lower bound for deployment; see secX-leftout}

\section{Conclusions}\label{se:conclusions}
We have introduced a formal model for the design and analysis of
coordination algorithms executed by networks composed of robotic agents.
In our discrete-time communication, continuous-time motion model, the
robotic agents evolve in the physical domain in continuous-time and have
the ability to exchange information (position and/or dynamic variables)
that affect their motion at discrete-time instants.  Under this framework,
motion coordination algorithms are formalized as feedback control and
communication laws. Drawing analogies with the classical theory on
distributed algorithms, we have defined two measures of complexity for this
formal notion of execution: the time and the mean communication complexity
of achieving a specific task.  We have defined the notion of re-scheduling
of a control and communication law and analyzed the invariance of the
proposed complexity measures under this operation.  These concepts and
results have been illustrated with various examples of different types of
robotic networks and different classes of feedback control and
communication laws.  Finally, we have computed the proposed complexity
measures for a variety of algorithms performing spatially-distributed tasks
such as rendezvous, pursuit, and deployment.

Throughout the paper, we have made compromises in order to arrive at a
simple yet useful model that is able to capture a wide class of algorithms.
However, it is important to acknowledge that our results on network, task
and complexity modeling, and our results on algorithm analysis have
numerous limitations and open up numerous avenues for future research.  Let
us at provide a concise list here: (1) modeling of asynchronous as opposed
to only synchronous networks (see
however~\cite{JL-ASM-BDOA:04b,JC-SM-TK-FB:02j,PF-GP-NS-PW:05,XD-AK:02});
(2) models and analysis of failures in the agents (arrivals/departures) and
in the communication links (see
however~\cite{JC-SM-FB:04h-tmp,ROS-RMM:03c,LM:04,WR-RWB:03}); (3)
probabilistic versions of the complexity measures, as opposed to only
worst-case analysis (see however~\cite{EK:02a}); (4) quantization and
delays in the communication channels (see however~\cite{FF-KHJ-AS-SZ:04}
and the literature on quantized control); (5) parallel, sequential and
hierarchical composition of control and communication laws.  Additionally,
our analysis results essentially consist of a time-complexity analysis of
some basic algorithms, but many more open algorithmic questions remain
unresolved including (6) communication complexity for omnidirectional
communication models; (7) analysis of known algorithms for flocking and
cohesion; (8) complexity analysis of tasks as opposed to algorithms.

\section*{Acknowledgments}
This material is based upon work supported in part by ONR YIP Award
N00014-03-1-0512, NSF SENSORS Award IIS-0330008, and NSF CAREER Award
CCR-0133869. Sonia Mart{\'\i}nez's work was supported in part by a
Fulbright Postdoctoral Fellowship from the Spanish Ministry of Education
and Science.

{\small
%\bibliographystyle{plainnat}
%\bibliography{alias,FB,New,Main}

}

\appendix
\section{Basic geometric notions}
\label{app:geometry}

Let $S\subset\real^d$, $d\in\natural$, be compact.  The \emph{circumcenter
  of $S$}, denoted by $\CircumC (S)$, is the center of the smallest-radius
sphere in $\real^d$ enclosing $S$. Given an integrable function
$\map{\phi}{S}{\real_+}$, the \emph{mass} of $S$ is
$\Mass(S)=\int_S\phi(q)dq$, and the \emph{centroid} of $S$ is
\begin{align*}
  \Centroid(S) = \frac{1}{\Mass(S)} \int_S q \phi (q)dq \, .
\end{align*}
A \emph{partition} of $S$ is a collection of subsets of $S$ with disjoint
interiors and whose union is $S$.  Given a set of $N$ distinct points
$\PP=\{p_i\}_{i\in\until{N}}$ in $S$, the \emph{Voronoi partition of $S$}
generated by $\PP$ (with respect to the Euclidean norm) is the collection
of sets $\{V_i(\PP)\}_{i\in\until{N}}$ defined by $V_i(\PP) = \{q\in S \; |
\; \Enorm{q-p_i}\leq \Enorm{q-p_j} \, , \; \text{for all} \; p_j\in\PP\}$.
We usually refer to $V_i(\PP)$ as $V_i$.  For a detailed treatment of
Voronoi partitions we refer to~\cite{MdB-MvK-MO:97,AO-BB-KS-SNC:00}.

\begin{figure}[htbp]
  \psfrag{G2}[cc][cc]{{\small $r$-disk graph}}%
  \psfrag{G4}[cc][cc]{{\small $r$-lim. Delaunay graph}}% 
  \centering
  \fbox{\includegraphics[width=.31\linewidth]{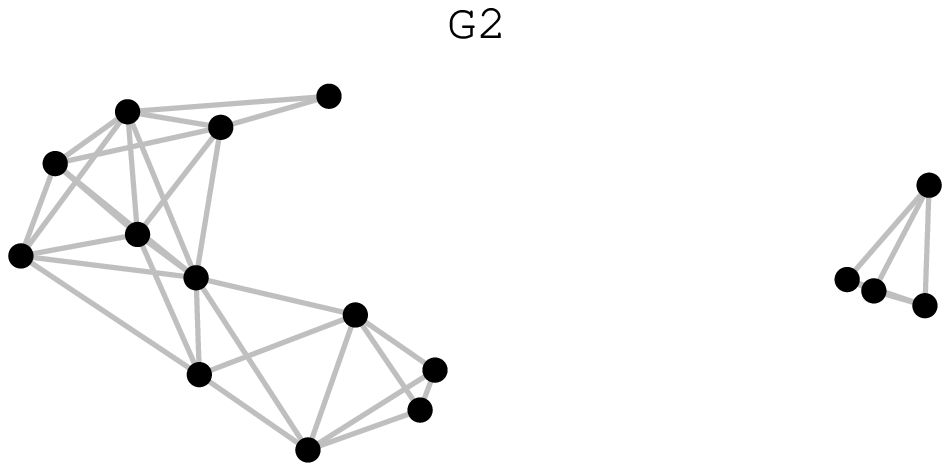}}%
  \fbox{\includegraphics[width=.31\linewidth]{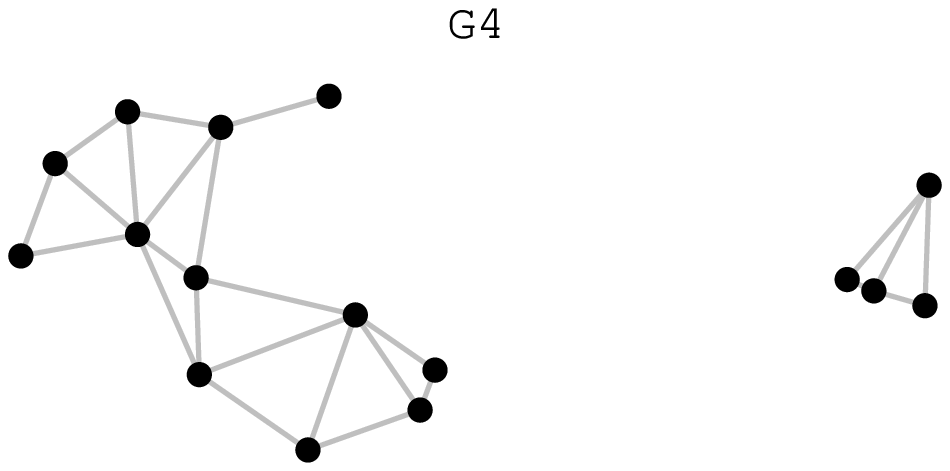}}%
  \caption{The $r$-disk and  $r$-limited Delaunay graphs in $\real^2$}
  \label{fig:example-graphs-2D}
\end{figure}

For $I=\until{N}$ and $S\subset\real^d$, a \emph{proximity edge map} is a
map of the form $\map{E}{S^N}{2^{I\times{I}\setminus\diag(I\times{I})}}$.
For $r\in\real_+$, we define the $r$-disk proximity edge map
$\map{\subscr{E}{\rdisk}}{(\real^d)^N}{2^{I\times{I}}}$ and the $r$-limited
Delaunay proximity edge map
$\map{\subscr{E}{\rLD}}{(\real^d)^N}{2^{I\times{I}}}$ by
\begin{subequations}
  \begin{align}
    \subscr{E}{\rdisk} (x_1,\dots,x_N) & = \setdef{(i,j)\in I \times I
      \setminus \diag (I \times I)}{ \Enorm{x_i - x_j} \le r} \, ,
    \label{eq:E-r}\\
    \subscr{E}{LD,r} (x_1,\dots,x_N) & = \setdef{(i,j)\in I \times I
      \setminus \diag (I \times I)}{\big(V_i \cap \cball{\tfrac{r}{2}}{x_i}
      \big) \intersection \big(V_j \cap \cball{\tfrac{r}{2}}{x_j} \big)
      \neq \emptyset } \, ,
    \label{eq:E-rD}
  \end{align}
\end{subequations} 
where $\{V_1,\dots,V_N\}$ is the Voronoi partition of $\real^d$ generated
by $\{x_1,\dots,x_N\}$. Illustrations of these concepts are given in
Fig.~\ref{fig:example-graphs-2D}.   
For $d=1$, the $r$-limited Delaunay proximity edge map has the following
intuitive characterization: two points are neighbors if and only if they
are within distance $r$ and no other point is between them.

As proved in~\cite{JC-SM-FB:03p-tmp}, the $r$-limited Delaunay graph and
the $r$-disk graph have the same connected components.  Additionally, the
$r$-limited Delaunay graph is ``computable'' on the $r$-disk graph in the
following sense: any node in the network can compute the set of its
neighbors in the $r$-limited Delaunay graph if it is given the set of its
neighbors in the $r$-disk graph.  This implies that any control and
communication law for a network with communication graph $\subscr{E}{\rLD}$
can be implemented on a analogous network with communication graph
$\subscr{E}{\rdisk}$.

\section{Tridiagonal Toeplitz and circulant matrices}
\label{app:toeplitz}%%
We briefly review some basic facts about certain classes of Toeplitz
matrices, see~\cite{CDM:01}.  For $N\geq 2$ and $a,b,c\in\real$, define the
$N \times N$ Toeplitz matrices $\tridmat_N(a,b,c)$ and $\circmat_N(a,b,c)$
by
\begin{align*}
  \tridmat_N(a,b,c) &=
  \begin{bmatrix}
    b & c & 0 & \dots & 0 \\
    a & b & c & \dots & 0 \\
    \vdots &  \ddots & \ddots & \ddots & \vdots \\
    0 & \dots &  a & b & c \\
    0 & \dots & 0 & a & b
  \end{bmatrix},
  \quad \text{and} \quad
  \circmat_N(a,b,c) = \tridmat_N(a,b,c) + 
  \begin{bmatrix}
    0 & \dots  & \dots & 0 & a \\
    0 & \dots  & \dots & 0 & 0 \\
    \vdots & \ddots & \ddots & \ddots & \vdots \\
    0 & 0 & \dots &  0 & 0  \\
    c & 0 & \dots &  0 & 0  
  \end{bmatrix}.
\end{align*}
The matrices $\tridmat_N$ and $\circmat_N$ are tridiagonal and circulant,
respectively.  The two matrices only differ in their $(1,N)$ and $(N,1)$
entries. Note our convention that $C_2(a,b,c)=
\begin{bmatrix}
  b & a+c \\ a+c & b
\end{bmatrix}$.
\begin{lemma}[Eigenvalues of tridiagonal Toeplitz and circulant matrices]
  \label{lem:eig-trid-circ}
  For $N\geq 2$ and $a,b,c\in\real$, the following statements hold:
  \begin{enumerate}
  \item for $ac\neq 0$, the eigenvalues and eigenvectors of
    $\tridmat_N(a,b,c)$ are
    \begin{equation*}
      b + 2 c \sqrt{\frac{a}{c}} \cos \left(
      \frac{\displaystyle i\pi}{\displaystyle N+1}\right), 
    \quad 
      \begin{bmatrix}
        \left(\frac{a}{c}\right)^{1/2} \sin \left(\frac{i\pi}{N+1} \right) \\
        \left(\frac{a}{c}\right)^{2/2}\sin \left(\frac{2i\pi}{N+1} \right)\\
        %%\left(\frac{a}{c}\right)^{3/2}\sin\left(\frac{3i\pi}{N+1} \right)\\ 
        \vdots \\
        \left(\frac{a}{c}\right)^{N/2}\sin \left(\frac{N i\pi}{N+1}\right)
      \end{bmatrix},   
    \quad i\in\{1,\dots,N\};
    \end{equation*}
  \item the eigenvalues and eigenvectors of $\circmat_N(a,b,c)$ are
    (note $\omega = \exp (\frac{2 \pi \sqrt{-1}}{N})$)
    \begin{align*}
      b + (a+c) \cos\left( \frac{i 2\pi}{N} \right) + \sqrt{-1} (c-a)
      \sin\left( \frac{i 2\pi}{N}\right), 
    \quad 
      \begin{bmatrix}
        1 \\
        \omega^i \\
        \vdots \\
        \omega^{(N-1)i}
      \end{bmatrix}
      \!,   
      \quad i\in\{1,\dots,N\}.
      \eqoprocend
    \end{align*}
  \end{enumerate}
\end{lemma}
\begin{proof}
  Both facts are discussed, for example, in~\cite[Example 7.2.5 and
  Exercise 7.2.20]{CDM:01}. Fact~(ii) requires some straightforward
  algebraic manipulations.
\end{proof}

\begin{remarks}
  \begin{enumerate}
  \item The set of eigenvalues of $\tridmat_N(a,b,c)$ is contained in the
    real interval $[b-2\sqrt{ac},b+2\sqrt{ac}]$, if $ac\geq 0$, and in the
    interval in the complex plane
    $[b-2\sqrt{-1}\sqrt{|ac|},b+2\sqrt{-1}\sqrt{|ac|}]$, if $ac \leq 0$.
  \item The set of eigenvalues of $\circmat_N(a,b,c)$ is contained in the
    ellipse on the complex plane with center $b$ and horizontal axis
    $2|a+c|$ and vertical axis $2|c-a|$.  
   
  \item Recall from \cite{CDM:01} that (1) a square matrix is normal if it
    has a complete orthonormal set of eigenvectors, (2) circulant matrices
    and real-symmetric matrices are normal, and (3) if a normal matrix has
    eigenvalues $\{\lambda_1,\dots,\lambda_n\}$, then its singular values
    are $\{|\lambda_1|,\dots,|\lambda_n|\}$.\oprocend
    %% a normal matrix has $n$ orthogonal eigenvectors
    %% a normal matrix is unitarily diagonalizable
    %% a normal matrix commutes with its tranjugate   
  \end{enumerate}  
\end{remarks}

We can now state the main result of this section.  

\begin{theorem}[Tridiagonal Toeplitz and circulant dynamical systems]
  \label{thm:convergence-trid-circ}
  Let $N\geq 2$, $\eps\in]0,1[$, and $a,b,c\in\real$.  Let
  $\map{x}{\naturalzero}{\real^N}$ and $\map{y}{\naturalzero}{\real^N}$ be
  solutions to
  \begin{equation*}
    x(\ell+1) = \tridmat_N(a,b,c) \, x(\ell),
    \quad \text{and}\quad
    y(\ell+1) = \circmat_N(a,b,c) \, y(\ell),
  \end{equation*}
  with initial conditions $x(0)=x_0$ and $y(0)=y_0$, respectively. The
  following statements hold: 
  \begin{enumerate}
  \item if $a=c\neq 0$ and $|b|+2|a|=1$, then
    $\lim_{\ell\to+\infty}x(\ell)=\mathbf{0}$, and the maximum time
    required for $\Enorm{x(\ell)}\leq \eps \Enorm{x_0}$ (over all initial
    conditions $x_0\in\real^N$) is $\Theta\big( N^2 \log\eps^{-1}\big)$;
    
  \item if $a\neq 0$, $c=0$ and $0<|b|<1$, then
    $\lim_{\ell\to+\infty}x(\ell)=\mathbf{0}$, and the maximum time
    required for $\Enorm{x(\ell)}\leq \eps \Enorm{x_0}$ (over all initial
    conditions $x_0\in\real^N$) is $O\big(N\log{N}+ \log\eps^{-1}\big)$;
    
  \item if $a\geq 0$, $c\geq 0$, $b>0$, and $a+b+c=1$, then
    $\lim_{\ell\to+\infty}y(\ell)=\yave\mathbf{1}$, where
    $\yave=\frac{1}{N}\mathbf{1}^Ty_0$, and the maximum time required for
    $\Enorm{y(\ell)-\yave\mathbf{1}}\leq \eps\Enorm{y_0-\yave\mathbf{1}}$
    (over all initial conditions $y_0\in\real^N$) is
    $\Theta\big(N^2\log\eps^{-1}\big)$.\oprocend
  \end{enumerate}
\end{theorem} 
\begin{proof}
  Let us prove fact~(i).  We start by bounding from above the eigenvalue
  with largest absolute value, that is, the largest singular value, of
  $\tridmat_N(a,b,a)$: 
  \begin{align*} 
    \max_{ i\in\{1,\dots,N\}} \left| b + 2 a \cos \left(
        \frac{\displaystyle i\pi}{\displaystyle N+1}\right) \right| \leq
    |b| + 2 |a| \max_{ i\in\{1,\dots,N\}} \left| \cos \left(
        \frac{\displaystyle i\pi}{\displaystyle N+1}\right)\right| \leq |b|
    + 2 |a| \cos \left( \frac{\displaystyle \pi}{\displaystyle N+1}\right).
  \end{align*} 
  Because $\cos(\frac{\pi}{N+1})<1$ for any $N\geq 2$, the matrix
  $\tridmat_N(a,b,a)$ is stable.  Additionally, for $\ell>0$, we bound from
  above the magnitude of the curve $x$ as
  \begin{align*} 
    \Enorm{x(\ell) } &=\Enorm{ \tridmat_N(a,b,a)^\ell x_0 } %% 
    \leq \left( |b|
      + 2 |a| \cos \left( \frac{\displaystyle \pi}{\displaystyle
          N+1}\right) \right)^\ell \, \Enorm{ x_0 } .  
  \end{align*} 
  In order to have $\Enorm{x(\ell)}<\eps\Enorm{x_0}$, it is sufficient that
  $\ell \log \left( |b| + 2 |a| \cos \left( \frac{\displaystyle
        \pi}{\displaystyle N+1}\right) \right) < \log \eps$, that is
  \begin{gather} \label{eq:ell-bound} \ell > \frac{\log\eps^{-1}} {- \log
      \left( |b| + 2 |a| \cos \left( \frac{\displaystyle \pi}{\displaystyle
            N+1}\right) \right) } .  
  \end{gather} 
  To show the upper bound, note that as $t\to{0}$ we have 
  \begin{equation*}
    -\frac{1}{\log(1-2|a|(1-\cos t))} = \frac{1}{|a| t^2} +
    O(1). %%\frac{1/a-3}{6} + O(x^2). 
  \end{equation*} 
  Now, assume without loss of generality that $ab>0$ and consider the
  eigenvalue $b+2a\cos(\frac{\pi}{N+1})$ of $\tridmat_N(a,b,a)$.  Note that
  $|b+2a\cos(\frac{\pi}{N+1})|=|b|+2|a|\cos(\frac{\pi}{N+1})$.  (If $ab<0$,
  then consider the eigenvalue $b+2a\cos(\frac{N\pi}{N+1})$.)  For $N>2$,
  define the unit-length vector
  \begin{equation}  \label{eq:vN} 
    \mathbf{v}_N = \sqrt{\frac{2}{N+1}} \begin{bmatrix}
      \sin\frac{\pi}{N+1} \\ \vdots \\ \sin\frac{N\pi}{N+1} \end{bmatrix}
    \in \real^{N}, 
  \end{equation} 
  and note that, by Lemma~\ref{lem:eig-trid-circ}(i), $\mathbf{v}_N$ is an
  eigenvector of $\tridmat_N(a,b,a)$ with eigenvalue
  $b+2a\cos(\frac{\pi}{N+1})$.  The trajectory $x$ with initial condition
  $\mathbf{v}_N$ satisfies $\Enorm{x(\ell)}=\left( |b| + 2 |a|
    \cos\left(\frac{\pi}{N+1}\right)\right)^\ell \Enorm{\mathbf{v}_N}$ and,
  therefore, it will enter $\oball{\eps\Enorm{\mathbf{v}_N}}{\mathbf{0}}$
  only when $\ell$ satisfies~\eqref{eq:ell-bound}. This completes the proof
  of fact~(i).
    
  Next we consider statement~(ii). Clearly, $\tridmat_N(a,b,0)$ is stable.
  For $\ell>0$, we compute
  \begin{gather*}
    \tridmat_N(a,b,0)^\ell = b^\ell \Big( I_N + \frac{a}{b}
    \tridmat_N(1,0,0) \Big)^\ell = b^\ell \sum_{j=0}^{N-1} \frac{\ell!}{j!
      (\ell-j)!}  \Big(\frac{a}{b}\Big)^{j} \tridmat_N(1,0,0)^{j}
  \end{gather*}
  because of the nilpotency of $\tridmat_N(1,0,0)$.  Now we can bound from
  above the magnitude of the curve $x$ as
  \begin{align*}
    \Enorm{x(\ell)} &= \Enorm{\tridmat_N(a,b,0)^\ell x_0 } \leq |b|^\ell
    \sum_{j=0}^{N-1} \frac{\ell!}{j! (\ell-j)!}
    \Big(\frac{a}{b}\Big)^{j}\, \bignorm{ \tridmat_N(1,0,0)^{j} x_0}{2}
    \\
    &\leq \textup{e}^{a/b} \ell^{N-1} \left| b \right|^\ell \norm{x_0}{2}.
  \end{align*}
  Here we used $\norm{\tridmat_N(1,0,0)^{j}x_0}{2}\leq\norm{x_0}{2}$ and
  $\max\setdef{\frac{\ell!}{(\ell-j)!}}{j\in\{0,\dots,N-1\}}\leq\ell^{N-1}$.
  Therefore, in order to have $\Enorm{x(\ell)}<\eps\Enorm{x_0}$, it
  suffices that $ \log(\textup{e}^{a/b}) + (N-1) \log \ell + \ell \log|b|
  \leq \log \eps$, that is
  \begin{gather*}
    \ell - \frac{N-1}{-\log|b|} \log \ell %%
    > \frac{\frac{a}{b} - \log\eps } {-\log|b|} .
  \end{gather*} 
  A sufficient condition for $\ell-\alpha\log\ell>\beta$, for
  $\alpha,\beta>0$, is that $\ell\geq 2\beta + 2\alpha
  \max\{1,\log\alpha\}$. For, if $\ell\geq 2\alpha$, then $\log\ell$
  is bounded from above by the line $\ell/2\alpha+\log\alpha$.
  Furthermore, the line $\ell/2\alpha+\log\alpha$ is a lower bound for
  the line $(\ell-\beta)/\alpha$ if $\ell\geq 2\beta + 2\alpha
  \log\alpha$.  In summary, it is true that
  $\Enorm{x(\ell)}\le\eps\Enorm{x(0)}$ whenever
  \begin{gather*}
    \ell\geq 2\frac{ \frac{a}{b} - \log\eps }
    {-\log|b|} + 2\frac{N-1}{-\log|b|} \max\left\{
      1,\log\frac{N-1}{-\log|b|} \right\}.
  \end{gather*}
  This completes the proof of the upper bound, that is, fact~(ii).
  
  The proof of fact~(iii) is similar to that of fact~(i). We analyze the
  singular values of $\circmat_N(a,b,c)$. It is clear that the eigenvalue
  corresponding to $i=N$ is equal to $1$; this is the largest singular
  value of $\circmat_N(a,b,c)$ and the corresponding eigenvector is
  $\mathbf{1}$.  In the orthogonal decomposition induced by the
  eigenvectors of $\circmat_N(a,b,c)$, the vector $y_0$ has a component
  $\yave$ along the eigenvector $\mathbf{1}$.  We now compute the second
  largest singular value:
  \begin{multline*}
    \max_{i\in\{1,\dots,N-1\}} \left| b +(a+c)
      \cos\left(\frac{i2\pi}{N}\right) + \sqrt{-1}(c-a)
      \sin\left(\frac{i2\pi}{N}\right) \right| \\
    = \left| 1-(a+c) \Big(1- \cos\big(\frac{2\pi}{N}\big) \Big) +
      \sqrt{-1}(c-a) \sin\left(\frac{2\pi}{N}\right) \right|.
  \end{multline*}
  Here $|\cdot|$ is the norm in $\complex$.  Because of the assumptions on
  $a,b,c$, the second largest singular value is strictly less than $1$.
  For $\ell>0$, we bound the distance of the curve $y(\ell)$ from
  $\yave\mathbf{1}$ as
  \begin{align*}
    \Enorm{ y(\ell) - \yave\mathbf{1}}
    &=\Enorm{ \circmat_N(a,b,c)^\ell y_0 - \yave\mathbf{1} }
    =\Enorm{ \circmat_N(a,b,c)^\ell \big(  y_0 - \yave\mathbf{1} \big)} 
    \\
    &\leq \left| 1-(a+c) \Big(1- \cos\big(\frac{2\pi}{N}\big) \Big) +
      \sqrt{-1}(c-a) \sin\left(\frac{2\pi}{N}\right) \right|^\ell \,
    \Enorm{ y_0 - \yave\mathbf{1} }.
  \end{align*}
  This proves that $\lim_{\ell\to+\infty}y(\ell)=\yave\mathbf{1}$.  Also,
  for $\alpha=a+c,\beta=c-a$ and as $t\to{0}$, we have
  \begin{equation*}
    -\frac{1}%
    {\log\Big( \big( 1-\alpha(1-\cos t) \big)^2+\beta^2\sin^2t \Big)^{1/2}}
    =
    \frac{2}{(\alpha-\beta^2)t^2} + O(1).
  \end{equation*}
  Here $\beta^2<\alpha$ because $a,c\in ]0,1[$.

  Now, consider the eigenvalues $\lambda_N = b +(a+c)
  \cos\left(\frac{2\pi}{N}\right) + \sqrt{-1}(c-a)
  \sin\left(\frac{2\pi}{N}\right)$ and $\overline{\lambda}_N = b
  +(a+c) \cos\left(\frac{(N-1)2\pi}{N}\right) + \sqrt{-1}(c-a)
  \sin\left(\frac{(N-1)2\pi}{N}\right)$ of $\circmat_N(a,b,c)$, and
  its associated eigenvectors (cf. Lemma~\ref{lem:eig-trid-circ}(ii))
  \begin{equation}
    \label{eq:vN-2}
    \mathbf{v}_N = 
    \begin{bmatrix}
      1 \\ w \\ \vdots \\ w^{N-1}
    \end{bmatrix} \in \complex^{N}, \quad
    \overline{\mathbf{v}}_N =
    \begin{bmatrix}
      1 \\ w^{N-1} \\ \vdots \\ w
    \end{bmatrix} \in \complex^{N} .
  \end{equation}
  Note that the vector $\mathbf{v}_N + \overline{\mathbf{v}}_N$
  belongs to $\real^N$. Moreover, its component $\yave$ along the
  eigenvector $\mathbf{1}$ is $0$.  The trajectory $y$ with initial
  condition $\mathbf{v}_N + \overline{\mathbf{v}}_N$ satisfies
  $\Enorm{y(\ell)} = \Enorm{ \lambda_N^\ell \mathbf{v}_N +
    \overline{\lambda}_N^\ell \overline{\mathbf{v}}_N} =
  |\lambda_N|^\ell \Enorm{\mathbf{v}_N + \overline{\mathbf{v}}_N}$
  and, therefore, it will enter $\oball{\eps\Enorm{\mathbf{v}_N +
      \overline{\mathbf{v}}_N}}{\mathbf{0}}$ only when
  \begin{equation*}
    \ell > \frac{\log\eps^{-1}} {- \log 
      \left| 1-(a+c) \Big(1- \cos\big(\frac{2\pi}{N}\big) \Big) +
        \sqrt{-1}(c-a) \sin\left(\frac{2\pi}{N}\right) \right| } .
  \end{equation*}
  This completes the proof of fact (iii).
\end{proof}
%% Mathematica code to verify the Taylor expansions:
%% f = -1/( Log[1-a(1-Cos[x])  ]);  Series[f,{x,0,1}]
%% f = -1/( Log[ Sqrt[ (1-a(1-Cos[x]) )^2  + b^2 Sin[x]^2]]); 
%% Series[f,{x,0,1}]

Next, we extend these results to another interesting set of matrices.  For
$N\geq 2$ and $a,b\in\real$, define the $N \times N$ matrices
$\atridmat_N^+(a,b)$ and $\atridmat_N^-(a,b)$ by
\begin{align*}
  \atridmat_N^{\pm}(a,b)&= \tridmat_N(a,b,a) \pm
  \begin{bmatrix}
    a & 0 &\dots &  \dots & 0 \\
    0 & 0 &\dots  &  \dots & 0 \\
    \vdots & \ddots & \ddots &\ddots & \vdots \\
    0 &  \dots &\dots &0 & 0 \\
    0 & \dots &\dots &0 & a
  \end{bmatrix}.
\end{align*}
If we define
\begin{align*}
  P_+ =
  \begin{bmatrix}
    1 & 1 & 0 & 0 &   \dots & 0\\
    1 & -1 & 1 & 0 &   \dots & 0\\
    1 & 0 & -1 & 1 &  \dots & 0\\
    \vdots && \ddots & \ddots & \ddots\\
    1 & 0 & \dots & 0 & -1 & 1\\
    1 & 0 & \dots & 0 & 0& -1
  \end{bmatrix}
  ,\enspace  \text{and} \quad
  P_- =
  \begin{bmatrix}
    1 & 1 & 0 & 0 &   \dots & 0\\
    -1 & 1 & 1 & 0 &   \dots & 0\\
    1 & 0 & 1 & 1 &  \dots & 0\\
    \vdots && \ddots & \ddots & \ddots\\
    (-1)^{N-2} & 0 & \dots & 0 & 1 & 1\\
    (-1)^{N-1} & 0 & \dots & 0 & 0& 1
  \end{bmatrix},
\end{align*}
then the following similarity transforms are satisfied:
\begin{equation}
  \label{eq:Pplus}
  \begin{split}
    \atridmat_N^+(a,b) &=
    P_+ \begin{bmatrix}
      b+2a & 0 \\
      0 & \tridmat_{N-1}(a,b,a)
    \end{bmatrix} P_+^{-1},
    \\
    \atridmat_N^-(a,b)      &=
    P_- \begin{bmatrix}
      b-2a & 0 \\
      0 & \tridmat_{N-1}(a,b,a)
    \end{bmatrix} P_-^{-1}.
  \end{split}
\end{equation}
To analyze the convergence properties of the dynamical systems determined
by $\atridmat_N^+(a,b)$ and $\atridmat_N^-(a,b)$, we recall that
$\mathbf{1}=(1,\dots,1)\in\real^d$, and we define $\mathbf{1}_- =
(1,-1,1,\dots,(-1)^{N-2},(-1)^{N-1}) \in\real^N$.
\begin{theorem}
  \label{thm:convergence-atrid}
  Let $N\geq 2$, $\eps\in]0,1[$, and $a,b\in\real$ with $a\neq 0$ and
  $|b|+2|a|=1$.  Let $\map{x}{\naturalzero}{\real^N}$ and
  $\map{z}{\naturalzero}{\real^N}$ be solutions to
  \begin{equation*}
    x(\ell+1) = \atridmat_N^+(a,b) \, x(\ell),
    \quad \text{and}\quad
    z(\ell+1) = \atridmat_N^-(a,b) \, z(\ell),
  \end{equation*}
  with initial conditions $x(0)=x_0$ and $z(0)=z_0$, respectively. The
  following statements hold: 
  \begin{enumerate}
  \item $\lim_{\ell\to+\infty} \big(x(\ell)-\xave(\ell)\mathbf{1}\big)
    =\mathbf{0}$, where
    $\xave(\ell)=(\frac{1}{N}\mathbf{1}^Tx_0)(b+2a)^\ell$, and the maximum
    time required for $\Enorm{x(\ell)-\xave(\ell)\mathbf{1}}\leq
    \eps\Enorm{x_0-\xave(0)\mathbf{1}}$ (over all initial conditions
    $x_0\in\real^N$) is $\Theta\big(N^2\log\eps^{-1}\big)$;
    
  \item $\lim_{\ell\to+\infty} \big(z(\ell)-\zave(\ell)\mathbf{1}_-\big)
    =\mathbf{0}$, where
    $\zave(\ell)=(\frac{1}{N}\mathbf{1}_-^Tz_0)(b-2a)^\ell$, and the
    maximum time required for $\Enorm{z(\ell)-\zave(\ell)\mathbf{1}_-}\leq
    \eps\Enorm{z_0-\zave(0)\mathbf{1}_-}$ (over all initial conditions
    $z_0\in\real^N$) is $\Theta\big(N^2\log\eps^{-1}\big)$.
  \end{enumerate}
\end{theorem} 
\begin{proof}
  We prove fact~(i) and remark that the proof of fact~(ii) is analogous.
  Consider the change of coordinates
  \begin{equation*}
    x(\ell)=P_+ \begin{bmatrix} \xave'(\ell) \\
      y(\ell) \end{bmatrix} = \xave'(\ell)\mathbf{1} + P_+
    \begin{bmatrix}
      0 \\ y(\ell)
    \end{bmatrix} ,
  \end{equation*}
  where $\xave'(\ell)\in\real$ and $y(\ell)\in\real^{N-1}$.  A quick
  calculation shows that $\xave'(\ell)=\frac{1}{N}\mathbf{1}^Tx(\ell)$, and
  the similarity transformation described in equation~\eqref{eq:Pplus}
  implies
  \begin{align*}
    y(\ell+1) &= \tridmat_{N-1}(a,b,a)\, y(\ell), 
    \enspace \text{and} \quad
    \xave'(\ell+1) = (b+2a)\xave'(\ell).
  \end{align*}
  Therefore, $\xave=\xave'$.  It is also clear that
  \begin{equation*}
    x(\ell+1)-\xave(\ell+1)\mathbf{1} = P_+
    \begin{bmatrix}
      0 \\ y(\ell+1)
    \end{bmatrix}   
    = \Big( P_+ \begin{bmatrix}
      0 & 0 \\
      0 & \tridmat_{N-1}(a,b,a)
    \end{bmatrix} P_+^{-1} \Big) (x(\ell)-\xave(\ell)\mathbf{1}).
  \end{equation*}
  Consider the matrix in parenthesis determining the trajectory
  $\ell\mapsto(x(\ell)-\xave(\ell)\mathbf{1})$. This matrix is symmetric,
  its singular values are $0$ and the singular values of
  $\tridmat_{N-1}(a,b,a)$, and its eigenvectors are $\mathbf{1}$ and the
  eigenvectors of $\tridmat_{N-1}(a,b,a)$ (padded with an extra zero).
  These facts are sufficient to duplicate, step by step, the proof of
  fact~(i) in Theorem~\ref{thm:convergence-trid-circ}.  Therefore, the
  trajectory $\ell\mapsto(x(\ell)-\xave(\ell)\mathbf{1})$ satisfies the
  same properties as those stated in
  Theorem~\ref{thm:convergence-trid-circ}(i).
\end{proof}
  
We conclude this section with some useful bounds.
\begin{lemma}\label{lem:bounds-P}
  Assume $x\in\real^N$, $y\in\real^{N-1}$ and $z\in\real^{N-1}$ jointly
  satisfy
  \begin{equation*}
    x = P_+ \begin{bmatrix} 0 \\ y \end{bmatrix},
    \qquad
    x = P_- \begin{bmatrix} 0 \\ z \end{bmatrix}.
  \end{equation*}
  Then $\frac{1}{2} \norm{x}{2} \leq \norm{y}{2} \leq (N-1) \norm{x}{2}$
  and $\frac{1}{2} \norm{x}{2} \leq \norm{z}{2} \leq (N-1) \norm{x}{2}$.\oprocend
\end{lemma}
\begin{proof}
  The proof is based on the coordinate expressions:
  \begin{gather*}
    x_1=y_1, \; x_2=y_2-y_1, \; \dots \; x_{N-1}=y_{N-1}-y_{N-2}, \; x_N=-y_{N-1}, \\
    y_1= x_1,\;  y_2= x_2+x_1,\; y_3=x_3+x_2+x_1,\; \dots \; y_{N-1}= x_{N-1}+\dots+x_1,
  \end{gather*}
  and
  \begin{gather*}
    x_1=z_1, \; x_2=z_2+z_1, \; \dots \; x_{N-1}=z_{N-1}+z_{N-2}, \; x_N=z_{N-1}, \\
    z_1=x_1,\;  z_2= x_2-x_1,\; z_3=x_3-x_2+x_1,\; \dots \; 
    z_{N-1}= x_{N-1}+\dots+(-1)^{N-1}x_1.
  \end{gather*}
\end{proof}

\end{document}